\documentclass[reqno, 11pt]{amsart}

\usepackage[paper=letterpaper,margin=.9in]{geometry}

\usepackage[title]{appendix}
\usepackage[utf8]{inputenc}
\usepackage{amsmath,amsthm,amssymb,amsfonts}
\usepackage[pdfborder={0 0 0}]{hyperref}
\hypersetup{colorlinks,linkcolor={blue},citecolor={blue},urlcolor={red}} 
\usepackage[pdftex,dvipsnames]{xcolor} 
\usepackage{tikz}
\usepackage{adjustbox}
\usetikzlibrary{arrows}
\usepackage{mathtools}
\usepackage{theoremref}
\usepackage{bm}
\usepackage{stmaryrd}
\SetSymbolFont{stmry}{bold}{U}{stmry}{m}{n}

\usepackage{bbm}
\usepackage{stmaryrd}
\usepackage{array}
\usepackage{enumitem}
\usepackage{stmaryrd}
\usepackage{comment}
\usepackage{graphicx}

\newtheorem{theorem}{Theorem}[section]
\newtheorem{lemma}[theorem]{Lemma}
\newtheorem{proposition}[theorem]{Proposition}
\newtheorem{corollary}[theorem]{Corollary}
\newtheorem{definition}[theorem]{Definition}

\newtheorem{remark}[theorem]{Remark}

\theoremstyle{definition}

\numberwithin{equation}{section}

\newcommand{\be}{\begin{equation}}
\newcommand{\ee}{\end{equation}}
\newcommand{\de}{\begin{align*}}
\newcommand{\fe}{\end{align*}}
\newcommand{\R}{\mathbb{R}}

\newcommand{\E}{\mathbb{E}}

\newcommand{\Z}{\mathbb{Z}}
\newcommand{\Pb}{\mathbb{P}}
\newcommand{\ind}{\mathbbm{1}}
\newcommand{\conf}{\text{Conf}(\mathfrak{X})}
\newcommand{\XX}{\mathfrak{X}}

\newcommand{\A}{\bm{\mathcal{A}}}
\newcommand{\B}{\bm{\mathcal{B}}}

\newcommand{\red}[1]{{#1}}

\definecolor{mintbg}{rgb}{.63,.79,.95}
\colorlet{lightmintbg}{mintbg!40}
\newcommand{\blue}[1]{\textcolor{blue}{#1}}

\title{The stochastic six-vertex model speed process}

\author[H.\ Drillick]{Hindy Drillick}
\address{H.\ Drillick,
	Department of Mathematics, Columbia University,
	\newline\hphantom{\quad \ \ H. Drillick}
	2990 Broadway, New York, NY 10027 USA
}
\email{hindy.drillick@columbia.edu}

\author[L.\ Haunschmid-Sibitz]{Levi Haunschmid-Sibitz}
\address{L.\ Haunschmid-Sibitz,
	KTH Royal Institute of Technology,
	\newline\hphantom{\quad \ L. Haunschmid-Sibitz}
	Lindstedtsvägen 25, SE - 100 44 Stockholm, Sweden
}
\email{levihs@kth.se}

\begin{document}

\begin{abstract}
For the stochastic six-vertex model on the quadrant $\Z_{\geq0}\times\Z_{\geq0}$ with step initial conditions and a single second-class particle at the origin, we show almost sure convergence of the speed of the second-class particle to a random limit. 
This allows us to define the stochastic six-vertex model speed process, whose law we show to be ergodic and stationary for the dynamics of the multi-class stochastic six-vertex process.
To prove this, we develop precise bounds on the fluctuations of the height function of the stochastic six-vertex model around its limit shape using methods from integrable probability. As part of the proof, we also obtain a novel geometric stochastic domination result that states that a second-class particle to the right of any number of third-class particles will at any fixed time be overtaken by at most a geometric number of third-class particles, which is of independent interest. 
\end{abstract}
\maketitle
\setcounter{tocdepth}{1}
\tableofcontents
\section{Introduction}
\subsection{Preface}

\begin{figure}[ht]
		\centering
		\begin{tabular}{|c|c|c|c|c|c|c|}
			\hline
			Type & I & II & III & IV & V & VI \\
			\hline
			\begin{tikzpicture}[scale = 1.5]
			%----
			\draw[fill][white] (0.5, 0) circle (0.05);
			\draw[thick][white] (0, 0) -- (1,0);
			\draw[thick][white] (0.5, -0.5) -- 
			(0.5,0.5);
			\node at (0.5, 0) {Configuration};
			\end{tikzpicture}
			&
			\begin{tikzpicture}[scale = 1.2]
			%----
			\draw[fill] (0.5, 0) circle (0.05);
			\draw[thick] (0, 0) -- (1,0);
			\draw[thick] (0.5, -0.5) -- (0.5,0.5);
			\end{tikzpicture}
			&
			\begin{tikzpicture}[scale = 1.2]
			%----
			\draw[thick][white] (0, 0) -- (1,0);
			\draw[thick][white] (0.5, -0.5) -- (0.5,0.5);
			\draw[fill] (0.5, 0) circle (0.05);
			\end{tikzpicture}
			&
			\begin{tikzpicture}[scale = 1.2]
			%----
			\draw[thick][white] (0, 0) -- (1,0);
			\draw[thick] (0.5, -0.5) -- (0.5,0.5);
			\draw[fill] (0.5, 0) circle (0.05);
			\end{tikzpicture}
			&
			\begin{tikzpicture}[scale = 1.2]
			%----
			\draw[thick][white] (0, 0) -- (0.5,0);
			\draw[thick][white] (0.5, 0) -- (0.5, 0.5);
			\draw[thick] (0.5, 0) -- (1, 0);
			\draw[thick] (0.5, -0.5) -- (0.5, 0);
			(0.5,0.5);
			\draw[fill] (0.5, 0) circle (0.05);
			\end{tikzpicture}
			&
			\begin{tikzpicture}[scale = 1.2]
			%----
			\draw[thick] (0, 0) -- (1,0);
			\draw[thick][white] (0.5, -0.5) -- (0.5,0.5);
			\draw[fill] (0.5, 0) circle (0.05);
			\end{tikzpicture}
			&
			\begin{tikzpicture}[scale = 1.2]
			%----
			\draw[thick] (0, 0) -- (0.5,0);
			\draw[thick] (0.5, 0) -- (0.5, 0.5);
			\draw[thick][white] (0.5, 0) -- (1, 0);
			\draw[thick][white] (0.5, -0.5) -- (0.5, 0);
			(0.5,0.5);
			\draw[fill] (0.5, 0) circle (0.05);
			\end{tikzpicture}
			\\
			\hline
			Weight 
			%(first parametrization)
			& 1 & 1 & $b_1$ & $1- b_1$ & $b_2$ & $1-b_2$\\
			\hline
		\end{tabular}
		\caption{The six allowed configurations for the stochastic six-vertex model}
		\label{fig:s6v}
	\end{figure}

  \begin{figure}[ht]
		\centering
			\begin{tikzpicture}[scale = 0.8]
		\foreach \x in {1, 2, 3, 4, 5, 6, 7}
		{\draw[lightmintbg] (\x, 0) -- (\x, 7.5);
			\draw[lightmintbg] (0, \x) -- (7.5, \x);
		}
		% \draw[lightmintbg, ->] (7, 0) -- (7.5, 0);
		% \draw[lightmintbg, ->] (0, 7) -- (0, 7.5);
		\foreach \x in {1, 2, 3, 4, 5, 6, 7}
		{
			\foreach \y in {1, 2, 3, 4, 5, 6, 7}
			\draw[fill] (\x, \y) circle (0.05);
		}
		\draw (0, 1) -- (4, 1);
		\draw (0, 2) -- (1, 2);
		\draw (0, 3) -- (2, 3);
		\draw (4, 3) -- (7.5, 3);
		\draw (0, 4) -- (4, 4);
		\draw (0, 5) -- (3, 5) -- (3, 7.5);
		\draw (0, 6) -- (7.5, 6);
		\draw (0, 7) -- (5, 7) -- (5, 7.5);
		\draw (1, 2) -- (1, 7.5);
		\draw (2, 3) -- (2, 7.5);
		\draw (3, 5) -- (3, 5.5);
		\draw (4, 1) -- (4, 3);
		\draw (4, 4) -- (4, 7.5);
		%\node at (0.6, 0.6) {$(0, 0)$};
		\foreach \x in {0,1,2,3,4,5,6} 
		\node at (\x+ 0.75, 0.75) {\footnotesize $\blue{0}$}; 
		\foreach \x in {0,1,2,3} 
		\node at (\x+ 0.75, 1.75) {\footnotesize $\blue{1}$}; 
		\foreach \x in {4,5,6} 
		\node at (\x+ 0.75, 1.75) {\footnotesize $\blue{0}$}; 
		\foreach \x in {0} 
		\node at (\x+ 0.75, 2.75) {\footnotesize $\blue{2}$};
		\foreach \x in {1,2,3} 
		\node at (\x+ 0.75, 2.75) {\footnotesize $\blue{1}$};
		\foreach \x in {4,5,6} 
		\node at (\x+ 0.75, 2.75) {\footnotesize $\blue{0}$};
		\foreach \x in {0} 
		\node at (\x+ 0.75, 3.75) {\footnotesize $\blue{3}$};
		\foreach \x in {1} 
		\node at (\x+ 0.75, 3.75) {\footnotesize $\blue{2}$};
		\foreach \x in {2,3,4,5,6} 
		\node at (\x+ 0.75, 3.75) {\footnotesize $\blue{1}$};
		\foreach \x in {0} 
		\node at (\x+ 0.75, 4.75) {\footnotesize $\blue{4}$};
		\foreach \x in {1} 
		\node at (\x+ 0.75, 4.75) {\footnotesize $\blue{3}$};
		\foreach \x in {2, 3} 
		\node at (\x+ 0.75, 4.75) {\footnotesize $\blue{2}$};
		\foreach \x in {4, 5, 6} 
		\node at (\x+ 0.75, 4.75) {\footnotesize $\blue{1}$};
		\foreach \x in {0} 
		\node at (\x+ 0.75, 5.75) {\footnotesize $\blue{5}$};
		\foreach \x in {1} 
		\node at (\x+ 0.75, 5.75) {\footnotesize $\blue{4}$};
		\foreach \x in {2} 
		\node at (\x+ 0.75, 5.75) {\footnotesize $\blue{3}$};
		\foreach \x in {3} 
		\node at (\x+ 0.75, 5.75) {\footnotesize $\blue{2}$};
		\foreach \x in {4,5,6} 
		\node at (\x+ 0.75, 5.75) {\footnotesize $\blue{1}$};
		\foreach \x in {0} 
		\node at (\x+ 0.75, 6.75) {\footnotesize $\blue{6}$};
		\foreach \x in {1} 
		\node at (\x+ 0.75, 6.75) {\footnotesize $\blue{5}$};
		\foreach \x in {2} 
		\node at (\x+ 0.75, 6.75) {\footnotesize $\blue{4}$};
		\foreach \x in {3} 
		\node at (\x+ 0.75, 6.75) {\footnotesize $\blue{3}$};
		\foreach \x in {4,5,6} 
		\node at (\x+ 0.75, 6.75) {\footnotesize $\blue{2}$};
%		\node at (-0.3, -0.3) {\footnotesize $(0, 0)$};
\end{tikzpicture}
		\caption{A possible sampling of the stochastic six-vertex model on the quadrant with step initial data. The height function is denoted in blue.}
		\label{fig:s6vSampling}
\end{figure}

The stochastic six-vertex model was first introduced by Gwa and Spohn in \cite{GwaSpohnStochasticSixVertex} as a specialization of the six-vertex model, which is a classical model from equilibrium statistical mechanics going back to \cite{pauling1935structure}.
Recently there has been a lot of interest in this model.
It is connected via a suitable limit degeneration to ASEP \cite{Aggarwal2017S6VtoASEP}, the Kardar-Parisi-Zhang equation \cite{CorwinTsai2017KPZlimit,CorwinGhosalShenTsai2020SPDElimit, Yier2020KPZlimit}, the stochastic telegraph equation \cite{BorodinGorin2019StochasticTelegraph, MR3951443, MR4193889} and lies in the (one-dimensional) KPZ universality class \cite{GwaSpohnStochasticSixVertex, BCGStochasticSixVertex, MR4561796, aggarwal2024scaling}.
Furthermore, it can be put into the more general setting of higher-spin vertex models, see \cite{CorwinPetrov2016higherspinLine, BorodinPetrov2018higherspinSymmetric}.

We will now define the stochastic six-vertex model on the quadrant $\Z_{\geq 0}\times\Z_{\geq 0}$ as a stochastic path ensemble, where each path can be thought of as a particle trajectory. The model depends on two parameters $b_1, b_2 \in [0,1]$ that we fix throughout the paper. We call the bottom and left edges incident to a vertex its \textbf{incoming} edges and the top and right edges its \textbf{outgoing} edges. We say that edges are \textbf{occupied} if they are part of a particle trajectory. Each edge will be occupied by at most one particle, and the model will satisfy a local conservation law, namely that the number of particles entering a vertex equals the number of particles exiting.

We first specify a boundary configuration on the incoming edges along the left and bottom boundaries of the quadrant, i.e., the incoming edges from the left at the vertices $\{0\}\times\Z_{\geq0}$ and from the bottom at the vertices $\Z_{\geq0}\times\{0\}$. This boundary configuration specifies a subset of edges along the boundary that are occupied, i.e., locations where particles enter the system. The most common boundary condition that we will work with is one where all incoming edges from the left boundary of the quadrant are occupied and all incoming edges from the bottom boundary are empty.
We will refer to this boundary condition as \textbf{step initial conditions}.
% in analogy with analogous initial conditions in the theory of interacting particle systems. 

Given the information about which incoming edges of a vertex are occupied, we will sample the outgoing edges using the weights in Figure \ref{fig:s6v}. Note that the weights of configurations with the same occupied incoming edges add up to $1$, so that the outgoing edges can be determined in a stochastic manner.  We will begin from $(0,0)$ and proceeding  in an antidiagonal way (i.e, ordered by $x+y$), we will tile each of the vertices in the interior of the quadrant with one of the six possible configurations in Figure \ref{fig:s6v}, according to their weights. This procedure defines a law on configurations of upright paths in $\Z_{\geq0}\times \Z_{\geq0}$. See Figure \ref{fig:s6vSampling} for a sampling of the stochastic six-vertex model with step initial conditions.

There is an alternative parameterization of the model by parameters $q,\kappa>0$ defined as
\[
    q\coloneqq\frac{b_1}{b_2}\qquad\kappa\coloneqq\frac{1-b_1}{1-b_2}\,.
\]
This parameterization will be quite useful to us, and these variables will appear in many formulas throughout the paper.

For a given configuration $\omega$ of the stochastic six-vertex model with step initial conditions, we define the height function $H(x,t) = H(x,t;\omega)$ for $x,t\in\R_{\geq0}$ by setting $H(x,0;\omega)=0$ for all $x$ and increasing $H$ whenever one crosses a path in the vertical direction, see Figure \ref{fig:s6vSampling}.

\begin{figure}\label{fig:initial}
    \centering
    \includegraphics[width=0.3\textwidth]{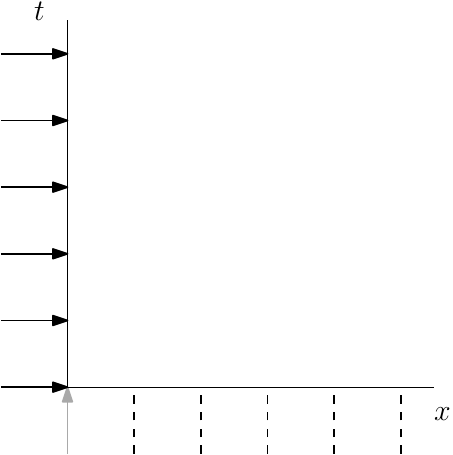}\hspace{0.1\textwidth}\includegraphics[width=0.3\textwidth]{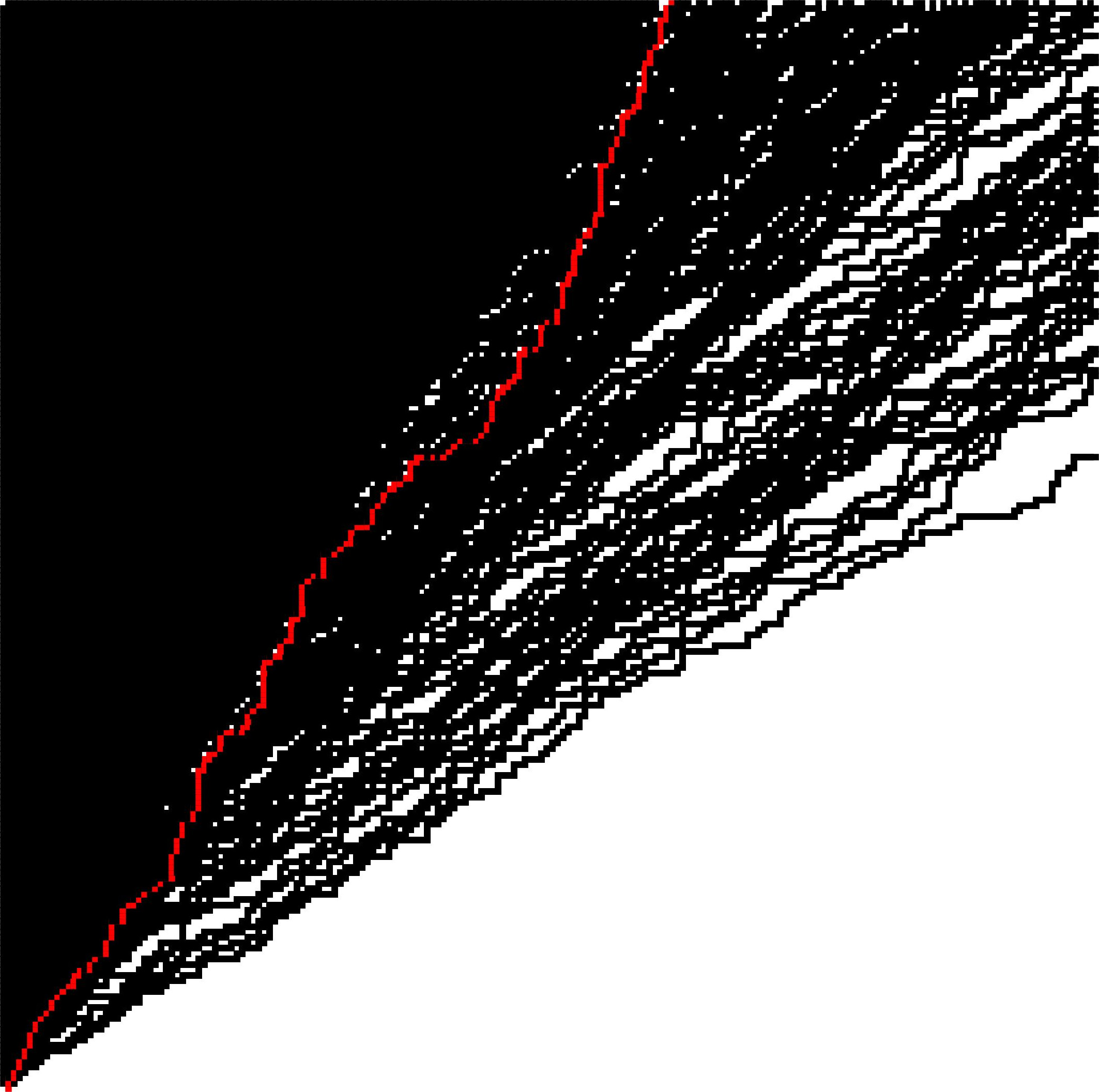}
    \caption{Left panel: step initial conditions with a second-class particle at the origin. Black arrows denote first-class particles, while the grey arrow denotes the second-class particle. Dashed lines denote holes.
    Right panel: a simulation of this process on a 200 by 200 square with $b_1=0.3$ and $b_2=0.6$ and with the second-class particle in red.}
    \label{fig:initial_conditions}
\end{figure}

The model exhibits two very different behaviors depending on whether $b_1$ or $b_2$ is larger.
If $b_1>b_2$, then particles prefer moving up to moving to the right.
Since the upper part of the quadrant is already packed, this leads to a sharp transition between a region with density $1$ and a region with density $0$, whose boundary stays close to the line $x=t$. This behavior is known as a \textbf{shock}.
On the other hand if $b_1<b_2$, particles want to move right more than up, and thus they spread out.
Three regions form: one above the line $x=\kappa^{-1}t$, where the density of particles is $1$, one below the line $x=\kappa t$, where the density is $0$, and one in between, where the density decreases continuously from $1$ to $0$ (see the right-hand side of Figure~\ref{fig:initial_conditions} for a simulation). The middle section is known as the \textbf{rarefaction fan}.
Both the shock and rarefaction fan regimes are interesting in their own rights, but our results concern the latter: from now on we always assume $b_1<b_2$.

We will now introduce the multi-class stochastic six-vertex model.
Instead of every edge being either occupied or unoccupied it will now be assigned a \textbf{class} from $\mathbb R\cup\{-\infty,\infty\}$.
The classes assigned to the two outgoing edges equal the classes of the incoming edges, and the weight of a vertex depends on the classes, see Figure \ref{fig:multiClasss6v}.
Intuitively if $i<j$ then a particle of class $i$ treats particles of class $j$ as \textbf{holes}.
The single-class stochastic six-vertex model can be obtained from the multi-class one by setting the class of unoccupied edges to $1$ and the class of occupied edges to $\infty$. In the vertex model literature the multi-class stochastic six-vertex model is also called the \emph{colored stochastic six-vertex model}. See \cite{Kulish:1981aa, MR806529, MR824090, kuniba2016stochastic, MR4518478} for previous works on this model from an algebraic perspective. 

\begin{figure}[t]
		\centering
		\begin{tabular}{|c|c|c|c|c|c|}
			\hline
			\begin{tikzpicture}[scale = 1.5]
			%----
			\draw[fill][white] (0.5, 0) circle (0.05);
			\draw[thick][white] (0, 0) -- (1,0);
			\draw[thick][white] (0.5, -0.5) -- 
			(0.5,0.5);
			\node at (0.5, 0) {Configuration};
			\end{tikzpicture}
			&
			\begin{tikzpicture}[scale = 1.2]
			%----
			\draw[thick, red] (0, 0) -- (1,0);
			\draw[thick, red] (0.5, -0.5) -- (0.5,0.5);
            \draw[fill] (0.5, 0) circle (0.05);
            \node at (-0.2, 0) {$i$};
            \node at (0.5, 0.7) {$i$};
            \node at (1.2, 0) {$i$};
            \node at (0.5, -0.7) {$i$};

			\end{tikzpicture}
			&
			% \begin{tikzpicture}[scale = 1.2]
			% %----
			% \draw[thick][blue] (0, 0) -- (1,0);
			% \draw[thick][blue] (0.5, -0.5) -- (0.5,0.5);
			% \draw[fill] (0.5, 0) circle (0.05);
			% \end{tikzpicture}
			% &
			\begin{tikzpicture}[scale = 1.2]
			%----
			\draw[thick, blue](0, 0) -- (1,0);
			\draw[thick, red] (0.5, -0.5) -- (0.5,0.5);
			\draw[fill] (0.5, 0) circle (0.05);
            \node at (-0.2, 0) {$j$};
            \node at (0.5, 0.7) {$i$};
            \node at (1.2, 0) {$j$};
            \node at (0.5, -0.7) {$i$};
			\end{tikzpicture}
			&
			\begin{tikzpicture}[scale = 1.2]
			%----
			\draw[thick][blue] (0, 0) -- (0.5,0);
			\draw[thick][blue] (0.5, 0) -- (0.5, 0.5);
			\draw[thick, red] (0.5, 0) -- (1, 0);
			\draw[thick, red] (0.5, -0.5) -- (0.5, 0);
			(0.5,0.5);
			\draw[fill] (0.5, 0) circle (0.05);
            \node at (-0.2, 0) {$j$};
            \node at (0.5, 0.7) {$j$};
            \node at (1.2, 0) {$i$};
            \node at (0.5, -0.7) {$i$};
			\end{tikzpicture}
			&
			\begin{tikzpicture}[scale = 1.2]
			%----
			\draw[thick, red] (0, 0) -- (1,0);
			\draw[thick][blue] (0.5, -0.5) -- (0.5,0.5);
			\draw[fill] (0.5, 0) circle (0.05);
            \node at (-0.2, 0) {$i$};
            \node at (0.5, 0.7) {$j$};
            \node at (1.2, 0) {$i$};
            \node at (0.5, -0.7) {$j$};
			\end{tikzpicture}
			&
			\begin{tikzpicture}[scale = 1.2]
			%----
			\draw[thick, red] (0, 0) -- (0.5,0);
			\draw[thick, red] (0.5, 0) -- (0.5, 0.5);
			\draw[thick, blue] (0.5, 0) -- (1, 0);
			\draw[thick, blue] (0.5, -0.5) -- (0.5, 0);
			(0.5,0.5);
			\draw[fill] (0.5, 0) circle (0.05);
            \node at (-0.2, 0) {$i$};
            \node at (0.5, 0.7) {$i$};
            \node at (1.2, 0) {$j$};
            \node at (0.5, -0.7) {$j$};
			\end{tikzpicture}
			\\
			\hline
			Weight 
			%(first parametrization)
			& 1  & $b_1$ & $1- b_1$ & $b_2$ & $1-b_2$\\
			\hline
		\end{tabular}
		\caption{The allowed configurations for the multi-class stochastic six-vertex model, where red lines represent class $i$ and blue lines represent class $j$ for $i < j$. }
		\label{fig:multiClasss6v}
	\end{figure}

Our main theorem concerns the following variant of step initial conditions, which we will call \textbf{step initial conditions with a second-class particle at the origin}.
All particles coming in from the left have class $1$, there is a single particle coming in from the bottom at $(0,0)$, and all other incoming particles from the bottom have class $\infty$, i.e. are holes, see Figure~\ref{fig:initial_conditions}.
By the conservative property of the model, for every $t$ there is exactly one $x$ such that the vertical arrow leaving $(x,t)$ has class $2$.
We call this $x$ the \textbf{position of the second-class particle at time $t$} and denote it by $\bm{X}_t$. Our main result states that the speed $\frac{\bm{X}_t}{t}$ of the second-class particle converges a.s. to a random limit:

\begin{theorem}\label{thm:main} Let $0<b_1<b_2 < 1$, and consider the stochastic six-vertex model with step initial positions with a second-class particle at the origin.
Let $\bm{X}_t$ be the position of the second-class particle at time $t$.
Then almost surely
\begin{equation}
    \lim_{t \to \infty} \frac{\bm{X}_t}{t} = \bm{U}
\end{equation}
where $\bm{U}$ is a continuous random variable taking values in $[\kappa^{-1}, \kappa]$ with density $\frac{\sqrt{\kappa}}{2(\kappa-1)}x^{-\frac32}$.
\end{theorem}
Even the weak convergence of the speed of the second-class particle has not been stated in the literature, to the best of the authors' knowledge.
However, it follows readily from the hydrodynamic limit proved in \cite{Aggarwal2020} in the same way as for ASEP using the arguments from \cite{FerrariKipnis1995}.
For the convenience of the reader, we adapt this argument to our setting in Appendix~\ref{ap:weak}. Strengthening the weak convergence to a.s. convergence is highly nontrivial and is the main contribution of this theorem.

For TASEP with step initial conditions, the weak convergence of the speed of a second-class particle at the origin was first proven in \cite{FerrariKipnis1995} and a.s. convergence was proven in \cite{MountfordGuiol05} (see also \cite{FerrariPimentel05} and \cite{FerrariMartinPimentel09} for alternative proofs). Almost sure convergence for second-class particles starting with suitable initial conditions was then proven for the Hammersley process \cite{MR2312944}, the totally asymmetric zero range process (TAZRP) \cite{Goncalves14}, and the Plancherel-TASEP particle system \cite{MR3306003}.
% For the Hammersley process, a.s. convergence of the speed of a second-class particle at the origin was proven in \cite{MR2312944} under suitable initial conditions, and for the totally asymmetric zero range process (TAZRP), this was proven in \cite{Goncalves14}. \red{The work \cite{MR3306003} proved the same result for the second-class particle in the Plancherel-TASEP particle system}. 
All of these proofs rely crucially on connections between the models under consideration and last passage percolation models (for example, TASEP can be coupled with exponential last passage percolation). 
Since this no longer holds for ASEP, new tools were required to prove the analogous result for ASEP under step initial conditions, and this was done in \cite{ACG2023asepspeed} using inputs from integrable probability as well as a coupling due to Rezakhanlou \cite{REzakhanlou1995Shocks}. Since for the stochastic six-vertex model, last passage methods also do not apply, our proof strategy for Theorem \ref{thm:main} is inspired by the ideas in  \cite{ACG2023asepspeed}. The speed of second-class particles for ASEP and the Hammersley process has also been studied for other classes of initial conditions in \cite{ MR3112921, ghosal2019limitingspeedsecondclass, MR4010933}.

We also derive a bound on the fluctuations around the limiting speed:
\begin{theorem}\label{thm:fluct}
    Let $\bm{X}_t$ be the position of the second-class particle at time $t$ as above and $\bm{U}$ its almost sure limit.
    Then for any $\delta>0$, almost surely we have that 
    \begin{equation}\label{eq:fluct}
        \lim_{t\to\infty}|\bm{X}_t-t\bm{U}|t^{-(\frac79+\delta)}=0\,.
    \end{equation}
\end{theorem}

The exponent $-\left(\frac{7}{9} + \delta\right)$ is not optimal. For both ASEP and the stochastic six-vertex model with stationary initial conditions, a second-class particle starting at the origin has been shown to fluctuate on the order of $t^{2/3}$. This was shown for ASEP in \cite{MR2318311, MR2630064} and for the stochastic six-vertex model in \cite{landonSosoe2023tail} (building on earlier ideas developed in \cite{Aggarwal2016CurrentFO}). However, the exact fluctuations of a second-class particle starting from step initial conditions are not known for either model. Because the fluctuations at stationarity are of the order $t^{2/3}$, the best exponent one could expect to achieve in \eqref{eq:fluct} is $-(\frac23+\delta)$, see also Remark \ref{rmk:sn}. Our proof techniques can be used for ASEP as well, where they would give an analogous result to Theorem \ref{thm:fluct} with the same exponent of $-\left(\frac{7}{9} + \delta\right)$.

Going beyond adding a single second-class particle into our model, we can consider initial conditions where each incoming particle has a different class in $\Z$.
Individually, each particle will have an asymptotic speed given by Theorem~\ref{thm:main}. 
By considering the joint speeds of all the particles simultaneously, we can construct the \textbf{stochastic six-vertex model speed process}.
Speed processes have previously been constructed and studied for TASEP \cite{AmirAngelValko2008TasepSpeed}, TAZRP \cite{ABGMtazrpspeed}, and ASEP \cite{ACG2023asepspeed}. To define the speed process, we first need to discuss how we can view the stochastic six-vertex model as a particle system, as was first done in \cite{GwaSpohnStochasticSixVertex}, see also \cite[Section 2.2]{BCGStochasticSixVertex}. 

\subsection{The stochastic six-vertex model as an interacting particle system}\label{sec:particlesystem}

Until this point, we have treated the stochastic six-vertex model as a measure on up-right paths.
However, it is also natural to consider it as a particle system with time as the vertical direction, as has been quite noticeable in the language we have been using and was already observed in \cite{GwaSpohnStochasticSixVertex}.
Let us now introduce notation that emphasizes this connection.
For a given configuration $\omega$ on $\Z_{\geq0}\times\Z_{\geq0}$, define $\eta_t(x)$ for $x\in\Z_{\geq0}$ by
\[
    \eta_t(x)=\begin{cases}
        1 &\text{if the incoming vertex at $(x,t)$ from below in $\omega$ is occupied}\\
        0&\text{else.}
    \end{cases}
\]

Defined like this $(\eta_t)_{t\in\Z_{\geq0}}$ is a Markov process with values in $\{0,1\}^{\Z_{\geq0}}$.
We call this a \textbf{stochastic six-vertex process}. 
The boundary conditions on the bottom give the initial condition $\eta_0$ and the boundary conditions on the left inject particles at specific times.
The transition probabilities of this process can be described as follows: Particles stay in place with probability $b_1$ and start moving to the right with probability $(1-b_1)$. If a particle starts moving, the amount it moves is the minimum of a $\text{Geo}(b_2)$ distributed random variable and the distance to the nearest particle to its right. If it moves to the location of the neighboring particle to the right, that other particle then starts moving, following the above-described rules. 
See \cite[Section 2.2]{BCGStochasticSixVertex} for these transition weights written out in more detail. 

We now define the height function in this setting and show that it generalizes the definition of $H(x,t)$ above for the case of step initial conditions. 
\begin{definition}[Height Function]\label{def:heightfunction}
For a given stochastic six-vertex process $(\eta_t)_{t\geq0}$, the height function $h_t(x)=h_t(x;\eta)$ is the unique function (up to a global shift) that satisfies
\begin{align}
    h_t(x;\eta)-h_t(x+1;\eta)&=\eta_t(x)\label{eq:heightgradienthorizontal}\\
    h_{t+1}(0;\eta)-h_{t}(0;\eta)&=\label{eq:heightgradientvertical}\begin{cases}
        1&\text{if there is an incoming arrow from the left at }(0,t)\\
        0&\text{else.}
    \end{cases}
\end{align}
\end{definition}
Since the height function is only unique up to a global shift, unless otherwise specified the choice of height function is made by setting $h_0(0)=0$, but in some places, it will be convenient to choose some other $h_0(0)$.
For a configuration $\omega$ of the stochastic six-vertex model with step initial conditions, one recovers the definition of $H(x,t)$ above, since by \eqref{eq:heightgradienthorizontal}, $h_0(x)=0$ for all $x$. 

\begin{definition}\label{def:s6vLine}
    As shown in \cite{Aggarwal2020}, these dynamics can be extended to processes $\eta_t:\Z\to\{0,1\}$.
We call this the \textbf{stochastic six-vertex process on the line}. 
\end{definition}
 
Given an initial condition $\eta_0:\Z\to\{0,1\}$ that satisfies $\eta_0(x)=\bm{1}_{x < 0}$, the restriction $(\eta_t(x))_{x,t\in\Z_{\geq0}}$ of the stochastic six-vertex process on the line to $x\geq0$ agrees with the process on the quadrant with step initial conditions.
It is this process that we will be considering in Sections~\ref{sec:couplings} to \ref{sec:hydroToLinear}.
The height function is still defined by \eqref{eq:heightgradienthorizontal} and \eqref{eq:heightgradientvertical}.

This extension is also compatible with the multi-class stochastic six-vertex process. While the single-class processes $\eta_t:\Z\to\{0,1\}$ have occupation variables in $\{0,1\}$ with $0$ encoding holes and $1$ encoding particles, we will let the multi-class processes have occupation variables in $\R\cup\{\infty\}$, with $\infty$ encoding holes and all other values encoding particles of different classes. In other words, we define the \textbf{multi-class stochastic six-vertex process on the line} as $\eta_t:\Z\to\ \R\cup\{\infty\}$, where $\eta_t(x) = i$ if at time $t$, there is a particle of class $i$ at position $x$. To avoid confusion, we will always specify in the text whether we are considering a single- or multi-class process.

We can now define the speed process whose existence will be obtained as a corollary of Theorem \ref{thm:main}. 
\begin{corollary}[Existence of the Speed Process]\label{cor:speeddef}
Consider the multi-class stochastic six-vertex model on the line with initial conditions $\eta_0(x)=x$ for all $x\in\mathbb Z$, i.e. at position $x$ there is a particle of class $x$. We call this \textbf{packed initial conditions}. Denote by $\bm{X}_t(x)$ the position of the unique particle of class $x$ at time $t$.
Then the process $\left(\frac{\bm{X}_t(x)}{t}\right)_{x\in\mathbb{Z}}$ converges a.s. as $t\to\infty$ to a process $U(x)$.
We call $U(x)$ the \textbf{stochastic six-vertex model speed process}. 
\end{corollary}

Now that the stochastic six-vertex model speed process is defined, we can study some of its properties. 
Speed processes were first studied in the work of \cite{AmirAngelValko2008TasepSpeed} where they demonstrated that the TASEP speed process has a highly nontrivial structure. They proved that the TASEP speed process is stationary with respect to the dynamics of the multi-class TASEP, and they computed some of its marginals. In particular they showed that there is a \emph{nonzero} probability that the speeds of two neighboring particles are identical and for each particle, there are infinitely many particles with the same speed. This phenomenon where multiple particles have the same speed is known as a \emph{convoy}. They also showed that many of the same results should hold for the ASEP speed process conditioned on its construction (which was carried out in \cite{ACG2023asepspeed}). See also \cite{MartinStationaryASEP} for more results about the ASEP speed process.

In Section \ref{sec:speedsymmetries}, we will prove that the stochastic six-vertex model speed process is ergodic and stationary with respect to the dynamics of the multi-class stochastic six-vertex model.
The uniqueness of multi-class stationary measures with a given marginal for ASEP is known, see \cite{liggett1976coupling,ferrari1991microscopic} and was already believed to hold for the stochastic six-vertex model as well, see \cite[Remark 7.9]{ANPstationarityBaxter}.
Since we were not able to find a proof of this in the literature we state it in Proposition~\ref{prop:uniqueStationaryMeasures} and provide a proof in Appendix~\ref{ap:unique}.
It then follows from a close examination of the construction of such stationary measures in \cite{ANPstationarityBaxter} that the ergodic stationary measures for the multi-class stochastic six-vertex model on the line are the same as the ergodic stationary measures for the multi-class ASEP. This then implies that the stochastic six-vertex speed process is related to the ASEP speed process by a deterministic map given by a pointwise composition with a specific map, see Proposition \ref{prop:deterministicmapping}. In particular, this implies that the same convoy structure studied in \cite{AmirAngelValko2008TasepSpeed, MartinStationaryASEP} is also present for the stochastic six-vertex model speed process.

There are also many avenues for further work on these processes.
In particular, the article \cite{busani2022scaling} shows that the suitably rescaled TASEP speed process converges weakly to a process known as the stationary horizon. The stationary horizon was first introduced in \cite{busani2023diffusivescalinglimitbusemann} and is expected to be a universal scaling limit for multi-class invariant measures of models in the KPZ universality class. 
Then in \cite{busani2024scalinglimitmultitypeinvariant}, they develop a more general framework to show convergence to the stationary horizon. In particular, they show that if a model converges to the directed landscape under suitable rescaling, then the stationary measures of the associated multi-class process converge to the stationary horizon at the level of finite-dimensional projections. 
In \cite{aggarwal2024scaling} they prove the convergence of the stochastic six-vertex model and ASEP to the directed landscape, and hence using the results from \cite{busani2024scalinglimitmultitypeinvariant}, they obtain as a corollary \cite[Corollary 2.14]{aggarwal2024scaling} that the stationary measures for the multi-class ASEP converge to marginals of the stationary horizon. As mentioned above, these stationary measures are the same as for the multi-class stochastic six-vertex model. It is still an open problem to prove the convergence of the ASEP and stochastic six-vertex model speed processes to the stationary horizon in the space $D(\R, C(\R))$.

In this paper, we deal with step and packed initial conditions for the (multi-class) stochastic six-vertex model. These initial conditions play the role of ``fundamental solutions" from which more general initial conditions can then be studied. Step initial conditions have long been of interest in the study of interacting particle systems and have been particularly important in the study of KPZ universality for the stochastic six-vertex model, see, e.g., \cite{BCGStochasticSixVertex,CorwinDimitrov2018Transversalfluctuations, MR4561796, dasLiaoMucciconi2024lower}. Packed initial conditions (which can be viewed as a coupling of different step initial conditions) have recently been the subject of \cite{aggarwal2024scaling} where they showed that stochastic six-vertex model and ASEP with packed initial conditions converge to the Airy Sheet. In \cite{aggarwal2024kpzfixedpointconvergence}, this convergence was extended to more general classes of initial conditions, using the convergence for packed initial conditions as a key input.
% Furthermore, it is precisely using packed initial conditions that we can obtain the speed process which yields the stationary measures for the multiclass model. 

\subsection{Proof Ideas}
The proof of the main theorem uses a variety of tools. We follow the general strategy developed in \cite{ACG2023asepspeed}, which requires certain model-specific inputs that have not yet been developed for the stochastic six-vertex model. 
In particular, we prove the following two key ingredients:
\begin{itemize}
    \item A novel geometric stochastic domination result that states that a second-class particle to the right of any number of third-class particles will at any fixed time be overtaken by at most a geometric number of third-class particles.
    \item Effective hydrodynamic estimates that quantify how close the height function of the stochastic six-vertex model started from step initial conditions will be to its limit shape.
\end{itemize}

These results will be used in the following way. We want to control the behavior of a single second-class particle. Hydrodynamic theory allows us to control the bulk behavior of many particles, so we augment our system by filling up all empty positions to the left of $\bm{X}_t$ with third-class particles. We then use our effective hydrodynamic estimates to control the union of the second- and third-class particles. Finally, we can revert this back to an estimate of the position of the second-class particle since we know that our second-class particle is to the left of at most a geometric number of the third-class particles. A similar argument can be made to bound the position of the second-class particle from the left. 

We now state these two results in detail. The first will be the content of Theorem~\ref{thm:orderinginequality} and the second, the content of Propositions \ref{prop:MeixnerTail} and \ref{prop:otherTail}.

\subsubsection{Controlling a Second-Class Particle by Third-Class Particles} The following theorem allows us to control the behavior of a single second-class particle by controlling the behavior of a large number of third-class particles inserted to the left of the second-class particle.

Recall that $q = \frac{b_1}{b_2}$. By $X \sim \mathrm{Geo}(q)$ we denote the law given by
\[
\mathbb P[X=k]=(1-q)q^k\text{ for }k\geq0\,.
\]
\begin{theorem}[Geometric Stochastic Domination]\label{thm:orderinginequality}
    Let $(\eta_t)_{t\geq 0}$ be a multi-class stochastic six-vertex process on the line with parameters $0<b_1<b_2<1$ and with the following initial conditions:
    \begin{itemize}
        \item There are some first-class particles (finitely or infinitely many).
        \item There is a single second-class particle.
        \item There are $M$ third-class particles, all to the left of the second-class particle.
    \end{itemize}
    Let $\bm{Z}_t(0)>\bm{Z}_t(1)>\dots>\bm{Z}_t(M)$ be the ordered positions of the second- and third-class particles at time $t$, i.e. $\eta_t(x)=2$ or $\eta_t(x)=3$ if and only if there is a $k$ such that $\bm{Z}_t(k)=x$, and further, the $\bm{Z}_t(k)$ are decreasing in k.
    Further, let $L_t$ be the number of third-class particles to the right of the second-class particle at time $t$.
    Then for any $t$, the law of $L_t$, conditioned on both $(\bm{Z}_s(k))_{0\leq k\leq M,0\leq s\leq t}$ and the space-time history of the first-class particles up to time $t$, is dominated by $\mathrm{Geo}(q)$.
\end{theorem}

A similar result for ASEP appears in \cite[Lemma 4.1]{MR2630064}. 
% That proof proceeds by showing that the number of overtaking particles in that context is given by a variant of a birth death chain $M(t)$ which has the geometric distribution $\pi$ as its reversible measure. In other words, if we let $M_{\text{stat}}$ denote the reversible birth death chain starting with initial conditions $M_{\text{stat}}(0) \sim  \pi$, then $M_{\text{stat}}(t) \sim \pi$ for all $t > 0$. The result then follows from the monotonicity of the basic coupling for ASEP which ensures that if we start this birth death chain with initial conditions $M(0) = 0 \leq M_{\text{stat}}(0)$, then $M(t) \leq M_{\text{stat}}(t)$ for all $t > 0$ if the coupling between the two processes is the one induced by the basic coupling. 
% However, in our setting, we no longer have this monotonicity property. In particular, it is possible to construct an initial condition for $N(t)$ that is stochastically dominated by a geometric random variable, yet at later times this no longer holds. Therefore, we need to work with the smaller class of configurations that can be written as a convex combination of the basis vectors $\mathbf{v}_i$, and this class now has the property that it maintains this ordering. 
The difference between the proofs of \cite[Lemma 4.1]{MR2630064} and Theorem~\ref{thm:orderinginequality} can be best understood by considering initial conditions in which the second-class particle is not necessarily the rightmost particle among the second- and third-class particles, but rather has a random position among the $(Z_0(k))_{0\leq k\leq M}$. In other words, $L_0$ is a random variable with support $\{0,\dots,M\}$. For simplicity we take $M = \infty$. 
One can then see that in both settings the stationary measure for $L_t$ is given by a $\text{Geo}(q)$ random variable.
For ASEP, certain monotonicity properties imply that for any initial condition on $L_0$ which is dominated by $\text{Geo}(q)$, $L_t$ will be dominated by $\text{Geo}(q)$ for any future time.
This is not true in our setting as these monotonicity properties no longer hold for the stochastic six-vertex model (see the discussion in Section \ref{sec:couplings}).
To circumvent this, in the proof of Theorem~\ref{thm:orderinginequality} we consider a smaller class of distributions that are preserved by the dynamics, which contains the deterministic initial condition $L_0=0$ and whose elements are all dominated by $\text{Geo}(q)$.
This class is given by distributions which can be obtained by taking the minimum of a $\text{Geo}(q)$ random variable and an independent random variable supported on non-negative integers.

    Let us also briefly compare this result with Rezakhanlou's coupling from \cite{REzakhanlou1995Shocks}, which was used to control a second-class particle in ASEP in \cite{ACG2023asepspeed}.
    In \cite{REzakhanlou1995Shocks} an auxiliary label process on the second and third-class particles is defined, which has the following properties.
    \begin{itemize}
        \item Every second- and third-class particle has a unique label from $0$ to $M$, which can change over time.
        \item The law of this labeling process at any fixed time is that of a uniform permutation, and it is stationary.
        \item It is coupled to the dynamics of the multi-class ASEP, such that at any time, the particle with label $1$ is to the left of the single second-class particle.
    \end{itemize}
    This allows us to control the second-class particle with a uniformly chosen third-class particle, see \cite[Proposition 5.4]{ACG2023asepspeed}.

    One can construct an analogous coupling for the stochastic six-vertex model,\footnote{Such a coupling was presented by Ivan Corwin at the 2022 PIMS-CRM Summer School in Probability.} but only for the case $b_1<\frac12$. This is essentially due to the fact that certain monotonicity properties that hold for ASEP, do not hold for the stochastic six-vertex model (see the discussion before Proposition \ref{prop:approxmono}). Theorem~\ref{thm:orderinginequality} takes a different approach and works for all $b_1<b_2$. There are two key differences between these approaches:
    \red{Firstly, Theorem~\ref{thm:orderinginequality} does not proceed by explicitly coupling $(L_t)_{t\geq0}$ and a process with Geo$(q)$ marginals.}
    Secondly, the bound in Theorem~\ref{thm:orderinginequality} is significantly stronger for large $M$.
    Intuitively, the result from \cite{REzakhanlou1995Shocks} shows that the number of third-class particles that \emph{do not} pass the second-class particle grows linearly in the number of third-class particles, while Theorem~\ref{thm:orderinginequality} shows that the number that \emph{do} pass is of order $1$.

    Since the statement of Theorem~\ref{thm:orderinginequality} is entirely insensitive to scaling time or space, it can be carried over to ASEP, with $q=\frac{b_1}{b_2}$ fixed.
    For ASEP this result could also be obtained from the censoring inequality \cite{PeresWinkler2013CensoringInequality} (as pointed out to the authors by Dominik Schmid). \red{While the censoring inequality appeared later, these ideas are also implicitly found in the arguments of \cite{MR2630064}.}
    
\subsubsection{Tail Bounds for the Height Function} In this subsection, we state effective hydrodynamic estimates for the fluctuations of the height function $H(x,t)$ of the stochastic six-vertex model with step initial conditions. To do so we first state the law of large numbers for $H$.

With probability one it holds that
\begin{equation}\label{eq:asconverge}
		\lim_{n \to \infty} \frac{H(\lfloor nx \rfloor, \lfloor ny\rfloor)}{n} = g(x, y), \qquad \forall x, y \in \mathbb{R}_{\geq 0},
		\end{equation}
where for $b_1 \leq b_2$, we have
\begin{equation}\label{eq:hlimit}
		g(x, y) = 
		\begin{cases}
  	    y - x & \text{if } \frac{x}{y} \leq \kappa^{-1}\\
		\frac{\big( \sqrt{x}-\sqrt{\kappa y}\big)^2}{\kappa - 1} &\text{if } \kappa^{-1} < \frac{x}{y} < \kappa \\
		0 &\text{if } \frac{x}{y} \geq \kappa
		\end{cases}
	\end{equation}
and for $b_1 \geq b_2$, we have 
\begin{equation*}
g(x, y) = 
\begin{cases}
0 &\text{if } x \geq y\\
y-x &\text{if }x \leq y.
\end{cases}
\end{equation*}
This was proven at the level of weak convergence in \cite{BCGStochasticSixVertex} and \cite{Aggarwal2020} and was strengthened to almost sure convergence in \cite{drillickLinS6V}.

Let $g(x):=g(x,1)$.
We prove the following two tail bounds on the fluctuations of the height function $H$ around its limit shape $g$.

\begin{proposition}[Lower Tail Bound]
\label{prop:MeixnerTail}
Fix $\varepsilon>0$. There exists a constant $c=c(\varepsilon)>0$ such that the following holds: For any $\mu\in[\kappa^{-1}+ \varepsilon,  \kappa - \varepsilon]$ and for any $T\geq1$, $s  \geq 0$, 
\begin{equation} \label{eq:tail_1}
\mathbb{P}\left[H(T\mu, T) \geq g(\mu)T + sT^{1/3}\right] \leq c^{-1} e^{-c s^\frac32},
\end{equation}
and $c$ can be chosen to weakly decrease in $\varepsilon$.
\end{proposition}

\begin{proposition}[Upper Tail Bound]
\label{prop:otherTail}
Fix $\varepsilon>0$. There exists a constant $c=c(\varepsilon)>0$ such that the following holds: For any $\mu\in[\kappa^{-1}+ \varepsilon,  \kappa - \varepsilon]$ and for any $T\geq1$, $s  \geq 0$,  
$$
\mathbb{P}\left[H(T\mu, T) \leq g(\mu)T - sT^{1/3}\right] \leq c^{-1} e^{-c s} ,
$$
and $c$ can be chosen to weakly decrease in $\varepsilon$.
\end{proposition}

\begin{remark}\label{rem:tailbounds}
    The power $T^\frac13$ on the left-hand side of Propositions~\ref{prop:MeixnerTail} and \ref{prop:otherTail} is optimal since on this scale the fluctuations of the height function have been shown to converge to the Tracy-Widom GUE distribution, see \cite[Theorem 1.2]{BCGStochasticSixVertex}.
    The optimal exponents on the right-hand side however, are expected to be $s^{3}$ for Proposition~\ref{prop:MeixnerTail} and $s^{3/2}$ for Proposition~\ref{prop:otherTail}, as was obtained for the longest increasing subsequence of a permutation in \cite{lowe2001moderateupper,lowe2002moderatelower}. The optimality of $s^{3/2}$ for the upper tail was recently confirmed in \cite{landonSosoe2023tail}, see the discussion below. The parameters $\mu_1$ and $\mu_2$ need to be bounded away from the edge of the rarefaction fan in order to obtain a uniform constant $c(\varepsilon)$.
\end{remark}

We call Proposition~\ref{prop:MeixnerTail} a ``lower tail'' bound even though it seemingly describes the upper tail of the random variable $H(T\mu, T)$ since it corresponds to the lower tail of the Tracy-Widom distribution. Similarly, we call Proposition~\ref{prop:otherTail} an ``upper tail'' bound. This better matches the usage in the literature of the terms ``upper" and ``lower" tails for models in the KPZ universality class.  The reason that the upper tail decays more slowly than the lower tail is that for the height function to be smaller than expected, we just need the position of the right-most particle in the stochastic six-vertex model to be small. On the other hand, for the height function to be larger than expected, we must have that the positions of many particles are large. Since this requires more deviations to occur, the probability decays more quickly.

We prove the lower tail bound by using an identity from  \cite{MR3760963} that expresses the $q$-Laplace transform of the height function in terms of an expectation with respect to the law of the Meixner ensemble. This identity allows us to bound the upper tail for the height function by the lower tail of the position of the smallest hole in the Meixner ensemble. The Meixner ensemble is a determinantal point process, so this tail can be expressed as a Fredholm determinant, which we then bound using Widom's trick \cite{widom2002convergence}. The upper tail bound is more straightforward. We directly express the $q$-Laplace transform of the height function in terms of a Fredholm determinant and use Fredholm determinant estimates from \cite{AB2019Aseps6vphase}.

The above tail bounds are in the ``moderate deviations" regime since we are considering fluctuations of order $T^{1/3}$, in contrast with the ``large deviations" regime, which considers fluctuations of order $T$. For TASEP, tail bounds in the moderate deviations regime have been obtained in \cite{MR3304747, basu2014last} (with some of the key ideas originating in \cite{BaikFerrariPeche2014twopointTASEP}). For ASEP, these bounds were obtained in \cite{ACG2023asepspeed}. 

\red{We briefly compare our proof methods to those used for ASEP in  \cite{ACG2023asepspeed}. In \cite{ACG2023asepspeed} the lower tail bound was proven by exploiting the connection between ASEP and the Laguerre ensemble to reduce the tail bound for ASEP to known tail bounds for the smallest hole of the Laguerre ensemble. Since analogous bounds were not already available for the smallest hole of the Meixner ensemble, we used Widom's trick to prove them. Our proof of the upper tail on the other hand is a relatively straightforward generalization of the proof of the upper tail in \cite{ACG2023asepspeed}.}

While our paper was in progress, \cite{landonSosoe2023tail} considered the stochastic six-vertex model under stationary initial conditions and obtained tail bounds in the moderate deviations regime. They also obtained \cite[Theorem 2.7]{landonSosoe2023tail} an upper tail bound with the optimal exponent of $s^{3/2}$ and optimal constants for step initial conditions by developing the Rains-EJS formula for the stochastic six-vertex model. Their proof differs from our proof of Proposition \ref{prop:otherTail} in that it does not rely on integrable methods, and instead uses probabilistic couplings, building on ideas developed in \cite{emrah2020righttailmoderatedeviationsexponential, MR4620410, MR4610276}. We still include our proof since the methods are substantially different. There is also an upcoming work \cite{GhosalSilva2024} that will prove tight upper and lower tail bounds with the optimal exponents in the moderate deviations regime using Riemann-Hilbert techniques. Going to the large deviations regime, the recent work \cite{dasLiaoMucciconi2024lower} obtains a large deviation principle for the lower tail of the stochastic six-vertex model under step initial conditions. 

Let us briefly summarize some other fluctuation results for the stochastic six-vertex model with step initial conditions. As mentioned earlier, the stochastic six-vertex model belongs to the KPZ universality class (see \cite{MR2930377} for a survey of this area), and in particular, this means that the height function exhibits fluctuations of scale $T^\frac13$ and correlations of scale $T^\frac23$ on a domain of size $T$. In \cite{BCGStochasticSixVertex} they proved that the one-point fluctuations of the height function around its limit shape are of the order $T^{1/3}$ and are asymptotically given by the Tracy-Widom distribution. To go beyond a one-point result, we can view the height function as a spatial process with space rescaled by $T^{2/3}$ and the fluctuations by $T^{1/3}$. This process was shown to be tight in \cite{CorwinDimitrov2018Transversalfluctuations} and its two-point distribution converges to the two-point distribution of the Airy process \cite{MR4561796}. 
Finally, the recent groundbreaking work of \cite{aggarwal2024scaling} proved the convergence of the height function when viewed as a four-parameter field to the directed landscape, fully confirming that the stochastic six-vertex model is in the KPZ universality class.

\subsubsection{Proof Sketch}\label{sec:sketch} We now sketch the proof of Theorem \ref{thm:main} using the above two ingredients. To show that the speed $\frac{\bm{X}_t}{t}$ converges a.s., we will introduce a sequence of times $S_n$ and prove that as long as we are not too close to the edge of the rarefaction fan, then with high probability, 
\begin{equation} \label{eq:controlSpeed}
\left|\frac{\bm{X}_{S_n}}{S_n}-\frac{\bm{X}_{S_{n+1}}}{S_{n+1}}\right| \leq S_n^{-\gamma}
\end{equation}
for some positive $\gamma$.

For this to imply convergence of the sequence $\frac{\bm{X}_{S_n}}{S_n}$, we need the right-hand side to be summable. For general times $S_n \leq t \leq S_{n+1}$, one can then use the monotonicity of $\bm{X}_t$ to bound $\left|\frac{\bm{X}_t}{t}-\frac{\bm{X}_{S_n}}{S_n}\right|$ as long as the sequence $S_n$ does not grow too quickly. We will take the sequence $S_{n+1}=S_n+T(S_n):=S_n+S_n^{\frac{7}{9}}$ and prove the necessary bounds on the probability of \eqref{eq:controlSpeed}, given certain events in Proposition \ref{prop:straightlinespoly}.

To prove Proposition \ref{prop:straightlinespoly}, we want to control the behavior of the second-class particle after some large initial time $S_0$. However, the effective hydrodynamic bounds in Propositions \ref{prop:MeixnerTail} and \ref{prop:otherTail} only allow us to control the behavior of a large number of particles, not of an individual one since they are mesoscopic statements as opposed to microscopic ones. Therefore, we fill up all empty positions to the left of $\bm{X}_S$ with third-class particles and control the union of the second- and third-class particles by Propositions \ref{prop:MeixnerTail} and \ref{prop:otherTail}. Letting  $T=T(S)=S^{\frac{7}{9}}$, Theorem~\ref{thm:orderinginequality} will guarantee that only a small number of these third-class particles will be to the right of $\bm{X}_{S+T}$ at time $S+T$, so that controlling the union of the second- and third-class particles gives us a bound on $\bm{X}_{S+T}$.

\begin{figure}
    \centering
    \includegraphics{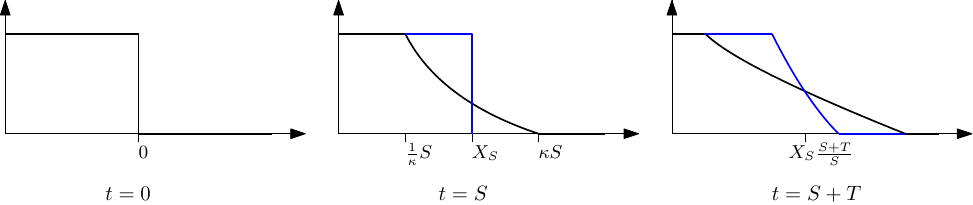}
    \caption{A sketch of the densities of the processes $\B^{(1)}$ in black at times $0,S$ and $S+T$ and $\B^\text{step}$ in blue at times $S$ and $S+T$.
    At time $S$ the process $\B^{(1,2,3)}$ is given exactly by the maximum of the two processes $\B^{(1)}$ and $\B^{\text{step}}$, while at time $S+T$ it is at least the maximum of $\B^{(1)}$ and $\B^{\text{step}}$.}
    \label{fig:densities}
\end{figure}

We split the proof of \eqref{eq:controlSpeed} into an upper and a lower bound, which are treated analogously.
Theorem~\ref{thm:orderinginequality} reduces the lower bound to showing that a large number of these second- and third-class particles are to the right of $\bm{X}_S+\frac{\bm{X}_S}{S}T-S^{1-\gamma}$ at time $S+T$.
To do so denote by $\B^{(1,2,3)}$ the augmented (single-class) stochastic six-vertex model containing the union of all first-, second- and third-class particles and by $\B^{(1)}$ the process with only the first-class particles.
Additionally, we introduce an auxiliary third process $\B^\text{step}$ which is started at time $S$ from the initial condition $\B^{\text{step}}_S(x)=\bm{1}_{x\leq \bm{X}_S}$.
At time $S$ these three processes satisfy
\begin{equation}
    \B^{(1,2,3)}_S(x)=\max(\B^{(1)}_S(x),\B^{\text{step}}_S(x))\,.
\end{equation}
The multi-class stochastic six-vertex process allows us to couple $\B^{(1,2,3)}$ and $\B^\text{step}$ such that at any later time $S+T$ it holds that $\B^{(1,2,3)}_{S+t}(x)\geq\B^\text{step}_{S+t}(x)$.
Since $\B^{(1,2,3)}$ and $\B^{(1)}$ are already coupled in such a way, this implies that for any $t\geq0$
\begin{equation}
    \B^{(1,2,3)}_{S+t}(x)\geq \max(\B^{(1)}_{S+t}(x),\B^{\text{step}}_{S+t}(x))\,.
\end{equation}
Note that this also couples $\B^{(1)}$ and $\B^{\text{step}}$ in some non-trivial way. See Figure \ref{fig:densities} for a sketch of the particle densities of the processes $\B^{(1)}$ and $\B^{\text{step}}$ at times $0, S$, and $S+T$. 

By using the effective hydrodynamic estimates together with a recent approximate monotonicity result from \cite{aggarwal2024scaling}, we show that with high probability $\B^{(1)}$ is still close to the hydrodynamic limit at time $S+T$, uniformly over all possible configurations of $\B^1_S$ when on a certain event $H_S$, which also occurs with high probability.
Since the process $\B^\text{step}$ is started from step initial conditions, it is also close to a hydrodynamic limit at time $S+T$, which is obtained by translating the hydrodynamic limit for standard step initial conditions.
By the coupling above
\begin{equation}
    \B^{(1,2,3)}_{S+T}(x)-\B^{(1)}_{S+T}(x)\geq\B^{\text{step}}_{S+T}(x)-\B^{(1)}_{S+T}(x).
\end{equation}
Using the hydrodynamic estimates for the two processes on the right-hand side, this gives a lower bound for the number of third-class particles to the right of $\frac{X_S}{S}(S+T)-S^{1-\gamma}$, as desired.

The proof of Theorem \ref{eq:fluct} uses similar arguments as above to bound $\left|\bm{X}_{S_n}-S_n\bm{U}\right|$. To bound   $\left|\bm{X}_{t}-s\bm{U}\right|$ for general times $S_n \leq t \leq S_{n+1}$, we will need to bound $|S_{n+1} - S_n|$ and it is here that we make crucial use of the exact form of the sequence $S_{n+1} = S_n + S_n^{\frac{7}{9}}$ to obtain the $\frac{7}{9}$ fluctuation exponent in Theorem \ref{thm:fluct}. 

This proof represents a significant streamlining of the strategy employed in \cite{ACG2023asepspeed}, which is possible in part because Theorem~\ref{thm:orderinginequality} is stronger than Rezakhanlou's coupling. In particular, the fact that Theorem~\ref{thm:orderinginequality} does not get worse with the number of particles allows us to fill in all empty positions to the left of the second-class particle with third-class particles.
In \cite{ACG2023asepspeed} only a small number of positions were filled, which made it necessary to deal with more complicated ``$\varphi$-distributed'' Bernoulli initial conditions and introduced a further approximation step. Furthermore, we are able to obtain Theorem~\ref{thm:fluct} due to taking a more refined time scale, and more careful analysis. See Remark \ref{rmk:sn} for further discussion of this. 
% Since Theorem~\ref{thm:orderinginequality} also holds for ASEP, these changes could be used to obtain Theorem~\ref{thm:fluct} for ASEP as well, with the same exponent.

\subsection{Summary} In summary, the purpose of this paper is four-fold:
\begin{itemize}
    \item We obtain a new tool in Theorem~\ref{thm:orderinginequality} for studying the stochastic six-vertex model and ASEP by adding higher class particles.
    Since this article was announced, this tool has already been used in \cite{aggarwal2024kpzfixedpointconvergence} to study the convergence of the stochastic six-vertex model and ASEP to the KPZ fixed point for general initial conditions.
    \item We give bounds on the deviations of the height-function of the stochastic six vertex model in the moderate deviation regime in Proposition~\ref{prop:MeixnerTail} and Proposition~\ref{prop:otherTail}. Our results contribute to a growing body of recent literature developing tools to prove lower tail bounds for \emph{positive temperature} KPZ models, see e.g., \cite{MR4094738, emrah2020righttailmoderatedeviationsexponential,MR4318370, MR4294287, MR4373176, MR4620410, MR4610276,  MR4709103, shen2024timecorrelationsinversegammapolymer, MR4771973}.
    \item We give bounds on the fluctuations of the second class particle started in step initial conditions in Theorems~\ref{thm:main} and \ref{thm:fluct}, in particular showing almost sure convergence of its speed.
    \item We show the existence of the stochastic six-vertex speed process, and show that this process is ergodic and stationary for the dynamics of the stochastic six-vertex model in Section~\ref{sec:speedsymmetries}.
\end{itemize}

\subsection{Structure}
In Section \ref{sec:couplings}, we recall some couplings and properties of the stochastic six-vertex model, including the approximate monotonicity result from \cite{aggarwal2024scaling} which is stated in Proposition \ref{prop:approxmono}. The two core ingredients are proved in Sections \ref{sec:orderinginequality} and \ref{sec:effectivehydrodynamics} respectively---in Section \ref{sec:orderinginequality} we prove Theorem~\ref{thm:orderinginequality} and in Section \ref{sec:effectivehydrodynamics} we prove Propositions \ref{prop:MeixnerTail} and \ref{prop:otherTail}.

These results are then used in Sections \ref{sec:linear}, \ref{sec:hydroToLinear} and \ref{sec:hydroevents} to prove the main theorem.
In order these sections show that
\begin{itemize}
    \item the main theorem follows if one can show that with high probability the second-class particle does not deviate too much from its current speed in a given time frame,
    \item which follows if one can show that the augmented progress with additional third-class particles does not deviate too much from its hydrodynamic limit with high probability,
    \item which follows from the effective hydrodynamics from Section~\ref{sec:effectivehydrodynamics} together with approximate monotonicity.
\end{itemize}
Finally in Section~\ref{sec:speedsymmetries} the existence of the speed process is deduced from Theorem~\ref{thm:main}, and we study some of its properties. 

\subsection{Notation}
Throughout the paper, many floor functions are dropped when we consider large integers.
We use
\[
\llbracket A,B\rrbracket=[A,B]\cap\Z
\]
for intervals of integers.
\red{All time variables are integer-valued, e.g., when we write that something holds for all times $t\geq0$, we mean it to hold for all non-negative integers.}

Our convention for geometric random variables is that a random variable $X\sim\text{Geo}(q)$ satisfies
\[
\mathbb P[X=k]=(1-q)q^k\,.
\]

We consider both single-class and multi-class processes by considering their occupation variables.
Single-class processes have occupation variables in $\{0,1\}$ with $0$ encoding holes and $1$ encoding particles, while multi-class processes have occupation variables in $\R\cup\{\infty\}$, with $\infty$ encoding holes and all other values encoding particles of different classes.

The parameters $b_1$ and $b_2$ are fixed throughout the paper and therefore all constants can depend on them freely even if this is not explicitly mentioned.

\subsection{Acknowledgements}

The authors thank Ivan Corwin and Amol Aggarwal for multiple discussions and guidance.
In particular, we want to thank Ivan Corwin for suggesting this problem to us and for his course at the 2022 PIMS-CRM Summer School in Probability, where this project was initiated and where he presented a version of Theorem~\ref{thm:orderinginequality} which was limited to $b_1\leq\frac{1}{2}$. We thank the summer school organizers for their hospitality and acknowledge the support from NSF DMS-1952466. We thank Dominik Schmid for pointing out the connection between our stochastic domination result and the censoring inequality in the setting of ASEP. We thank Promit Ghosal for sharing his forthcoming work \cite{GhosalSilva2024} with us and Evan Sorensen for explaining the connection between the stochastic six-vertex model speed process and the stationary horizon. 
Part of this work was compiled during a research visit of the second author to Columbia University in the fall of 2023. HD’s research was supported by the NSF Graduate Research Fellowship under Grant No. DGE-2036197, Ivan Corwin’s NSF grants DMS-1811143, DMS-1664650, DMS-1937254, and DMS-2246576, the W.M. Keck Foundation Science and Engineering Grant on “Extreme diffusion", as well as the Fernholz Foundation. This material is partly based upon work supported by the National Science Foundation under Grant No. DMS-2424139, while HD was in residence at the Simons Laufer Mathematical Sciences Institute in Berkeley, California, during the Fall 2025 semester.

\section{The basic coupling and some properties of the stochastic six-vertex model}\label{sec:couplings}

We consider the following construction of the single-class stochastic six-vertex model, which also allows us to couple multiple stochastic six-vertex models with varying boundary conditions.
We will first state it on the quadrant.
\begin{definition}[Basic Coupling]\label{def:basiccoupling}
    Consider two independent families $(\chi^1(x,t))_{x,t\geq 0}$ and  $(\chi^2(x,t))_{x,t\geq 0}$ of i.i.d. Bernoulli$(b_1)$ and Bernoulli$(b_2)$ random variables respectively.
    Given such random variables, we can sample the stochastic six-vertex model in the following way.
    If at a given vertex there are either two incoming arrows or no incoming arrows then there is only a single possible outcome.
    If there is a single incoming vertical arrow at $(x,t)$ and $\chi^1(x,t)=1$, then the outgoing arrow is vertical. If $\chi^1(x,t)=0$, then the outgoing arrow is horizontal.
    Similarly, if there is a single incoming vertical arrow and $\chi^2(x,t)=1$, then the outgoing arrow is vertical. If $\chi^2(x,t)=0$, then the outgoing arrow is horizontal. 

    Given boundary conditions on the left and bottom edge of $\Z_{\geq0}\times\Z_{\geq0}$ the random variables $(\chi^1(x,t))_{x,t\geq 0}$ and  $(\chi^2(x,t))_{x,t\geq 0}$ uniquely define a configuration, which can be obtained by updating the vertices along the anti-diagonal lines $\{(x,t):x+t=k\}$ with increasing $k$.
    Note also that the order of updates does not matter.
    Using the same $(\chi^1(x,t))_{x,t\geq 0}$ and  $(\chi^2(x,t))_{x,t\geq 0}$ for different boundary conditions gives a coupling of stochastic six-vertex models, which we call the \textbf{basic coupling}.
\end{definition}
Some useful properties of this coupling were recently developed in \cite{aggarwal2024scaling}. 
Before we discuss the properties of this coupling, let us show how it can be used to define the stochastic six-vertex-process on the line in a way that is similar to both the construction in \cite[Section 2.1]{Aggarwal2020} using a different coupling of the stochastic six-vertex model and to the graphical construction of Harris for ASEP on $\Z$ in \cite{Harris1978Graphical}.
\begin{proposition}[Extension to $\Z$]
    The construction in Definition~\ref{def:basiccoupling} can be extended to the domain $\Z\times\Z_{\geq0}$.
    More specifically given two independent families $(\chi^1(x,t))_{x\in\Z,t\geq 0}$ and  $(\chi^2(x,t))_{x\in\Z,t\geq 0}$ of i.i.d. Bernoulli$(b_1)$ and Bernoulli$(b_2)$ and any boundary conditions on the incoming edges of $\Z\times\{0\}$, there is almost surely a unique configuration on $\Z\times\Z_{\geq0}$ that is coherent with the boundary conditions and that at each vertex satisfies the rules in Definition~\ref{def:basiccoupling}, i.e. if there is only one incoming arrow, the configuration at the vertex $(x,t)$ is given by the values of $\chi^1(x,t)$ and $\chi^2(x,t)$.
    Furthermore, the law of this unique configuration is given by the stochastic six-vertex model.
\end{proposition}
\begin{proof}
    We will construct the configuration line by line.
    Consider first the random variables $\chi^1(x,0)$ and $\chi^2(x,0)$.
    We call a vertex $(x,0)$ such that $\chi^1(x,0)=\chi^2(x,0)=0$ a \textbf{cut-vertex}.
    Almost surely, there are infinitely many cut-vertices both to the left and to the right of the origin since each vertex $(x,0)$ has an independent positive probability of $(1-b_1)(1-b_2)$ to be a cut-vertex.
    Notice that at a cut-vertex, the outgoing horizontal edge is occupied if and only if the incoming vertical edge is occupied, and the outgoing vertical edge is occupied if and only if the incoming horizontal edge is occupied.
    Therefore, if $(x_0,0)$ and $(x_1,0)$ with $x_0<x_1$ are cut-vertices, the configuration of all vertices $(x,0)$ with $x_0<x\leq x_1$ is determined by the incoming arrows at these vertices and the Bernoulli variables $\chi^1(x,0)$ and $\chi^2(x,0)$ for $x_0\leq x\leq x_1$.
    Therefore on the probability $1$ event that there are cut-vertices infinitely far to the left, the configuration is uniquely determined.
\end{proof}
Again, using the same Bernoulli random variables for different initial conditions gives a coupling of stochastic six-vertex processes. Let us now consider several properties of this coupling starting with attractivity. As mentioned in the introduction, we will use the notation $(\eta_t(x))_{x\in\Z,t\geq0}$ for the occupation variables, i.e. $\eta_t(x)=1$ if the vertical incoming edge is occupied.
The initial conditions are then given by a function $\eta_0(x):\Z\to\{0,1\}$.
\begin{lemma}[Attractivity]\label{prop:attractivity}
    Given a collection of initial conditions $\eta_0^k$ for $k=1,\dots,n$, such that $\eta_0^i(x)\leq\eta_0^j(x)$ for $i\leq j$ and all $x\in\Z$, under the basic coupling it will hold that $\eta_t^i(x)\leq\eta_t^j(x)$ for all $t\in\Z_{\geq0}$ and $x\in\Z$.
\end{lemma}
\begin{proof}
    Let us consider $\eta^i$ and $\eta^j$.
    Assume that the desired property is true until updating a specific vertex.
    If at this vertex the incoming arrows are identical for $\eta^i$ and $\eta^j$, by the coupling the outgoing arrows will also be identical.
    If they are not, since the property holds for all the previous steps, either there are two incoming arrows in $\eta^i$ or no incoming arrows in $\eta^j$.
    In either case, the outgoing arrows will also still satisfy the desired condition.
\end{proof}
\begin{remark}
    Note that the basic coupling with initial conditions $\eta_0^k$ for $k=1,\dots,n$, such that $\eta_0^i(x)\leq\eta_0^j(x)$ for $i\leq j$ and all $x\in\Z$, exactly corresponds to the $n+1$-class stochastic six-vertex model with classes $\{1,\dots,n,\infty\}$ in the following way.
    Define
    \[\eta^\text{mult}_t(x)= \min\{i\in \{1,\dots n\}:\eta^i_t(x)=1\}\,,
    \]
    where the convention is used that the minimum of the empty set is $\infty$.
    By considering the possible situations at a single vertex, one easily checks that $\eta^\text{mult}_t$ is a multi-class stochastic six-vertex process.
\end{remark}

The attractivity property also has the following analogue for the multi-class process. 
\begin{lemma}[Merging]\label{lem:merging}
    Let $(\eta_t)_{t\in\Z_{\geq0}}$ be a multi-class stochastic six-vertex model with classes in $\Z\cup\{-\infty,\infty\}$, i.e. $\eta_t:\Z\to\Z\cup\{-\infty,\infty\}$.
    Then for any weakly increasing function $\phi:\Z\cup\{-\infty,\infty\}\to\Z\cup\{-\infty,\infty\}$, the process $(\phi\circ\eta_t)_{t\in\Z_{\geq0}}$ is also a multi-class stochastic six-vertex model.
\end{lemma}
\begin{proof}
    This is an immediate consequence of the weights in Figure~\ref{fig:multiClasss6v} only depending on the incoming classes $i$ and $j$ via their ordering.
    Consider a vertex for which an update is about to be performed.
    If the two incoming classes $i$ and $j$ are equal, they will also be equal after applying the map, and in either case, there is exactly one outcome which then of course has probability $1$.
    If the two incoming classes $i$ and $j$ are different, i.e. $i<j$ (note that we do not assume whether $i$ is the horizontal or vertical incoming arrow), then either $\phi(i)<\phi(j)$ or $\phi(i)=\phi(j)$.
    In the first case, there are two possible outcomes for both a vertex with incoming arrows $i$ and $j$ and a vertex with incoming arrows $\phi(i)$ and $\phi(j)$ and the probabilities match, since the relative order of the incoming arrows is the same.
    In the second case there are two possible outcomes before applying $\phi$ but only one outcome after applying $\phi$.
    Since the two possibilities before applying $\phi$ are complementary, their probabilities sum up to $1$, which is the probability of the one possible outcome after applying $\phi$.
\end{proof}

\red{Another natural property that one would like to have for a coupling of particle systems is monotonicity. We say that a coupling of the stochastic six-vertex model is \emph{monotone} if when we couple together two processes $\eta^1_t,\eta^2_t$  whose initial height functions are ordered, then this ordering is preserved for all time under the coupling. In other words if $h_0(x; \eta^1) \leq h_0(x; \eta^2)$ for all $x \in \mathbb Z$, then almost surely we have $h_t(x; \eta^1) \leq h_t(x; \eta^2)$ for all $t> 0$ and for all $x \in \mathbb Z$.  However, while the basic coupling for ASEP is monotone, the basic coupling for the stochastic six-vertex model is \emph{not monotone} \footnote{There is a monotone coupling for the stochastic six-vertex model (see e.g., \cite{Aggarwal2020}); however, this coupling cannot handle more than two initial conditions simultaneously, i.e., it is not a \emph{global monotone} coupling.}. This leads to an extra layer of complexity when proving results for the stochastic six-vertex model.} As monotonicity properties will be crucial in our proof, we will use a recent result from \cite{aggarwal2024scaling} which gives an approximate form of monotonicity for the basic coupling. Recall that given a stochastic six-vertex process $(\eta_t(x))_{x\in\Z,t\geq0}$, there is a height function $h_t(x;\eta)$ defined up to a global shift defined in Definition~\ref{def:heightfunction} The following proposition is Lemma D.3 of \cite{aggarwal2024scaling}.

\begin{proposition}[Approximate Monotonicity]\label{prop:approxmono}
   Consider two single-class initial conditions $\eta^1:\Z\to\{0,1\}$ and $\eta^2:\Z\to\{0,1\}$ both with at most $N$ particles .
    Further consider height functions $h_t(x;\eta^1)$ and $h_t(x;\eta^2)$ satisfying $h_t(x;\eta^1)=h_t(x;\eta^2)=0$ for $x$ large enough.
    If $M\geq (log N)^2$ and $|h_0(x;\eta^1)-h_0(x;\eta^2)|<K$ for all $x\in\Z$, and $t\geq0$, then with probability at least $1-c^{-1}e^{-cM}$, and for all $x\in\Z$ it holds that
    \[
        |h_t(x;\eta^1)-h_t(x;\eta^2)|\leq K+M\,.
    \]
\end{proposition}
\begin{remark}
In \cite{aggarwal2024scaling} this is stated without the absolute value.
However, the basic coupling has the following property:
If $(\eta^1,\eta^2)$ are two stochastic six-vertex processes coupled using the basic coupling so are $(\eta^2,\eta^1)$. Note that this is a property that the monotone coupling in \cite[Proposition 2.6]{Aggarwal2020} does not have.
Additionally, the conditions on $\eta^1$ and $\eta^2$ are symmetric, and therefore the statement with the absolute value follows from the statement without the absolute value by a simple union bound.
\end{remark}

Another property that we will need is a special case of \cite[Lemma D.4]{aggarwal2024scaling}, and the proof is quite similar to \cite[Proposition 2.17]{Aggarwal2020}.
\begin{proposition}[Finite Speed of Discrepancies]\label{prop:finitespeed}
    There exists a constant $c=c(b_2)>0$ depending only on $b_2$ such that the following holds.
    Consider two particle configurations $\eta_0$ and $\xi_0$ with height function $h_0(x;\eta)$ and $h_0(x;\xi)$which are equal on some interval $\llbracket A,B \rrbracket$.
    Then, under the basic coupling, with probability at least $1-c^{-1}e^{-cT}$ it holds that $h_t(x;\eta)=h_t(x;\xi)$ for all $t\leq T$ and all $x\in\left\llbracket A+\frac{2T}{1-b_2}+1,B\right\rrbracket$. 
\end{proposition}

Using Propositions~\ref{prop:approxmono} and  \ref{prop:finitespeed} together, we can show that given two initial conditions with height functions close on an interval, the height functions will stay close on a smaller interval for some time.
\begin{lemma}[Approximate Monotonicity on Intervals] \label{lem:finitiemonotonicity}
There exists a constant $c=c>0$, depending only on $b_1,b_2\in(0,1)$, such that the following holds.
Consider two particle configurations $\eta_0$ and $\xi_0$ with height functions $h(x;\eta_0)$ and $h(x;\xi_0)$ such that for $x\in\llbracket A,B \rrbracket$ we have $|h_0(x;\xi)- h_0(x;\eta)|\leq K$.
Let $M\geq \log(B-A)^2$.
Then we can couple them such that with probability at least $1-c^{-1}(e^{-cT}+e^{-cM})$ it holds that $|h_T(x;\xi)- h_T(x;\eta)|\leq 3K+M$ for all $x\in\llbracket A+\frac{2T}{1-b_2}+1,B\rrbracket $.
\end{lemma}
\begin{proof}
This will follow from Propositions \ref{prop:approxmono} and \ref{prop:finitespeed}.
Let $\widetilde{\eta}_0$ be the particle configuration obtained from $\eta_0$ by setting
\begin{equation}
    \widetilde{\eta}_0(x)=
    \begin{cases}
        0 &\text{if } x<A\\
        \eta(x)&\text{if } x\in\llbracket A,B \rrbracket\\
        0&\text{if } x>B\,,
    \end{cases}
\end{equation}
and define $\widetilde{\xi}_0$ in the same way.
Couple $\eta,\xi, \widetilde{\eta}$ and $\widetilde{\xi}$ all with one basic coupling (i.e. all using the same i.i.d. Bernoulli random variables).
Let the height functions $h_0(x;\widetilde{\eta})$ and $h_0(x;\widetilde{\xi})$ be chosen such that $h_0(B;\widetilde{\eta})=h_0(B;\widetilde{\xi})=0$, i.e. $h_0(x;\widetilde{\eta})=h_0(x;\eta)-h_0(B,\eta)$ for $x\in \llbracket A,B \rrbracket$ and the same for $\xi$.
Note that $h_0(x;\eta)-h_0(B,\eta)$ is a height function for $\eta_0$, and therefore by applying Proposition~\ref{prop:finitespeed} twice, once for $\eta$ and once for $\xi$, and a union bound, we obtain that
\begin{equation}\label{eq:finitesuccess}
h_t(x,\widetilde{\eta})=h_t(x,\eta)-h_0(B,\eta)\text{ and } h_t(x,\widetilde{\xi})=h_t(x,\xi)-h_0(B,\xi)
\end{equation}
holds for all $t\leq T$ and all $x\in\llbracket A+\frac{2T}{1-b_2},B\rrbracket$ with probability at least $1-c^{-1}e^{-cT}$.

Further note that at time $0$, for all $x$
\begin{equation}
     |h_0(x;\widetilde{\eta})-h_0(x;\widetilde{\xi})|=|(h_0(x;\eta)-h_0(B,\eta))-(h_0(x;\xi)-h_0(B;\xi))|\leq 2K\,.
\end{equation}
Therefore we can apply Proposition \ref{prop:approxmono} to $\widetilde{\xi}$ and $\widetilde{\eta}$ since they are coupled with the basic coupling.
Indeed both $\widetilde{\xi}$ and $\widetilde{\eta}$ have at most $B-A$ particles each, so we will have with probability at least $1-c^{-1}e^{-cM}$ that 
\begin{equation}\label{eq:monotonicitysuccess}
|h_T(x;\widetilde{\eta})-h_T(x;\widetilde{\xi})|\leq 2K+M\text{ for all }x\in\Z.
\end{equation}
By a union bound, with probability at least $1-c^{-1}(e^{-cT}+e^{-cM})$ both events \eqref{eq:finitesuccess} and \eqref{eq:monotonicitysuccess} take place. 
On this event it holds for all $x\in\llbracket A+\frac{2T}{1-b_2}+1,B\rrbracket$ that 
\begin{align*}
|h_T(x,\eta)-h_T(x,\xi)|&\leq |h_T(x,\eta)-h_T(x,\widetilde{\eta})+h_T(x,\widetilde{\xi})-h_T(x,\xi)+h_T(x,\widetilde{\eta})-h_T(x,\widetilde{\xi})|\\
&\leq  |h_T(x,\eta)-h_T(x,\widetilde{\eta})-h_T(x,\widetilde{\xi})+h_T(x,\xi)|+|h_T(x,\widetilde{\eta})-h_T(x,\widetilde{\xi})|\\
&=|-h_0(B,\eta)+h_0(B,\xi)|+|h_T(x,\widetilde{\eta})-h_T(x,\widetilde{\xi})|\\
&\leq K+2K+M=3K+M\,,
\end{align*}
where we are used a triangular inequality, \eqref{eq:finitesuccess}, \eqref{eq:monotonicitysuccess} and that $|-h_0(B,\eta)+h_0(B,\xi)|<K$ by the assumption on the height functions at time $0$.
\end{proof}
\begin{remark}
    The factor $3$ in the term $3K+M$ in the previous step is an artifact of Proposition~\ref{prop:approxmono} being stated only for height functions which are $0$ far enough to the right.
    This restriction could be removed, which would remove the factor $3$.
    However, for our purposes the above is sufficient.
\end{remark}

This property will be used in Proposition~\ref{prop:hydroEvent}, to show that if a stochastic six-vertex process $\eta$ is close to its hydrodynamic limit at time $S$, it will still be close to its hydrodynamic limit at time $S+T$ with high probability, even conditioned on its full configuration at time $S$.

Finally, the stochastic six-vertex model has the following two symmetries which are often used together. 
\begin{proposition}[Particle-Hole Inversion]\label{prop:particle-hole} 
    If we interchange all particles and holes in a stochastic six-vertex process, we obtain another stochastic six-vertex process, but with $b_1$ and $b_2$ swapped. 
\end{proposition}
\begin{proposition}[Space Inversion]\label{prop:space}
    If we exchange the two coordinate axes in a stochastic six-vertex process, we obtain another stochastic six-vertex process, but with $b_1$ and $b_2$ swapped. 
\end{proposition}
\begin{proof}
    Both of these can be seen by looking at what happens to the six configurations in Figure \ref{fig:s6v} under this inversion. 
\end{proof}
Using both of these symmetries on the quadrant, which is symmetric with respect to the line $x=t$, we obtain a symmetry of one stochastic six-vertex model with itself.
In particular one can see that the law of the stochastic six-vertex model started from step initial conditions on the quadrant is invariant after applying both inversions.
Furthermore, the step initial condition with a single particle coming in at the origin from the left is dual to step initial conditions with a single particle coming in at the origin from the bottom.
Therefore it suffices to prove the main theorem for this kind of initial condition.

\section{Number of overtaking third-class particles}\label{sec:orderinginequality}

The purpose of this section is to prove Theorem~\ref{thm:orderinginequality} which will allow us to control an individual second-class particle by controlling a large number of third-class particles.

\begin{figure}
    \centering
    \includegraphics[width=0.8\linewidth]{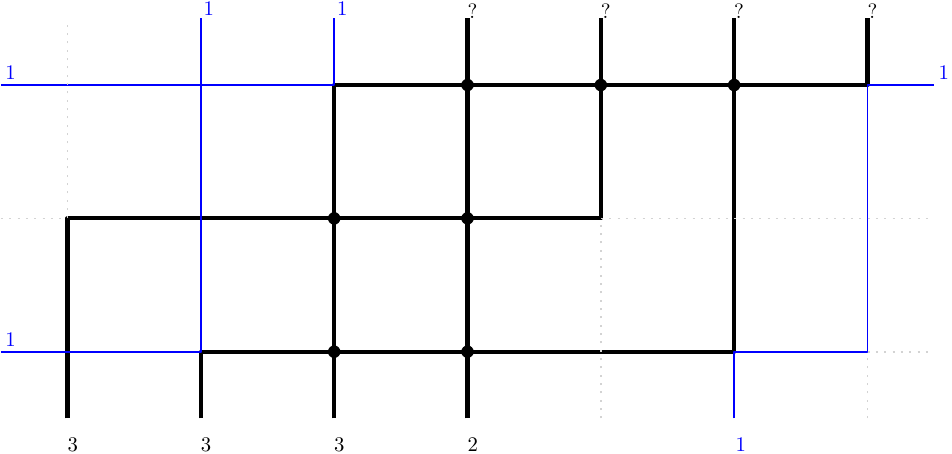}
    \caption{A sample configuration for the proof of Theorem~\ref{thm:orderinginequality}.
    The black lines depict the paths of the second-and third-class particles, and the blue lines depict the paths of the first-class particles.
    The classes of the incoming particles are given below the black lines, and correspond to $L_0=0$.
    The black dots mark the vertices at which both incoming particles are second- and third-class particles.
    The law of the position of the second class particle at the top of this configuration $L_t$ is given by $P_0P_1P_2P_1P_2P_0P_1e_0$.}
    \label{fig:configuration}
\end{figure}

\begin{proof}[Proof of Theorem \ref{thm:orderinginequality}]
\red{
Fix $t > 0$. Let $\bm{Z}=(\bm{Z}_{s}(k))_{0\leq s\leq t,0\leq k\leq M}$. As stated, we will prove the statement conditioning on both $\bm{Z}$ and the space-time history of the first-class particles up to time $t$, i.e., we will show that for every possible realization of $\bm{Z}$ and the positions of the first-class particles, the law of $L_t$ conditioned on this realization is dominated by $\text{Geo}(q)$.
We will first condition on $\bm{Z}$. See Figure \ref{fig:configuration} for a realization of $\bm{Z}$ (given by the black edges in the figure).
Given $\bm{Z}$, there is a finite number $W=W(\bm{Z})$ of vertices $(x_i,t_i)_{1\leq i\leq W}$ at which the incoming edges are both second- or third-class particles.
Let the ordering of these vertices be row by row, from left to right, i.e. $t_i$ is weakly increasing, and $x_i<x_j$, when $i<j$ and $t_i=t_j$.
For a specific example, these vertices are marked in Figure~\ref{fig:configuration}.
The value of $L_t$ is determined by $\bm{Z}$ together with the values of the Bernoulli random variables $(\chi^1(x_i,t_i),\chi^2(x_i,t_i))_{1\leq i\leq W(\bm{Z})}$ (recall Definition~\ref{def:basiccoupling}).
Indeed, the event $\{\bm{Z}=Z\}$ for any configuration $Z$ is independent of $(\chi^1(x_i,t_i),\chi^2(x_i,t_i))_{1\leq i\leq W(Z)}$ and therefore the law of $(\chi^1(x_i,t_i),\chi^2(x_i,t_i))_{1\leq i\leq W(\bm{Z})}$ is still that of i.i.d. Bernoulli random variables with parameters $b_1$ for $\chi^1$ and $b_2$ for $\chi^2$ after conditioning on $\bm{Z}$.
Conceptually, this is because which of the higher class particles is the second-class particle does not matter for the interactions between the first and higher class particles.}

\red{Note that after conditioning on $\bm{Z}$, $L_t$ is now independent of the space-time history of the first-class particles up to time $t$, so we can condition on that as well without changing the law of $L_t$.}

\red{
We now introduce an auxiliary process $(\widetilde{L}_i)_{0\leq i\leq W}$ that encodes the vertex by vertex updates to $L$.
First, let $\widetilde{L}_0=0$.
For each $1 \leq i \leq W$, let $k_i$ be such that $\bm{Z}_{t_i}(k_i)=x_i$.
Now obtain $\widetilde{L}_i$ as follows:
\begin{itemize}
    \item If $\widetilde{L}_{i-1}=k_i$ then $$\widetilde{L}_i= \begin{cases} k_i & \text{if } \chi^1(x_i,t_i)=0 \\ k_i+1 & \text{if } \chi^1(x_i,t_i)=1 \end{cases}.$$
    \item If $\widetilde{L}_{i-1}=k_i+1$ then 
    $$\widetilde{L}_i= \begin{cases} k_i & \text{if } \chi^2(x_i,t_i)=1 \\ k_i+1 & \text{if } \chi^2(x_i,t_i)=0 \end{cases}.$$
    \item if $\widetilde{L}_{i-1}$ is neither $k_i$ nor $k_i+1$, then $\widetilde{L}_{i}=\widetilde{L}_{i-1}$
\end{itemize}
Intuitively, $\widetilde{L}_i$ records the number of third-class particles that are to the right of the second-class particle after all vertices below row $t_i$ and all vertices in row $t_i$ to the left of $(x_i, t_i)$ have been sampled.}

\red{
Note that with this definition, $\widetilde{L}_W=L_t$.
Indeed following the trajectory of the second-class particle, the vertices among the $(x_i,t_i)_{1\leq i\leq W}$ visited by the second-class particle are exactly the vertices where the first two options above occur, i.e. $\widetilde{L}_{i-1}$ is either $k_i$ or $k_{i+1}$. $\widetilde{L}$ increases or decreases whenever, at such a vertex, the second class particle goes straight up or straight to the right respectively.
Finally one can see that the number of third class particles to the right of the second class particle at time $t$ is exactly the number of vertices where it has encountered a third-class particle coming from the left and crossed it, minus the number of times it has encountered a third-class particle coming from below and crossed it.
}

\red{
Now let us see how the law of $\widetilde{L}_i$ evolves in $i$.
Recall that we are conditioning on $\bm{Z}$ and therefore consider $(x_i,t_i)_{1\leq i\leq W}$ as fixed.
We will identify laws on $\{0,1,\dots,M\}$ with vectors in $\mathbb R^{M+1}$ and write $(e_i)_{0\leq i\leq M}$ for the standard coordinate basis of this space.
The law of $\widetilde{L}_0$ is given by the coordinate vector $e_0$ since $L_0$ is deterministically $0$.
By the above description, the law of $\widetilde{L}_i$ is obtained from the vector corresponding to the law of $\widetilde{L}_{i-1}$ through the following matrix:
}
\begin{equation}
        P_{k_i}=
        \begin{pmatrix}
            1 & &&&& 0\\
            & \ddots &&&& \\
            && 1-b_1 & b_2 &&\\
            && b_1 & 1- b_2 &&\\
            &&&&\ddots&\\
            0& &&&& 1
        \end{pmatrix}.
    \end{equation}
\red{
Here the $2\times2$ block is in the rows and columns corresponding to the basis vectors $e_{k_i}$ and $e_{k_i+1}$, which are the $(k_i+1)$th and $(k_i+2)$th rows and columns, since we index starting at $0$.
Therefore the law of $L_t=\widetilde{L}_W$ is given by
\begin{equation}\label{eq:Pprod}
        \left(P_{k_W}\dots P_{k_2}P_{k_1}\right) e_0\,.
\end{equation}
}

       \red{Given this description of $L_t$, we will now show that $L_t$ is stochastically dominated by the law of a geometric random variable.  To do so, we introduce a new basis} $\left(\mathbf{v}_i\right)_{i =0}^{M}$.
    Let $\mathbf{v}_i$ be the vector corresponding to the law of the random variable $\min(i,G)$ where $G\sim\text{Geo}(q)$, i.e. 
    \begin{equation}
    (\mathbf{v}_i)_k:= \Pb \left[\min(i,G) = k\right] = 
        \begin{cases}
            (1-q)q^{k} &\text{if } k<i\\
            q^{i}&\text{if } k=i \\
            0,&\text{if $k>i$}.
        \end{cases}
    \end{equation}
    
    The key observation is that this basis satisfies the following relation with the matrices $P_k$ for all $0\leq k\leq M-1$ and $0\leq j\leq M$:
    \begin{equation}\label{eq:PvRelation}
    P_k\mathbf{v}_j=
    \begin{cases}
    \mathbf{v}_j&\text{if $j\neq k,k+1$}\\
    (1-b_2)\mathbf{v}_k+b_2\mathbf{v}_{k+1}&\text{if $j=k$}\\
    b_1\mathbf{v}_k+(1-b_1)\mathbf{v}_{k+1}&\text{if $j=k+1$}\,.
    \end{cases}
\end{equation}
\red{This means that $P_k\mathbf{v}_j$ is a \emph{convex} combination of basis elements $\mathbf{v}_k$ which as we will see below implies that $P_k\mathbf{v}_j$ can be interpreted as the minimum of a geometric random variable with some other independent random variable.}

To see \eqref{eq:PvRelation}, first recall that $q = \frac{b_1}{b_2}$ so that $qb_2 = b_1$. We now check each of the three cases in \eqref{eq:PvRelation}:

\emph{1. $j \neq k, k+1$:}
Since $P_k$ is equal to the identity matrix in all rows except $k$ and $k+1$ $(P_k\mathbf{v}_j)_i=(\mathbf{v}_j)_i$ for $i\neq k,k+1$.
For $j < k$, we have $(\mathbf{v}_j)_k=(\mathbf{v}_j)_{k+1}=0$ and therefore also $(P_k\mathbf{v}_j)_i=(\mathbf{v}_j)_i$ for $i =k,k+1$. For $j > k+1$, we have $(\mathbf{v}_j)_{k+1}=q(\mathbf{v}_j)_k$ and therefore:
\begin{align}
    (P_kv_j)_k=(1-b_1)(\mathbf{v}_j)_k+b_2(\mathbf{v}_j)_{k+1}=(1-b_1+qb_2)(\mathbf{v}_j)_k=(\mathbf{v}_j)_k
\end{align}
and 
\begin{align}
    (P_k\mathbf{v}_j)_{k+1}=b_1(\mathbf{v}_j)_k+(1-b_2)(\mathbf{v}_j)_{k+1}=\left(\frac{b_1}{q}+1-b_2\right)(\mathbf{v}_j)_{k+1}=(\mathbf{v}_j)_{k+1}.
\end{align}

\emph{2. $j=k$:}
We have $(\mathbf{v}_j)_k=q^{k}$ and $(\mathbf{v}_j)_{k+1}=0$. Therefore,
\begin{align}
    \left(\begin{matrix}
  (P_k\mathbf{v}_j)_k\\
  (P_k\mathbf{v}_j)_{k+1}
\end{matrix}\right)=&
    \left(\begin{matrix}
  1-b_1 & b_2\\
  b_1 & 1-b_2
\end{matrix}\right)
\left(\begin{matrix}
  q^{k}\\
  0
\end{matrix}\right)\\=&
\frac{b_1}{q}\left(\begin{matrix}
  (1-q)q^{k}\\
  q^{k+1}
\end{matrix}\right)+
\left(1-b_1-\frac{b_1(1-q)}{q}\right)\left(\begin{matrix}
  q^{k}\\
  0
\end{matrix}\right)\\=&
b_2\left(\begin{matrix}
  (1-q)q^{k}\\
  q^{k+1}
\end{matrix}\right)+
(1-b_2)\left(\begin{matrix}
  q^{k}\\
  0
\end{matrix}\right).
\end{align}
This suffices since $(\mathbf{v}_k)_j=(\mathbf{v}_{k+1})_{j}$ for $j$ different from $k,k+1$. \\

\emph{3. $j=k+1$:} The calculation is similar to the one above and we omit the details.

Returning to the law of $L_t$, by \eqref{eq:PvRelation}, we can write it as a convex combination
\begin{equation} \sum_{i=0}^M\lambda_t(i)\mathbf{v}_i
\end{equation}
of the $\mathbf{v}_i$ for some coefficients $\lambda_t(i)$. Therefore, the vector $\lambda_t$ is the law of a random variable on $\{0,1,\dots,N\}$.
Letting $X$ be a random variable with this law independent of $G\sim\text{Geo}(q)$, we see that $L_t$ is equal in distribution to $\min(X, G)$, and therefore is dominated stochastically by $G$.
This proves the statement.

\red{Note that $M$ being finite allowed us to consider only finitely many vertices in the proof; however, this restriction can easily be removed by noting that we only really need to consider vertices $(x_i, t_i)$ where $x_i \geq \bm{Z}_t(0)$, and there are only finitely many such vertices, even for $M = \infty$.} 
\end{proof}

\begin{remark}
    The proof shows that for fixed t the law of $L_t$ is equal in distribution to the law of the minimum between a geometric random variable and a process $X_t$, which behaves in the same way as $L_t$, except that $b_1$ and $b_2$ are exchanged.
    This seems to be some kind of duality statement.
    It would be interesting to see if this is a specific case of some more general duality.
\end{remark}

We can also obtain a dual statement to Theorem~\ref{thm:orderinginequality}:
\begin{corollary}\label{cor:dualordering}
Let $(\eta_t)_t$ be a multi-class stochastic six-vertex process with the following initial conditions:
    \begin{itemize}
        \item There are some first-class particles (finitely or infinitely many).
        \item There are $M$ second-class particles.
        \item There is a single third-class particle, to the left of all second-class particles.
    \end{itemize}
    Let $L_t$ be the number of second-class particles to the left of the third-class particle.
    Then conditioned on the paths of the first-class particles, and the joint paths of the second- and third-class particles, for any $t\geq0$ the random variable $L_t$ is stochastically dominated by $\text{Geo}(q)$.    
\end{corollary}
\begin{proof}
    In the initial configuration, there are four classes of particles: $\{1,2,3,\infty\}$, (recall that holes are considered particles of class $\infty$).
    We invert the order of classes so that particles of class $1$ become holes, holes become particles of class $1$, and the second and third-class particles swap class.
    Doing this and swapping the $x$ and $t$ coordinates,
    we obtain a stochastic six-vertex model with the same parameters $b_1$ and $b_2$ by Propositions \ref{prop:particle-hole} and \ref{prop:space}.
    This is now a stochastic six-vertex process on the domain $\Z_{\geq 0}\times\Z$, i.e. the right half plane, which can be defined in the same way as the stochastic six-vertex model on the line.
    The boundary conditions obtained after these transformations satisfy the hypothesis of Theorem~\ref{thm:orderinginequality}, with the third-class particles being above the second-class particle.
    The proof then goes through without any changes.
\end{proof}

\begin{remark}
    While we stated Theorem~\ref{thm:orderinginequality} and Corollary \ref{cor:dualordering} for the stochastic six-vertex process on the line, they can also be stated for the stochastic six-vertex model on domains whose boundary is a down-right path. Since the proof takes a vertex-by-vertex approach, it will carry through with minimal changes.
\end{remark}

\section{From linear trajectories to the proof of the main theorem}\label{sec:linear}
We can now begin to prove the main theorem. Let us recall the setup of Theorem~\ref{thm:main}.
We start a stochastic six-vertex process from step initial conditions with a second-class particle at the origin. In this section, and in Sections \ref{sec:hydroToLinear} and \ref{sec:hydroevents}, we always view the stochastic six-vertex model as a particle system on the line (see Definition \ref{def:s6vLine}).
On the line, the above initial conditions are given by a first-class particle at every position $x<0$, a second-class particle at position $0$, and holes at positions $x>0$.
Denote by $(\A_t)_{t\geq0}$ the single-class stochastic six-vertex process given by the first-class particles in this process and by $\bm{X}_t$ the position of the second-class particle at time $t$.
These processes are started from the initial conditions
\[
\A_0(x)=\bm{1}_{x<0}\quad\text{and}\quad X_0=0\,,
\]
and $(\A,\bm{X})$ contains all the information of the multi-class process.
Let $\mathcal F_s$ denote the $\sigma$-algebra generated by $(\A,\bm{X})$ until time $s$, for $s\in\mathbb Z_{\geq 0}$.

Let us now define some events, which will be vital to the proof of the main theorem.

\begin{definition}\label{def:events}
    Fix positive integers $S$ and $T$.
    We define the following $\mathcal F_S$-measurable event, which depends on some $\varepsilon > 0$:
    \begin{equation}\label{eq:bulkSlope}
        P_S =\red{P_S^\varepsilon =} \left\{\frac{\bm{X}_S}{S}\in[\kappa^{-1} + \varepsilon,\kappa -\varepsilon]\right\}.
    \end{equation}
    We also define the following $\mathcal F_{S+T}$-measurable events, which depend on some $\gamma \in [0,1]$:
    \begin{align}
        E_S^\geq &=\left \{\bm{X}_{S+T}-\bm{X}_S\geq \frac{\bm{X}_S}{S}T-S^{1-\gamma}\right\} \label{eq:defE}\\
        E_S^\leq &= \left\{\bm{X}_{S+T}-\bm{X}_S\leq \frac{\bm{X}_S}{S}T+S^{1-\gamma}\right\} \nonumber.
    \end{align}
    Finally, we let $E_S:=E_S^\geq\cap E_S^\leq$.
\end{definition}

On the event $P_S$, the speed $\frac{\bm{X}_S}{S}$ is bounded strictly away from the edge of the rarefaction fan so that the effective hydrodynamic bounds in Corollary \ref{cor:effective_hydrodynamics} will apply. On event $E_S$, we control how much the speed of the second-class particle at time $S+T$ deviates from the speed at time $S$. In the following proposition, we choose appropriate time steps $S$ and $T$ and show the existence of a high-probability $\mathcal{F}_S$-measurable \emph{hydrodynamic event} $H_S$ upon which $E_S$ will hold with high probability at time $S+T$.
The precise definition of this event will be given in Proposition~\ref{prop:hydroEvent}.
We call $H_s$ the hydrodynamic event since it is the event upon which at time $S$ the height function of $\A$ has not deviated too much from the hydrodynamic limit. We will show that the same is true at time $S+T$ with high probability, for any possible configuration in $H_S$.

\begin{proposition}\label{prop:straightlinespoly}
    For any integer $S \geq 1$, let $T=S^\beta$ for some $\beta\in(\frac23,1)$.
    For any positive $\alpha<\frac\beta2-1/3$ and for any $\varepsilon\in(0,\frac14)$, there is a $c=c(\varepsilon,\alpha)>0$ and an $\mathcal{F}_S$-measurable event $H_S=\red{H_S^\varepsilon}$ such that for all $S \geq 1$  and for $\gamma=\frac56-\frac\beta2-\alpha$ we have
    \begin{equation}
        \mathbb P[H_S]\geq 1- c^{-1}e^{-cS^{\alpha}}\quad\mathbb P[E_S|\mathcal F_S]\geq (1-c^{-1}e^{-cT^\alpha})\bm{1}_{P_s\cap H_s}.
    \end{equation}
\end{proposition}

Note that since $T < S$, the bound on $\Pb [H_S]$ remains true if we replace $S$ with $T$.  
\red{We suppress the dependence of $P_S^\varepsilon$ and $H_S^\varepsilon$ in most cases, since many of the following considerations are for fixed $\varepsilon$.}

\begin{remark}
The parameters $\alpha,\beta$, and $\gamma$ have the following meaning.
The parameter $\beta$ controls the size of the time steps, with $\beta$ closer to $1$ giving larger time steps.
The parameter $\gamma$ determines our control on the trajectory of the particle, with larger $\gamma$ giving tighter bounds.
For each $\beta$ we can prove the statement of the proposition for each $\gamma<\frac56-\frac{\beta}{2}$ and the difference between the two is $\alpha$.
Therefore, small $\alpha$ are of interest.
There is a trade-off for the value of $\beta$.
Bigger $\beta$ gives better control along the sequence of time steps.
To see this note that $S^{1-\gamma} =T^{\frac16\beta^{-1}+\frac12+\alpha\beta^{-1}}$, and this exponent is minimized for $\beta=1$. 
However, between the time steps, bigger $\beta$ leads to less control, since we use a rough bound based on monotonicity between time steps.
We will ultimately set $\beta=\frac79$, where the value $\beta=\frac79$ is obtained by an optimization balancing these two effects.
\end{remark}

Before proving this proposition, let us see how it implies Theorems \ref{thm:main} and \ref{thm:fluct}.
This largely follows \cite[Section 4]{ACG2023asepspeed}, especially Proposition~\ref{prop:intersectionprobability}, except that there a different choice of time sequence was made, see Remark~\ref{rmk:sn}.

\begin{definition}\label{def:timesteps}
    For $S_0\geq 2$ define $(S_m)_{m\geq1}$ and $(T_m)_{m\geq0}$ as follows. Let $T(S) = S^{\beta}$. 
    Let $T_m=T(S_m)$ and $S_{m+1}=S_m+T_m$.
    Note that $T(S)$ is strictly increasing for all $S\geq 1$.
\end{definition}

\begin{lemma}\label{lem:boundSm}
   There exists constants $z_-=z_-(\beta)>0$ and $z_+=z_+(\beta)>0$, such that $z_-m^{\frac{1}{1-\beta}} \leq S_m\leq z_+m^{\frac{1}{1-\beta}} $ for all $m\geq1$.
\end{lemma}

\begin{proof}
   This follows from a straightforward induction argument, using the fact that $z_-x^\frac{1}{1-\beta}+(z_-x^\frac{1}{1-\beta})^\beta-z_-(x+1)^\frac{1}{1-\beta}$ is both positive at $x=2$ and strictly increasing on $[2,\infty]$ for $z_-$ small enough (but still positive).
    The upper bound also follows by induction, for $z^+$ bigger than $1$ and such that the inequality is true for $S_1$.
\end{proof}

\begin{lemma}\label{lem:bounds}
    For each $\varepsilon$ there exists a $D=D(\varepsilon,\beta,\alpha)$ such that for $S_0\geq D$ the following holds
\begin{equation}
    \sum_{m\geq0}S_m^{-\gamma}\leq\varepsilon/9;
    \quad\sum_{m\geq0}c^{-1}e^{-cT_m^\alpha}<\varepsilon/2;
    \quad\kappa\frac{T_m}{S_m}<\varepsilon/9\textrm{,~for~all~} m\geq0\,.
\end{equation}
\end{lemma}

\begin{proof}
    By Lemma \ref{lem:boundSm} and the monotonicity mentioned in Definition~\ref{def:timesteps}, $S_m\geq \max(S_0,z_-S^\frac{1}{1-\beta})$.
    The term $(z_-m^\frac{1}{1-\beta})^{-\gamma}=z_-^{-\gamma}m^{\frac{-\tfrac56+\tfrac\beta2+\alpha}{1-\beta}}$ is summable since the exponent is less than $-1$ by the assumption $\alpha<\frac{\beta}{2}-\frac13$.
    Therefore by the dominated convergence theorem as $D\to\infty$ this sum goes to $0$.
    The second and third sums are convergent since $S_m$ and $T_m$ grow faster than linear in $m$, and therefore the sums can be bounded by geometric series.
    By the same dominated convergence argument, the statement follows.    
    Finally $\kappa\frac{T_m}{S_m}<\varepsilon/4$ simply follows from the fact that $\frac{T_m}{S_m}$ goes to $0$ as $S_0$ goes to $\infty$. 
\end{proof}

\begin{definition}
    Define $\varepsilon_m$ by setting $\varepsilon_0=\varepsilon\in(0,\frac14)$ and
    \begin{equation*}
        \varepsilon_{m+1}=\varepsilon_m-S_m^{-\gamma}
    \end{equation*}
    and note that for $S_0>D(\varepsilon,\beta,\alpha)$ all $\varepsilon_m$ are positive (in fact they are greater than $\frac{8\varepsilon}{9}$) by Lemma \ref{lem:bounds}.
    Define further the event
    \begin{equation}
        L^\varepsilon_{S_0}(k)=
        \bigcap_{m=0}^{k-1}P^{\varepsilon_m}_{S_m}\cap H^{\varepsilon_m}_{S_m} \cap E_{S_m}
    \end{equation}
    and let $L^\varepsilon_{S_0}=L^\varepsilon_{S_0}(\infty)$.
\end{definition}

Since the constant $c$ in Proposition~\ref{prop:straightlinespoly} can be taken to be weakly decreasing in $\varepsilon$, we can assume that the statement of Proposition~\ref{prop:straightlinespoly} holds with the same $c$ for all $\varepsilon_m$ in the definition above.

\begin{proposition}\label{prop:intersectionprobability}
    There exists a constant $d$ such that for all $\varepsilon\in(0,1/4)$ there is a constant $D$ such that for all $S_0\geq D$ the probability of  $L^\varepsilon_{S_0}$ is at least $1-(d+1)\varepsilon$.
\end{proposition}
\begin{proof}
    Note first that the event $P^{\varepsilon_m}_{S_m}\cap E_{S_m}$ is contained in the event $P^{\varepsilon_{m+1}}_{S_{m+1}}$, since on the event $E_{S_m}$ it holds that $\left|\frac{\bm{X}_{S_{m+1}}}{S_{m+1}}-\frac{\bm{X}_{S_m}}{S_m}\right |\leq S_m^{-\gamma}$.
    Using this we obtain
    \begin{align}\label{eq:probL}
        \mathbb P[L^\varepsilon_{S_0}(k)]&=
        \mathbb P[P^{\varepsilon_0}_{S_0}]
        -\sum_{m=0}^{k-1}\mathbb P[L_{S_0}^{\varepsilon_0}(m)\cap P_{S_m}^{\varepsilon_m}\cap(H_{S_m}^{\varepsilon_m})^c]\\
        &-\sum_{m=0}^{k-1}\mathbb P[L_{S_0}^{\varepsilon_0}(m)\cap P_{S_m}^{\varepsilon_m}\cap H_{S_m}^{\varepsilon_m}\cap(E_{S_m}^{\varepsilon_m})^c]
        -\sum_{m=1}^{k-1}\mathbb P[L_{S_0}^{\varepsilon_0}(m)\cap(P_{S_m}^{\varepsilon_m})^c]\\
        &\geq \mathbb P[P^{\varepsilon_0}_{S_0}]-\sum_{m\geq0}\mathbb P[P^{\varepsilon_m}_{S_m}\cap (H^{\varepsilon_m}_{S_m})^c]-\mathbb P[P^{\varepsilon_m}_{S_m}\cap H^{\varepsilon_m}_{S_m}\cap (E_{S_m})^c]\,,
    \end{align}
    where we used the fact noted above to observe that $\mathbb P[L_{S_0}^{\varepsilon_0}(m)\cap(P_{S_m}^{\varepsilon_m})^c]=0$.
    By Proposition \ref{prop:weak}, $\frac{\bm{X}_t}{t}$ converges in law to a continuous random variable on $[\kappa^{-1},\kappa]$,  with density $\frac{\sqrt{\kappa}}{2(\kappa-1)}x^{-\frac32}$.
     This density is bounded, and therefore there is a constant $d$ such that for all $S_0$ large enough $\mathbb P[P^{\varepsilon_0}_{S_0}]\geq 1-d\varepsilon$, for all $\varepsilon_0$.
    Note further that by Proposition \ref{prop:straightlinespoly}, both $\mathbb P[P^{\varepsilon_m}_{S_m}\cap (H^{\varepsilon_m}_{S_m})^c]$ and $\mathbb P[P^{\varepsilon_m}_{S_m}\cap H^{\varepsilon_m}_{S_m}\cap (E_{S_m})^c])$ are less than $c^{-1}e^{-cT_m^\alpha}$, and thus by Lemma~\ref{lem:bounds} the right-hand side of \eqref{eq:probL} is at least $1-(d+1)\varepsilon$.
\end{proof}
\begin{proof}[Proof of Theorem \ref{thm:main}]
    We already have weak convergence to the desired distribution from Proposition \ref{prop:weak}, so it remains to prove a.s. convergence. We will show that with probability at least $1-(d+1)\varepsilon$ the difference between the limit superior and the limit inferior of $\frac{\bm{X}_t}{t}$ is less than $\varepsilon$.
    This immediately implies that this probability is indeed $1$ since these events form a decreasing family as $\varepsilon$ goes to $0$.
    Since this holds for any $\varepsilon$, the conclusion follows as the limit superior and the limit inferior must then be equal with probability $1$.

   \textbf{Claim:} Fix $\varepsilon > 0$. There exists $D$ such that for $s, s' > D$, with probability at least $1-(d+1)\varepsilon$:
    \begin{equation}
        \left|\frac{\bm{X}_{s}}{s}-\frac{\bm{X}_{s'}}{s'}\right|\leq \varepsilon\,.
    \end{equation}
    
    \textbf{Proof of Claim:} Let $D$ be large enough such that both Lemma~\ref{lem:bounds} and Proposition~\ref{prop:intersectionprobability} hold.
    Let $S_0=D$.
    For $S_0<s<s'$ let $m$ be the largest integer such that $S_m<s$ and $m'$ be the smallest integer such that $s'<S_{m'}$.
    By Proposition~\ref{prop:intersectionprobability} the event $L^\varepsilon_{S_0}$ holds with probability at least $1-(d+1)\varepsilon$.
    Assume now that the event $L^\varepsilon_{S_0}$ takes place, so that in particular $E_{S_m}$ takes place for all $m\geq0$.
    Then
    \begin{align}
        \left|\frac{\bm{X}_{s}}{s}-\frac{\bm{X}_{S_m}}{S_m}\right|
        &\leq \left|\frac{\bm{X}_{s}}{s}-\frac{\bm{X}_{S_{m+1}}}{s}\right|+\left|\frac{\bm{X}_{S_{m+1}}}{s}-\frac{\bm{X}_{S_{m+1}}}{S_{m+1}}\right|+\left|\frac{\bm{X}_{S_{m+1}}}{S_{m+1}}-\frac{\bm{X}_{S_m}}{S_m}\right|\\
        &\leq \frac{|\bm{X}_{S_m}-\bm{X}_{S_{m+1}}|}{S_m}+\frac{|\bm{X}_{S_{m+1}}|}{S_{m+1}}\frac{S_{m+1}-S_m}{S_m}+S_m^{-\gamma}\\
        &\leq 2\left(\kappa \frac{T_m}{S_m}+S_m^{-\gamma}\right)
 \leq \tfrac49\varepsilon\,,
    \end{align}
    where in the second inequality we used monotonicity of $\bm{X}_s$ for the first two terms, and for the third term we used that on $E_{S_m}$ it holds that $|\frac{\bm{X}_{S_{m+1}}}{S_{m+1}}-\frac{\bm{X}_{S_m}}{S_m}|\leq S_m^{-\gamma}$.
    For the third inequality we used the event $E_{S_m}$ to control $|\bm{X}_{S_m}-\bm{X}_{S_{m+1}}|$ and that $\frac{\bm{X}_{S_m}}{S_m}<\kappa$, since $P_{S_m}^{\varepsilon_m}$ holds.
    The final inequality then follows from Lemma~\ref{lem:bounds}.
    By the same argument $\left|\frac{\bm{X}_{s'}}{s'}-\frac{\bm{X}_{S_{m'}}}{S_{m'}}\right|\leq \frac{4}{9}\varepsilon$.
    Finally, it holds that
    \begin{equation}
        \left|\frac{\bm{X}_{S_m}}{S_m}-\frac{\bm{X}_{S_{m'}}}{S_{m'}}\right|\leq\sum_{n=m}^{m'-1} \left|\frac{\bm{X}_{S_{n}}}{S_{n}}-\frac{\bm{X}_{S_{n+1}}}{S_{n+1}}\right|\leq\sum_{n=m}^{m'-1}S_n^{-\gamma}\leq \tfrac19\varepsilon
    \end{equation}
    since $E_{S_n}$ holds for all $m\leq n\leq m'$ and Lemma~\ref{lem:bounds}.
    It follows that with probability at least $1-(d+1)\varepsilon$,
    \begin{equation}
        \left|\frac{\bm{X}_{s}}{s}-\frac{\bm{X}_{s'}}{s'}\right|\leq \varepsilon\,,
    \end{equation}
    for all $s,s'> S_0$ which was the claim.

    This now implies that with probability at least $1-(d+1)\varepsilon$, 
    \begin{equation}
       \left| \limsup_{t \to \infty} \frac{\bm{X}_t}{t} - \liminf_{t \to \infty} \frac{\bm{X}_t}{t}\right| < \varepsilon,
    \end{equation}
    and the proof of Theorem~\ref{thm:main} follows.
    
\end{proof}

Note that the conditions on $\beta$ and $\alpha$ together with the definition of $\gamma$ imply
\begin{equation}\label{eq:betaplusgamma}
    \beta+\gamma>1
\end{equation}

\begin{proof}[Proof of Theorem \ref{thm:fluct}]
For time steps given by $T(S)=S^\beta$, again for each $\varepsilon$ there exists a $D(\varepsilon,\beta,\alpha)$ such that the event $L^\varepsilon_{S_0}$ occurs with probability at least $1-(d+1)\varepsilon$ for all $S_0>D$. Letting $C$ be a constant that can depend on $\beta$ and $\alpha$ and that can change from line to line, we note that on this event we have  
    \begin{align}
    \left|\bm{X}_{S_n}-S_n\bm{U}\right|=&
    \left|\bm{X}_{S_n}-S_n\left(\frac{\bm{X}_{S_n}}{S_n}+\sum_{m\geq n}\frac{\bm{X}_{S_{m+1}}}{S_{m+1}}-\frac{\bm{X}_{S_m}}{S_m}\right)\right|\\
    =&S_n\left|\sum_{m\geq n}\frac{\bm{X}_{S_{m+1}}}{S_{m+1}}-\frac{\bm{X}_{S_m}}{S_m}\right|\\
    \leq&S_n\sum_{m\geq n}S_m^{-\gamma}\\
    \leq&CS_n\sum_{m\geq n}m^{-\gamma\frac{1}{1-\beta}} \label{eq:pseries}\\
    \leq&CS_nn^{-\gamma\frac{1}{1-\beta}+1}\\
    \leq&CS_n^{2-\gamma-\beta}. \label{eq:Snbound}
\end{align}
In \eqref{eq:pseries} and \eqref{eq:Snbound}, we used Lemma \ref{lem:boundSm} to bound $S_n$ from above and below by constants times $n^\frac{1}{1-\beta}$ and we used the fact that the $-\gamma\frac{1}{1-\beta}\leq -1$ by \eqref{eq:betaplusgamma} to ensure that $p$-series in \eqref{eq:pseries} converges and can be bounded.
It remains to check the behavior times $s$ that are not of the form $S_n$.
For $s$ between $S_n$ and $S_{n+1}$ we get
\begin{align}
    |\bm{X}_s-s\bm{U}|
    &\leq|\bm{X}_s-\bm{X}_{S_n}|+|\bm{X}_{S_n}-S_n\bm{U}|+|S_n\bm{U}-s\bm{U}|\\
    &\leq |\bm{X}_{S_{n+1}}-\bm{X}_{S_n}|+|\bm{X}_{S_n}-S_n\bm{U}|+\kappa|S_{n+1}-S_n|\\
    &\leq C(S_n^{1-\gamma}+S_n^{2-\gamma-\beta}+S_n^\beta)\leq C(s^{2-\gamma-\beta}+s^\beta) \label{eq:intertimes},,
\end{align}
where in the second inequality, we used monotonicity of $\bm{X}_s$ in the first term and that $U$ is bounded by $\kappa$ in the third term. For the third inequality, we used that the event $E_n$ holds for the first term, \eqref{eq:Snbound} for the second term, and the definition of $S_{n+1}$ for the third term. Setting $\beta=\frac79$ gives $\gamma=\frac49-\alpha$, which gives $2-\gamma-\beta=\frac{7}{9}+\alpha$, therefore showing that on an event of probability $1-(d+1)\varepsilon$ the limit of $|\bm{X}_S-s\bm{U}|s^{-\frac79-2\alpha}$ is $0$, which concludes the proof.
\end{proof}

\begin{remark} \label{rmk:sn}
    While we used the sequence $S_{n+1}= S_n + S_n^{\beta}$, which grows polynomially, in \cite{ACG2023asepspeed}, they used the sequence $S_{n+1}= S_n+\frac{S_n}{\log S_n},$ which grows like $e^{\sqrt{n}}$. The latter sequence would have sufficed to prove Theorem~\ref{thm:main}.
    Furthermore, this time scale gives
    \[
    |X_{S_n}-S_nU|\leq S_n^{\frac23+}
    \]
    which is the expected order of fluctuations.
    However, it grows too quickly to prove the finer statement in Theorem~\ref{thm:fluct} as in \eqref{eq:intertimes}, we crucially used the fact that $|S_{n+1} - S_n| \leq s^{\beta}$ to bound the fluctuations between the times $S_n$. 
\end{remark}

\section{From hydrodynamic events to linear trajectories}\label{sec:hydroToLinear}

The purpose of this section is to prove Proposition \ref{prop:straightlinespoly}.
We write the full proof only for $E^\geq_S$, since $E^\leq_S$ can be treated similarly.
We will outline the proof for $E^\leq$ at the end of Section~\ref{sec:hydroevents}.

Let us now proceed with the proof for the process $E^\geq$.
We couple the process $(\A,\bm{X})$ to a new multi-class stochastic six-vertex process $\B$ by filling in every position to the left of $\bm{X}_S$ with particles.
Then Theorem~\ref{thm:orderinginequality} allows us to control the position of $\bm{X}_{S+T}$ by controlling a large number of these additional particles.

\begin{definition}
    Define a new multi-class process $(\B_t)_{t\geq0}$, with first-, second- and third-class particles, depending on $\mathcal A_S$ in the following way:
    \begin{itemize}
        \item It has the same parameters $b_1$ and $b_2$.
        \item At time $0$ each site $j\in\mathbb Z$ in $\B_0$ is occupied by a first-class particle if it is occupied by a particle in $\A_S$,
        \item the site $\bm{X}_S$ is occupied by a second-class particle,
        \item and each site to the left of $\bm{X}_S$ not occupied by a particle in $\A_S$ is also occupied by a third-class particle.
    \end{itemize}
    Further let $\bm{M}$ be the number of third-class particles in $\B$, which is finite, since it is at most $\bm{X}_S$.
    Further, let $\bm{Z}_t(0)>\bm{Z}_t(1)>\dots>\bm{Z}_t(\bm{M})$ denote the ordered positions of the second-class particle and the third-class particles in $\B_t$.
    At time $0$ we have $\bm{Z}_0(0)=\bm{X}_S$, but at later times, $\bm{X}_{S+t}$ can be any of the positions $\bm{Z}_t$.
\end{definition}

We will also consider the two single-class processes obtained by merging the second-class and third-class particles in $\B$ with either the holes or the first-class particles.
That is, let $(\B^{(1)}_t)_{t\geq0}$ be the single-class process given by just the first-class particles in $\B$ and let $(\B^{(1,2,3)}_t)_{t\geq0}$ be the single-class process of the first-, second- and third-class particles in $\B$ forgetting their classes. 
The triplet $(\B^{(1)}_t,\B^{(1,2,3)}_t,\bm{X}_{S+t})_{t\geq0}$ contains all the information of $\B$.

The process $\B$ depends both on $\A_S$ and $\bm{X}_S$.
We will often use a union bound over possible values of $\bm{X}_S$.
To do so it will be convenient to introduce a version of $\B$ in which the position of the second-class particle is replaced by a deterministic position $X$.
For $X$ in the interval $[S(\kappa^{-1}+\varepsilon),S(\kappa-\varepsilon)]\cap\mathbb Z$, define $\B^X$ as the process obtained from $\A_S$ by adding a second-class particle at $X$ if that position is empty and filling all empty positions $k<X$ with third-class particles (note that $\B^X$ may not have a second-class particle if there is a particle at position $X$ in $\A_S$).
Let further $\B^{(1,2,3),X}$ be the single-class process of the first- and second-class particles in $\B^X$.
Note that it is possible to couple all the processes $\B^X$, $\B$, $\A$ and $\bm{X}$ such that the first-class particles in $\B^X_t$ and $\B_t$ are given by $\A_{S+t}$, and such that substituting $\bm{X_S}$ for $X$ it holds that
\begin{equation}\label{eq:coupling}
\B^{\bm{X}_S}=\B\,.
\end{equation}
Note however that the law of $\B$ conditioned on $\bm{X}_S=X$ is not given by $\B^X$, since $\bm{X}_S$ and the first-class particles in $\B$ are non-trivially correlated.
However, by showing that certain events happen for all $\B^X$ with exponentially small probability, one can use \eqref{eq:coupling} to show that they also happen for $\B$ with exponentially small probability.

\begin{lemma}\label{lem:domination}
Let $\bm{G}\sim\mathrm{Geo}(q)$ be independent of $\mathcal F_S$ and \red{$(\bm{Z}_t)_{0\leq t\leq T}$}.
Then for any $y\in\Z$ and any $S, T \geq 1$ it holds that
\begin{equation}\label{eq:stochdom}
    \mathbb{P}[\bm{X}_{S+T}\geq \bm{X}_S+y|\mathcal F_S]\geq \mathbb{P}[\bm{Z}_T(\bm{G}\wedge \bm{M})\geq \bm{X}_S+y|\mathcal F_S]\,.
\end{equation}
\end{lemma}
\begin{proof}
    \red{Recall that $\mathcal F_S$ is the sigma-algebra generated by both the process $\A$ and the position of the second-class particle $\bm{X}$ until time $S$}.
    Denote the number of third-class particles that are to the right of the second-class particle in $\B_T$ by $\bm{K}$, such that $\bm{X}_{S+T}=\bm{Z}_T(\bm{K})$.
    By Theorem~\ref{thm:orderinginequality}, the law of $\bm{K}$ conditioned on $\mathcal{F}_S$ and $(\bm{Z}_t)_{0\leq t\leq T}$ is dominated by $\textrm{Geo}(q)$.
    Therfore $\bm{K}$ can be coupled to a random variable $\bm{G}^*\sim\textrm{Geo}(q)$ independent of $\mathcal F_S$ and $\bm{Z}_T$ such that $\bm{K}\leq \bm{G}^*$ almost surely.
    Thus we obtain
    \begin{align}
    \mathbb{P}[\bm{X}_{S+T}\geq \bm{X}_S+y|\mathcal F_S]
    &=\mathbb{P}[\bm{Z}_T(\bm{K})\geq \bm{X}_S+y|\mathcal F_S]\\
    &\geq \mathbb{P}[\bm{Z}_T(\bm{G}^*\wedge \bm{M})\geq \bm{X}_S+y|\mathcal F_S],\label{eq:firstineq}
    \end{align}
    where in \eqref{eq:firstineq} we used $\bm{K}\leq\bm{G}^*\wedge \bm{M}$ and the ordering of $\bm{Z}_T$.
    Note that the right-hand side does not depend on the coupling between $\bm{K}$ and $\bm{G}^*$ since $\bm{G}^*$ is independent of $\mathcal F_S$ and $\bm{Z}$. Therefore, we can replace $\bm{G}^*$ by $\bm{G}$ in \eqref{eq:firstineq}.
\end{proof}
Let $\bm{L}$ be defined as
\begin{equation}\label{eq:defL}
    \bm{L}=\#\left\{\textrm{second- and third-class class particles in }\B_T\textrm{ to the right of }\frac{\bm{X}_S}{S}(S+T)-S^{1-\gamma}\right\}.
\end{equation}
Using Lemma~\ref{lem:domination}, we can reduce the proof of Proposition \ref{prop:straightlinespoly} to the following lemma which states that $\bm{L}$ is of order at least $S^\frac13$ with high probability.

\begin{lemma}\label{lem:particlesStayClose}
    For any positive $\varepsilon<\frac14$, and for $T,\alpha$ and $\gamma$ as in Proposition \ref{prop:straightlinespoly}, there is a constant $c>0$ and an $\mathcal{F}_S$-measurable event $H_S$ such that for all $S \geq 1$
    \begin{equation}
        \mathbb{P}(H_S)\geq 1-c^{-1}e^{-cS^{\alpha}}
    \end{equation}
    and
    \begin{equation}\label{eq:particlesStayClose}
        \mathbb P[\bm{L}\geq S^{\frac13}|\mathcal{F}_S]\geq (1-c^{-1}e^{-cT^\alpha})\bm{1}_{P_S\cap H_S}\,.
    \end{equation}
\end{lemma}
Before proving this lemma let us see how it implies Proposition \ref{prop:straightlinespoly}.

\begin{proof}[Proof of Proposition~\ref{prop:straightlinespoly}]
    As said above, we only prove the statement for $E^\geq_S$.
    Condition on $\mathcal F_S$ and assume that $H_S\cap P_S$ holds.
    Let $\bm{G}\sim\textrm{Geo}(q)$ be independent of $\bm{Z}$ and $\mathcal F_S$, as above. 
    Define the events
    \begin{equation}
        N_S= \left\{\bm{Z}_T(\bm{G}\wedge \bm{M})\geq\frac{\bm{X}_S}{S}(S+T)+S^{1-\gamma}\right\}\quad\text{and}\quad G_S=\left\{\bm{L}>S^{\frac13}\right\}
    \end{equation}
    and recall the definition of the event $E^\geq_S$ is
    \[
    E_S^\geq =\left \{\bm{X}_{S+T}-\bm{X}_S\geq \frac{\bm{X}_S}{S}T-S^{1-\gamma}\right\}\,.
    \]
    Setting $y =\frac{\bm{X}_S}{S}T-S^{1-\gamma}$, it follows from Lemma~\ref{lem:domination} that $\mathbb P[E^\geq_S|\mathcal F_S]\geq\mathbb P[\red{N}_S|\mathcal F_S]$.
    By the distribution of $\bm{G}$, 
    \[\mathbb P[N_S | G_S, \mathcal F_S ]  \geq \mathbb P[\bm{G} \leq S^{1/3}] \geq  1-q^{S^{\frac13}}.\]
    Combining this with the statement of Lemma \ref{lem:particlesStayClose} gives
    \begin{align*}
        \mathbb P[E^\geq_S|\mathcal F_S]&\geq \mathbb P[N_S | \mathcal F_S] \\
        & \geq \mathbb P[N_S | G_S, \mathcal F_S]\Pb[G_S|\mathcal F_S]\\
        &\geq \left(1-q^{S^{\frac13}}\right)\left(1-c^{-1}e^{-cT^\alpha}\right)\bm{1}_{H_S\cap P_S}\\
        &\geq (1-c^{-1}e^{-cT^\alpha})\bm{1}_{H_S\cap P_S}\,.
    \end{align*}
    The final inequality is obtained by decreasing $c$ so that we can absorb the term $q^{S^{\frac13}}$ (since $\alpha\leq1/3$ and $S\geq T$).
\end{proof}

\section{From effective hydrodynamics to hydrodynamic events}\label{sec:hydroevents}
The purpose of this section is to prove Lemma~\ref{lem:particlesStayClose}.
To this end let $(\B^{\text{step},X})_{t\geq0}$ be a stochastic six-vertex process started from step initial conditions shifted by $X$, i.e. $\B^{\text{step},X}_{0}(x) = \bm{1}_{x \leq X}$.
Clearly, at time $0$ we have that $\B^{(1,2,3),X}_0(x)\geq \B^{\text{step},X}_0(x)$ for all $x\in\Z$, so by attractivity (Proposition \ref{prop:attractivity}) we can couple them so that this holds for any later time as well.
Note that $\B^{(1,2,3),X}$ and $\B^{(1)}$ are already coupled such that $\B^{(1,2,3),X}_t(x)\geq \B^{(1)}_t(x)$ for any time $t$ by their relation to the multi-class process $\B^X$.
This gives a coupling of the three processes $\B^{(1,2,3),X},\B^{(1)},\B^{\text{step},X}$ such that at all times $t$ it holds that 
\[
\B^{(1,2,3),X}_t(x)\geq\max(\B^{(1)}_t(x),\B^{\text{step},X}_t(x)) \text{ for all } x\in\Z.
\] 
The process $\B^{\text{step},X}$ is a stochastic six-vertex process started from step initial conditions, thus we can use Corollary~\ref{cor:effective_hydrodynamics} to control its height function.
Understanding $\B^{(1)}$ is more intricate since its initial conditions are given by $\A_S$. 
We will not be able to control $\B^{(1)}$ for all values of $\A_S$, and instead find an $\mathcal F_S$-measurable event $H_S$, which we call the \emph{hydrodynamic event}, that holds with high probability.
On this event, with high probability, $\B_T^{(1)}$ is close to the hydrodynamic profile of a stochastic six-vertex process started from step initial conditions and evaluated at time $S+T$.

To simplify notation throughout this section, we will define $h_t([X, Y]; \eta):= h_t(X; \eta) - h_t(Y; \eta).$ Let us now define the following shorthand notation for the event that a process is close to its hydrodynamic limit.

\begin{definition}
  For a single-class stochastic six-vertex process $(\A_t)_{t\geq0}$ started from step initial conditions, a time $t$ and $\alpha,\varepsilon>0$ define the event $C^{\alpha,\varepsilon}_t(\A)$ as 
    \begin{equation}\label{eq:hydrostayclose}
        \left\{\left|
    h_{t}([X,Y];\A)
    -t
    (g(\tfrac{X}{t})-g(\tfrac{Y}{t}))\right|
    \leq t^{\frac{1}{3}+\alpha}\text{, for all $X$ and $Y$ in $[t(\kappa^{-1}+\varepsilon),t(\kappa-\varepsilon)]$}\right\}\,.
    \end{equation}
    Note that this event does not depend on the choice of height function $h(x,\A)$ since it only concerns height function differences.
\end{definition}
Using this definition, Corollary~\ref{cor:effective_hydrodynamics} with $s=t^{\alpha}$ gives $\mathbb{P}[C^{\alpha,\varepsilon}_t(\A)]\geq 1-c^{-1}t^2e^{-ct^\alpha}$, with the constant $c>0$ depending on $\varepsilon$.
\begin{proposition}\label{prop:hydroEvent}
    For $T=T(S)$ as in Proposition~\ref{prop:straightlinespoly}, for $\alpha,\varepsilon>0$, and for $S\geq1$, there is a constant $c=c(\varepsilon)$ and an $\mathcal F_S$-measurable event $H_S$ that holds with probability $1-c^{-1}e^{-cS^\alpha}$ such that
\begin{equation}
    \mathbb P
    \left[C^{\alpha,\varepsilon}_{S+T}(\A)|\mathcal F_S\right]\geq \bm{1}_{H_S}(1-c^{-1}e^{-cS^\alpha})\,.
\end{equation}
\end{proposition}
\begin{proof}
    By making $c$ smaller the right-hand side of \eqref{eq:hydrostayclose} can be made non-positive for $S$ small, so it suffices to consider $S$ large enough.
    Let $(\widetilde{\A}_t)_{0\leq t\leq S}$ be an independent copy of $(\A_t)_{0\leq t\leq S}$, i.e. a stochastic six-vertex process started from step initial conditions and run until time $S$.
    After time $S$ we will couple these two processes, such that they are no longer independent.
    
    Define $H_S=C^{\frac{a}{2},\frac{\varepsilon}{2}}_S(\A)$.
    This event has the desired probability by Corollary~\ref{cor:effective_hydrodynamics} and is $\mathcal F_S$-measurable.
    Further, let $\widetilde{H}_S=C^{\frac{a}{2},\frac{\varepsilon}{2}}_S(\widetilde{\A})$.
    On the intersection of $H_S$ and $\widetilde{H}_S$ it holds that
    \begin{equation} \label{eq:hydroH}
       \left|
    h_{S}([X, Y];\A)
    -S
    \left(g\left(\tfrac{X}{S}\right)-g\left(\tfrac{Y}{S}\right)\right)\right|
    \leq S^{\frac{1}{3}+\frac{\alpha}{2}}\text{, for all $X$ and $Y$ in $[S(\kappa^{-1}+\varepsilon/2),S(\kappa-\varepsilon/2)]$}   
    \end{equation}
    
    and
    \begin{equation}\label{eq:hydroHtilde}
        \left|
    h_{S}([X,Y];\widetilde{\A})
    -S
    \left(g\left(\tfrac{X}{S}\right)-g\left(\tfrac{Y}{S}\right)\right)\right|
    \leq S^{\frac{1}{3}+\frac{\alpha}{2}}\text{, for all $X$ and $Y$ in $[S(\kappa^{-1}+\varepsilon/2),S(\kappa-\varepsilon/2)]$}.
    \end{equation}
    These events do not depend on the choice of height function and therefore we can choose the height functions which satisfy $h_S(S(\kappa-\frac\varepsilon2),\A)=h_S(S(\kappa-\frac\varepsilon2),\widetilde{\A})=0$.
    With this choice of height function and setting $Y=S(\kappa-\frac\varepsilon2)$ in \eqref{eq:hydroH} and \eqref{eq:hydroHtilde},   it follows that on the intersection $H_S\cap \widetilde{H}_S$ we have 
    \[
    |h_S(X;\A)-h_S(X;\widetilde{\A})|\leq 2S^{\frac{1}{3}+\frac\alpha2}\text{ for }X\in[S(\kappa^{-1}+\varepsilon/2),S(\kappa-\varepsilon/2)]\,.
    \]
    Thus we can couple $(\widetilde{\A}_t)_{t\geq S}$ with $(\A_t)_{t\geq S}$ via the coupling given in Lemma~\ref{lem:finitiemonotonicity} with $M=S^{\frac13}$ and $K=2S^{\frac{1}{3}+\frac\alpha2}$.
    Define the event that this coupling succeeds as
    \begin{equation}
        D=\left\{\left|h_{S+T}(X;\widetilde{\A})-h_{S+T}(X;\A)\right|<6S^{\frac13+\frac\alpha2}+S^{\frac13}\text{ for } X\in[S(\kappa^{-1}+\tfrac\varepsilon2)+\tfrac{2T}{1-b_2}+1,S(\kappa-\tfrac\varepsilon2)]\right\}
    \end{equation}
    Then the statement of Lemma~\ref{lem:finitiemonotonicity} gives
    \[
    \mathbb P[D|H_S\cap \widetilde{H}_S,\mathcal F_S]\geq 1-c^{-1}e^{-cT}-c^{-1}e^{-cS^{\frac{1}{3}}}\,,
    \]
    which implies that
    \begin{equation}\label{eq:couplingprobability}
    \mathbb P[D|\mathcal F_S]\geq\mathbb P[D|H_S\cap \widetilde{H}_S,\mathcal F_S]\mathbb{P}[H_S\cap \widetilde{H}_S|\mathcal F_S]\geq \bm{1}_{H_S}(1-c^{-1}e^{-cS^\alpha})\,,
    \end{equation}
    where we used that $\widetilde{H}_S$ is independent of $\mathcal F_S$ and both $c^{-1}e^{-cT}$ and $c^{-1}e^{-cS^{\frac{1}{3}}}$ were absorbed into $c^{-1}e^{-cS^\alpha}$ by decreasing $c$.
    
    Consider now the event
    \begin{equation}
        D\cap C^{\frac{a}{2},\frac{\varepsilon}{2}}_{S+T}(\widetilde{\A})\,.
    \end{equation}
    On this event, it holds that
    \begin{multline}
    |h_{S+T}([X, Y];\A)-(S+T)(g(\tfrac{X}{S+T})-g(\tfrac{Y}{S+T}))|\leq13S^{\frac13+\frac\alpha2}+2S^{\frac13}
    \\\text{for all $X$ and $Y$ in }[S(\kappa^{-1}+\tfrac\varepsilon2)+\tfrac{2T}{1-b_2}+1,S(\kappa-\tfrac\varepsilon2)]\,,
    \end{multline}
    by repeated use of the triangle inequality.
    For large enough $S$ we have 
    \[
    [S(\kappa^{-1}+\varepsilon),S(\kappa-\varepsilon)] \subset[S(\kappa^{-1}+\tfrac\varepsilon2)+\tfrac{2T}{1-b_2}+1,S(\kappa-\tfrac\varepsilon2)]\,
    \]
    since $T = S^\beta$ for some $\beta<1$.
    Furthermore, for large enough $S$
    \[
    13S^{\frac13+\frac\alpha2}+2S^{\frac13}\leq S^{\frac{1}{3}+\alpha}\,.
    \]
    Therefore, for large enough $S$, the event $C_{S+T}^{\alpha,\varepsilon}(\A)$ contains $D\cap C^{\frac{a}{2},\frac{\varepsilon}{2}}_{S+T}(\widetilde{\A})$.
    Using \eqref{eq:couplingprobability} and noting that $C^{\frac{a}{2},\frac{\varepsilon}{2}}_{S+T}(\widetilde{\A})$ is independent of $\mathcal F_S$ we obtain that
    \begin{align}
    \mathbb P[C_{S+T}^{\alpha,\varepsilon}(\A)|\mathcal F_S]&\geq \mathbb P[D \cap C^{\frac{a}{2},\frac{\varepsilon}{2}}_{S+T}(\widetilde{\A})|\mathcal{F}_S]\\
    &\geq (\Pb[D|\mathcal F_S]-\Pb[C^{\frac{a}{2},\frac{\varepsilon}{2}}_{S+T}(\widetilde{\A})^c|\mathcal{F}_S])\vee0\\
    &\geq(1-c^{-1}e^{-cS^\alpha}-c^{-1}(S+T)^2e^{-c(S+T)^\alpha})\bm{1}_{H_S}\\
    &\geq(1-c^{-1}e^{-cS^\alpha})\bm{1}_{H_S}
    \,,
    \end{align}
    where $(c^{-1}(S+T)^2e^{-c(S+T)^\alpha})\bm{1}_{H_S}$ is absorbed by increasing $c$.
\end{proof}

We can now prove Lemma~\ref{lem:particlesStayClose}.

\begin{proof}[Proof of Lemma \ref{lem:particlesStayClose}]
We can decrease $c$ such that \eqref{eq:particlesStayClose} is trivial for $S$ small and therefore it suffices to consider large enough $S$.
First fix some $X\in[S(\kappa^{-1}+\varepsilon),S(\kappa-\varepsilon)]\cap\mathbb Z$.
By Proposition \ref{prop:hydroEvent}  we have 
\[
\mathbb P[C_{S+T}^{\alpha,\varepsilon/2}(\A)|\mathcal F_S]\geq\bm{1}_{H_S}(1-c^{-1}e^{-cS^\alpha})\,.
\]
The event $C_{S+T}^{\alpha,\varepsilon/2}(\A)$ states a bound for all pairs of points in $[(S+T)(\kappa^{-1}+\varepsilon/2),(S+T)(\kappa-\varepsilon/2)]$.
For $S$ large enough both $\frac{X}{S}(S+T)-S^{1-\gamma}$ and $\frac{X}{S}(S+T)$ are in this interval since $\frac{X}{S}\in[\kappa^{-1}+\varepsilon,\kappa-\varepsilon]$ and $\frac{S^{1-\gamma}}{S+T}\leq \varepsilon/2$ for $S$ large enough.
Since the law of $\B^{(1)}_T$ is equal to the law of $\A_{S+T}$, we have
\begin{multline}\label{eq:hydro1}
    \mathbb P
    \left[\left|
    h_T([\tfrac{X}{S}(S+T)-S^{1-\gamma},\tfrac{X}{S}(S+T)];\B^{(1)})
    -(S+T)(g(\tfrac{X}{S}-\tfrac{S^{1-\gamma}}{S+T})-g(\tfrac{X}{S}))\right|
    \leq (S+T)^{\tfrac{1}{3}+\alpha}
    \bigg|\mathcal F_S\right]\\\geq (1-c^{-1}e^{-cS^\alpha})\bm{1}_{H_S}\,.
\end{multline}

The process $\B^{\text{step},X}$ is also started from step initial data translated by $X$.
For large enough $S$ both $\frac{X}{S}(S+T)-S^{1-\gamma}-X$ and $\frac{X}{S}(S+T)-X$ are in $[T(\kappa^{-1}+\varepsilon/2),T(\kappa-\varepsilon/2)]$.
Indeed $\frac{X}{S}(S+T)-X=\frac{X}{S}T$, and $\frac{S^{1-\gamma}}{T}<\varepsilon/2$ since $T = S^{\beta}$ and $\beta + \gamma >1$ (see \eqref{eq:betaplusgamma}).
Thus we can apply Theorem~\ref{thm:effective_hydrodynamics} to obtain:
\begin{multline}\label{eq:hydro2}
    \mathbb P
    \left[\left|
    h_T([\tfrac{X}{S}(S+T)-S^{1-\gamma},\tfrac{X}{S}(S+T)];\B^{\text{step},X})
    -T(g(\tfrac{X}{S}-\tfrac{S^{1-\gamma}}{T})-g(\tfrac{X}{S}))\right|
    \leq T^{\tfrac{1}{3}+\alpha}
    \right]\geq 1-c^{-1}e^{-cT^\alpha}\,.
\end{multline}

Let us compare the limit shape terms in \eqref{eq:hydro1} and \eqref{eq:hydro2} by bounding their difference
\begin{equation}\label{eq:int1S6V}
    \Delta=T\left(g\left(\frac{X}{S}-\frac{S^{1-\gamma}}{T}\right)-g\left(\frac{X}{S}\right)\right)-\left(S+T\right)\left(g\left(\frac{X}{S}-\frac{S^{1-\gamma}}{S+T}\right)-g\left(\frac{X}{S}\right)\right)\,,
\end{equation}
By the explicit form of $g(x)=\frac{(\sqrt{x}-\sqrt{\kappa})^2}{\kappa -1}$, one can easily see that on the interval $[\kappa^{-1},\kappa]$ both the first and second derivative of $g$ are uniformly bounded and in particular $g''(x)\geq C$ for some $C$ depending only on $\kappa$, for all $x\in[\kappa^{-1},\kappa]$.
Considering the second-order Taylor expansion of $g$ at $\frac{X}{S}$, we have 
\begin{align*}
    T(g(\tfrac{X}{S}-\tfrac{S^{1-\gamma}}{T})-g(\tfrac{X}{S}))
    &=-T(g'(\tfrac{X}{S})\tfrac{S^{1-\gamma}}{T}+\tfrac{g''(\tfrac{X}{S})}{2}(\tfrac{S^{1-\gamma}}{T})^2+O(\tfrac{S^{1-\gamma}}{T}))\\
    &=-g'(\tfrac{X}{S})S^{1-\gamma}+\tfrac{g''(\tfrac{X}{S})}{2}(\tfrac{S^{2-2\gamma}}{T})+O(\tfrac{S^{3-3\gamma}}{T^2}).
\end{align*}
Similarly for the other term in \eqref{eq:int1S6V}
we obtain
\[
(S+T)(g(\tfrac{X}{S}-\tfrac{S^{1-\gamma}}{S+T})-g(\tfrac{X}{S}))=-g'(\tfrac{X}{S})S^{1-\gamma}+\tfrac{g''(\tfrac{X}{S})}{2}(\tfrac{S^{2-2\gamma}}{S+T})+O(\tfrac{S^{3-3\gamma}}{(S+T)^2})
\]
For large enough $S$ the error terms are smaller than the second order terms, since $\frac{S^{1-\gamma}}{T}\to 0$ and we obtain
\begin{equation}\label{eq:particlebound}
    \Delta\geq CS^{2-2\gamma}\left(\frac{1}{T}-\frac{1}{S+T}\right)\geq \frac{C}{2}S^{2-2\gamma}T ^{-1}\geq \frac{C}{2}S^{1/3+2\alpha}
\end{equation}
where the linear terms cancel each other out and we use that $T=S^\beta$ and $\gamma=\frac56-\frac\beta2-\alpha$.

Let $\bm{L}^X$ be defined as
\[
\bm{L}^X=\#\{\textrm{second- and third-class particles in }\B^X_T\textrm{ to the right of }\frac{X}{S}(S+T)-S^{1-\gamma}\}
\]
and note that
\[
\bm{L}^X=\sum_{x\geq \frac{\bm{X}_S}{S}(S+T)-S^{1-\gamma}}\B^{(1,2,3),X}_{T}(x)-\B^{(1)}_T(x)\geq\sum_{x= \frac{\bm{X}_S}{S}(S+T)-S^{1-\gamma}}^{\frac{X}{S}(S+T)}\B^{\text{step},X}_{T}(x)-\B^{(1)}_T(x)\,,
\]
since $\B^{(1,2,3),X}_{T}(x)\geq\B^{\text{step},X}_{T}(x)$ for every $x\in\Z$.
If the events in (\ref{eq:hydro1}) and \eqref{eq:hydro2} take place this sum can be bounded from below by
\begin{multline}
    \sum_{x= \tfrac{\bm{X}_S}{S}(S+T)-S^{1-\gamma}}^{\tfrac{X}{S}(S+T)}\B^{\text{step},X}_{T}(x)-\B^{(1)}_T(x)\\
    =h_T([\tfrac{X}{S}(S+T)-S^{1-\gamma},\tfrac{X}{S}(S+T)];\B^{\text{step},X})-h_T([\tfrac{X}{S}(S+T)-S^{1-\gamma},\tfrac{X}{S}(S+T)];\B^{(1)}))\\
    \geq CS^{\frac13+2\alpha}-(S+T)^{\frac13+\alpha}-T^{\tfrac13+\alpha}\geq S^\frac13\,,
\end{multline}
for $S$ large enough.
Therefore
\[
\mathbb{P}[\bm{L}^X\geq S^\frac{1}{3}|\mathcal F_S]\geq\bm{1}_{H_S}(1-c^{-1}e^{-cT^{\alpha}})\,,
\]
where we use that the event \eqref{eq:hydro2} is independent of $\mathcal F_S$, since it only depends on $\B^{\text{step},X}$, which is only coupled to $\B$ after time $S$.
Using a union bound we obtain
\[
\Pb[\bm{L}^X\geq S^\frac{1}{3}\text{ for all } X\in[S(\kappa^{-1}+\varepsilon),S(\kappa+\varepsilon)]\cap\mathbb Z | \mathcal{F}_S]\geq \bm{1}_{H_S}(1-Sc^{-1}e^{-cT^{\alpha}})
\] 
and we further can absorb $S$ into $c^{-1}e^{-cT^{\alpha}}$ by decreasing $c$.
By the definition of $\bm{L}$ in \eqref{eq:defL}, the definition of $P_S$ in \eqref{eq:bulkSlope}, and the observation \eqref{eq:coupling} we have that
\begin{align*}
\Pb[\bm{L}\geq S^\frac{1}{3}|\mathcal F_S] & \geq  \Pb[\bm{L}^X\geq S^\frac{1}{3}\text{ for all } X\in[S(\kappa^{-1}+\varepsilon),S(\kappa+\varepsilon)]\cap\mathbb Z | \mathcal{F}_S]\bm{1}_{P_S}\\
&\geq \bm{1}_{H_S\cap P_S}(1-c^{-1}e^{-cT^{\alpha}})
\end{align*}
as desired.
\end{proof}

Let us now sketch what needs to be changed for $E^\leq$.
We need to show that on the hydrodynamic event (which is the same) the second-class particle does not deviate too much to the right.
To do so we will delete all particles to the right of the second-class particle, which corresponds to a multi-class particle system $(\B_t)_{t\geq0}$ with the following initial conditions
\begin{itemize}
    \item A first-class particle in every position that is occupied by a first-class particle in $\A_S$ and is to the left of $\bm{X}_S$,
    \item A second-class particle in every position that is occupied by a first-class particle in $\A_S$ that is to the right of $\bm{X}_S$ and
    \item A third-class particle in the position $\bm{X}_S$.
\end{itemize}
This process satisfies the conditions of  Corollary~\ref{cor:dualordering} and by an argument analogous to Lemma~\ref{lem:particlesStayClose} it suffices to show that there are a large number of second-class particles in $\B_T$ to the left of $\frac{\bm{X}_S}{S}(S+T) +S^{1-\gamma}$. 
Denoting by $\B^{(1)}$ the process of the first-class particles in $\B$, by $\B^{(1,2,3)}$ the process of the first, second- and third-class particles, and by $\B^{\text{step}}$ a stochastic six-vertex process with step initial conditions translated to position $\bm{X}_t$.
At time $0$ we have 
\[
\B^{(1)}_0(x)\leq\min (\B^{(1,2,3)}_0(x),\B^{\text{step}}_0(x))\,,
\]
so $\B^{\text{step}}$ can be coupled to $\B^{(1)}$ such that $\B^{(1)}_t(x)\leq \B^{\text{step}}_t(x)$ at all later times $t$ as well.
Note that $\B^{(1)}_t(x)\leq \B^{(1,2,3)}_t(x)$ already holds by definition.
By Proposition~\ref{prop:hydroEvent}, $\B^{(1,2,3)}_T=\A_{S+T}$ is close to the hydrodynamic limit at time $S+T$ with high probability.
Since $\B^{\text{step}}$ is also a stochastic six-vertex model started from step initial conditions we can use Theorem~\ref{thm:effective_hydrodynamics} to say that it is also close to the hydrodynamic limit at time $T$ translated by $\bm{X}_T$ with high probability.
These two results, together with a union bound over all possible values of $\bm{X}_t$ and a calculation similar to \eqref{eq:particlebound}, yield the desired result.

\begin{remark}
    For ASEP one could have simply used particle-hole duality to obtain the proof for $E^\leq$ as a corollary of the proof for $E^\geq$. After exchanging first-class particles with holes, and reversing space, one again obtains a multi-class ASEP, and the events $E^\geq$ and $E^\leq$ are exchanged.
    For the stochastic six-vertex model, this is not the case.
    Applying the particle-hole duality for the stochastic six-vertex model exchanges the two axes and therefore maps the event $E^\geq$ into an event that concerns the times at which the second-class particle hits positions $S$ and $S+T$.
    This is clearly not the same as the event $E^\leq$.
    A different choice of $E^\geq$ and $E^\leq$ such that they are symmetric with respect to the particle-hole symmetry of the stochastic six-vertex model could be considered.
\end{remark}

\section{Symmetry and stationarity of the speed process}\label{sec:speedsymmetries}

In this section, we will prove Corollary \ref{cor:speeddef} and discuss various properties of the stochastic six-vertex model speed process. 

\begin{proof}[Proof of Corollary \ref{cor:speeddef}]
    By the color merging property, the law of $(\bm{X}_t(x))_{t\geq0}$ is equal to that of $(x+\bm{X}_t(0))_{t\geq0}$.
    By Theorem~\ref{thm:main} the speed $\frac{\bm{X}_t(0)}{t}$ converges almost surely, and therefore so does each speed $\frac{\bm{X}_t(x)}{t}$.
    Since there are countably many particles, this also implies that almost surely all of the speeds converge.
\end{proof}

An immediate consequence of the construction is the ergodicity of the speed process.
In this section \textbf{ergodic} always refers to ergodicity with respect to translations of $\Z$, i.e. of space.
\begin{proposition}[Ergodicity for the Speed Process]
The stochastic six-vertex speed process is ergodic.
\end{proposition}
\begin{proof}
    This is immediately inherited from the fact that the process $(\bm{X}_t(x))_{x\in\mathbb Z,t\geq 0}$ can be constructed by sampling i.i.d. pairs of $\text{Bernoulli}(b_1)$ and $\text{Bernoulli}(b_2)$ random variables at every vertex, which are ergodic under the shift.
\end{proof}

To obtain stationarity of the speed process we need the following symmetry, which is a special case of \cite[Corollary 7.1.]{BorodinBufetov2019ColorSymmetry}
\begin{proposition}[Color Position Symmetry in Finite Domains]\label{prop:finitesym}
    Consider the stochastic six-vertex model on an $M\times N$ box, with particles of class $1$ to $M+N$ coming in on the left and lower boundaries such that from the top left to the bottom right the classes are in increasing order. 
    Enumerate the outgoing positions along the top and right boundary with $\{1,\dots,M+N\}$ in descending order, first from left to right along the top and then from top to bottom along the right edge.
    Denote by $\pi$ the (random) permutation of $\{1,\dots,M+N\}$ obtained by letting $\pi(x)$ equal the class of the particle at position $x$.
    Then $\pi$ and $\pi^{-1}$ are equal in law.
\end{proposition}
\begin{proof}
    This follows from \cite[Corollary 7.1]{BorodinBufetov2019ColorSymmetry} by specializing the Ferrer diagram $S$ to a rectangle and using the fact that rectangles are invariant under point reflections.
\end{proof}

We now extend this to the stochastic six-vertex model on the line.
\begin{proposition}[Color Position Symmetry on the Line]
    Consider the random bijection $\pi_N:\Z\to\Z$ obtained by running the stochastic six-vertex model from packed initial conditions until time $N$ (i.e. on a box of infinite width and height $N+1$) and letting $\pi_N(x)$ be the class of the particle exiting at the vertex $(x,N)$.
    Then $(\pi_N(x))_{x\in\Z}$ and $(-\pi_N^{-1}(-x))_{x\in\Z}$ are equal in law.
\end{proposition}
\begin{proof}
    Consider the box $\llbracket-M,M\rrbracket\times \llbracket0,N\rrbracket$.
    Consider the boundary conditions consisting of the incoming arrows from the left with $\{-M-N-1,\dots,-M-1\}$, in increasing order from top to bottom, and incoming arrows from the bottom with classes $\{-M,M\}$ from the bottom, again in increasing order.
    Enumerate the outgoing positions on the top and right boundary with $\{-M-N-1,\dots,M\}$, again in clockwise order, i.e. starting with ${-M-N-1}$ in the bottom right corner and ending with $M$ in the top left corner.
    Again let $\pi_{M,N}(x)$ be the class of the outgoing particle at position $x$.
    Note that there is the following relation between $\pi_{M,N}$ and $\pi_N$ for all $x,y\in{[-M,M]}$
    \begin{equation}
        \mathbb P(\pi_N(x)=y)=\mathbb P(\pi_{M,N}(x)=-y)\,.
    \end{equation}
    Indeed, by using the merging property (Proposition \ref{lem:merging}), one can see that the trace of the particle of class $y$ inside the box $\llbracket-M,M\rrbracket\times \llbracket0,N\rrbracket$ is the same for both models since all particles coming from the left have a smaller class than $y$.
    Thus the probability to exit the box through a specific vertex along the top edge is the same in both models.
    The negative sign on the right-hand side is due to the outgoing boundary positions being enumerated in descending order.
    Using this twice along with Proposition \ref{prop:finitesym}, we obtain
    \begin{equation}
        \mathbb P(\pi_N(x)=y)=\mathbb P(\pi_{M,N}(x)=-y)=\mathbb P(\pi_{M,N}(-y)=x)=\mathbb P(\pi_N(-y)=-x)\,,
    \end{equation}
    which proves the statement.
\end{proof}
Now note that $\pi_N^{-1}(x)= \bm{X}_{N}(x)$, since it is the position of the particle of class $x$ at time $N$.
We can use this to prove that the speed process is stationary with respect to the dynamics of the multi-class stochastic six-vertex process.
\begin{proposition}[Stationarity of the Speed Process]\label{prop:stationarity}
Let $U$ be sampled from the stochastic six-vertex speed process.
Consider the multi-class stochastic six-vertex model with initial conditions given by $(-U(-x))_{x\in\Z}$.
This process is stationary.
\end{proposition}
\begin{proof}
    Start with packed initial conditions and run the process until times $N$ and $N+1$.
    Since $\pi_N(x)$ equals in distribution $-\pi_N^{-1}(-x)$, which equals $-\bm{X}_N(-x)$, we know that both $(\pi_N(x)/N)_{x\in\Z}$ and $(\pi_{N+1}(x)/N)_{x\in\Z}$ converge in law to $(-U(-x))_{x\in\Z}$ by Corollary \ref{cor:speeddef}.
    Let $\mu_N$ be the law of $(\pi_N(x)/N)_{x\in\Z}$ on the space $\R^\Z$ and $\nu_N$ the law of $(\pi_{N+1}(x)/N)_{x\in\Z}$. The laws $\mu_{N+1}$ and $\nu_N$ only differ by multiplying the corresponding random variables with a factor $\frac{N}{N+1}$.
    Since the dynamics do not change under monotone relabeling of classes, and division by $N$ is such a monotone relabeling, we have for any bounded function $f$ on $\R^\Z$ 
    \begin{equation}
        \int f(\eta)d\nu_N(\eta)=\int P_1f(\eta)d\mu_N(\eta)\,
    \end{equation}
    where $P_1$ is the one-step evolution operator for the process.
    Since both $\mu_N$ and $\nu_N$ converge to the law of $(-U(-x))_{x\in\Z}$, this proves the statement.
\end{proof}

The following result is known for multi-class TASEP and ASEP respectively proven in \cite{liggett1976coupling} and \cite{ferrari1991microscopic} respectively.
However, for the stochastic six-vertex model, no proof of this seems to be in the literature, see also \cite[Remark 7.9]{ANPstationarityBaxter}. We prove this in Appendix \ref{ap:unique}.
\begin{proposition}[Uniqueness of Stationary Translation-invariant Measures]\label{prop:uniqueStationaryMeasures}
    For $\lambda_k\in(0,1)$ with $\sum_{k=0}^n\lambda_k=1$, there is a unique ergodic stationary measure for the $n$-class stochastic six-vertex process on the line with $\mathbb P(\eta_0(x)=k)=\lambda_k$.
\end{proposition}
The existence can be derived abstractly from a compactness argument.
Recently in \cite{ANPstationarityBaxter} such measures have also been constructed in a way that is amenable to calculating marginals.

One can conclude the following:
\begin{proposition}\label{prop:ASEPcorrespondencemeasures}
   The ergodic stationary measures for the multi-class stochastic six-vertex model (on the line) are the ergodic stationary measures for the multi-class ASEP.
\end{proposition}
\begin{proof}
    This follows from Proposition \ref{prop:uniqueStationaryMeasures} together with the observation that the stationary measures constructed for both ASEP and the stochastic six-vertex model in \cite{ANPstationarityBaxter} are identical.
    This can be seen by noticing that the stationary measures in \cite[Theorem 3.3]{ANPstationarityBaxter} when specialized to ASEP, as done in \cite[Section 4.2]{ANPstationarityBaxter} and considered on the line instead of the cylinder, are exactly the stationary measures in \cite[Section 7.3]{ANPstationarityBaxter} for the stochastic six-vertex model.
\end{proof}
This result implies the following connection between the stochastic six-vertex speed process and the ASEP speed process.
\begin{proposition}\label{prop:deterministicmapping}
    Let $f$ be the unique increasing map from $[-\kappa,-\kappa^{-1}]$ to $[-1,1]$, which maps a random variable with density $\frac{\sqrt{\kappa}}{2(\kappa-1)}|x|^{-\frac32}\bm{1}_{x\in[-\kappa,-\kappa^{-1}]}$ to a uniform random variable in $[-1,1]$. If $U$ is sampled according to the stochastic six-vertex speed process.
    Then
    \[
    (f(-U(-x)))_{x\in\Z}
    \]
    has the same law as the ASEP speed process.
\end{proposition}
\begin{proof}
    Note first that since $f$ is increasing, and the dynamics only considers the relative ordering of the labels, $(f(-U(-x)))_{x\in\Z}$ is still stationary for the stochastic six-vertex process.
    Given any $(\lambda_k)_{k=1,\dots,n}$ which satisfy $\sum_{k=1}^n\lambda_k=1$, consider the increasing map $\phi_{\lambda}$ from $[-1,1]$ to $\{1,\dots,n\}$ such that $\phi_{\lambda}^{-1}(k)$ is an interval of length $2\lambda_k$.
    By the merging property (Lemma \ref{lem:merging}), $(\phi_\lambda\circ f(-U(-x)))_{x\in\Z}$ is a stationary measure for the $n$-class stochastic six-vertex process.
    By the definition of $\phi_\lambda$ and $f$ it also satisfies
    \[
    \mathbb P[\phi_\lambda\circ f(-U(0))=k]=\lambda_k\,.
    \]
    Ergodicity is also inherited from $U$.
    By Proposition~\ref{prop:uniqueStationaryMeasures} there is only one such measure, and by Proposition~\ref{prop:ASEPcorrespondencemeasures} this is the same as the unique ASEP stationary measure with the same densities.
    If one denotes by $\widetilde{U}$ a sample of the ASEP speed process, this implies the equality in law
    \[
    (\phi_\lambda\circ f(-U(-x)))_{x\in\Z}\stackrel{d}{=}(\phi_\lambda\circ \widetilde{U}(x))_{x\in\Z}
    \]
    This equality is satisfied for all $\phi_\lambda$, and this set of functions is sufficiently large to determine all the marginals of the speed processes, as was outlined in \cite{AmirAngelValko2008TasepSpeed}.
    It follows that
    \[
    (f(-U(-x)))_{x\in\Z}\stackrel{d}{=}(\widetilde{U}(x))_{x\in\Z}
    \]
    as desired.
\end{proof}

\red{\begin{remark}
   Note that this result implies for example that  the sets $\{x:(U(x) = U(0)\}$ and $\{x:(\widetilde{U}(x) = \widetilde{U}(0)\}$ have the same law. 
   Therefore, the results about the convoy structures of the ASEP speed process in \cite{AmirAngelValko2008TasepSpeed} and \cite{MartinStationaryASEP} hold here as well. 
\end{remark}}

\section{Effective hydrodynamic estimates}\label{sec:effectivehydrodynamics}

The purpose of this section is to prove Propositions \ref{prop:MeixnerTail} and \ref{prop:otherTail}. Before doing that, we combine them to prove the following theorem:  

Recall that $H(X,T)$ refers to the height function of a stochastic six-vertex model on the quadrant with step initial conditions and that $g(x)=g(x,1)$ is the limit shape of the height function (see \eqref{eq:asconverge}).

\begin{theorem} \label{thm:effective_hydrodynamics}
For any $\varepsilon>0$, there exists $c=c(\varepsilon)>0$ such that the following holds. For any $\mu_1,\mu_2\in[\kappa^{-1} + \varepsilon,  \kappa - \varepsilon]$, and for any $T\geq1$, $s  \in [0, T]$,
\begin{align}\label{eq:effective_hydrodynamics}
&\mathbb{P}\left[\left| H(T\mu_1, T) - H(T\mu_2, T)  - \left( g(\mu_1) - g(\mu_2)\right)T\right |  \geq  sT^{1/3}\right]  \leq c^{-1}e^{-cs}\,.
\end{align}
Furthermore, the constant $c$ can be chosen to weakly decrease in $\varepsilon$.
\end{theorem}

\begin{proof}[Proof of Theorem \ref{thm:effective_hydrodynamics}]
For any $\mu\in[\kappa^{-1}+\varepsilon,\kappa-\varepsilon]$ we have the following two bounds from Propositions \ref{prop:MeixnerTail} and \ref{prop:otherTail}, respectively. 
There exists a $c$ (that will change from line to line) such that
\begin{align*}
\mathbb{P}\left[H(T\mu, T) \geq g(\mu)T + sT^{1/3}\right] &\leq c^{-1} e^{-c s^\frac23} \\
\mathbb{P}\left[H(T\mu, T) \leq g(\mu)T - sT^{1/3}\right] &\leq c^{-1} (e^{-c s} 
+ e^{-cT}) \leq 2c^{-1}e^{-c s}.
\end{align*}
Combining these two bounds, we obtain
\begin{align*}
\mathbb{P}\left[\left| H(T\mu, T) - g(\mu)T\right |  \geq  sT^{1/3}\right] &\leq c^{-1} e^{-c s}\,.
\end{align*}
It follows from a union bound that 
\begin{align*}
&\mathbb{P}\left[\left| H(T\mu_1, T) - H(T\mu_1,T)  - \left(g(\mu_1) - g(\mu_2)\right)T\right |  \geq  sT^{1/3}\right] \\
& \leq \mathbb{P}\left[\left| H(T\mu_1, T )   - g(\mu_1)T\right| \geq  \frac{s}{2}T^{1/3}\right] +\mathbb{P}\left[\left| H(T\mu_2, T )   - g(\mu_2)T\right| \geq  \frac{s}{2}T^{1/3}\right]  \\
& \leq c^{-1}e^{-cs}.
\end{align*}
This finishes the proof of Theorem~\ref{thm:effective_hydrodynamics}.
The constant $c$ can be chosen to be weakly decreasing in $\varepsilon$ since this is the case for both Proposition~\ref{prop:MeixnerTail} and Proposition~\ref{prop:otherTail}.
\end{proof}

We immediately obtain the following corollary of Theorem \ref{thm:effective_hydrodynamics}:

\begin{corollary}\label{cor:effective_hydrodynamics}
For any $\varepsilon>0$, there exists $c=c(\varepsilon)>0$ such that the following holds. For any $T \geq 1$ and for any $s  \in [0, T]$,
\begin{align}
&\mathbb{P}\left[\max_{\mu_1,\mu_2\in[\kappa^{-1} + \varepsilon ,  \kappa - \varepsilon]} \left| H(T\mu_1, T) - H(T\mu_2, T)  - \left( g(\mu_1) - g(\mu_2)\right)T\right |  \geq  sT^{1/3}\right]  \leq c^{-1}T^2e^{-cs},
\end{align}
and $c$ can be chosen to weakly decrease in $\varepsilon$.
\end{corollary}

\begin{proof}
    Notice that there are only finitely many $\mu_i$ satisfying $\kappa^{-1} + \varepsilon \leq \mu_i \leq  \kappa - \varepsilon$ and such that $T \mu_i $ is an integer. In fact, there are at most $\kappa T$ of them, giving at most $O(T^2)$ possible pairs $(\mu_1, \mu_2)$. Taking a union bound of \eqref{eq:effective_hydrodynamics} over all such pairs yields the result.
\end{proof}

 Finally, we can quickly extend Proposition \ref{prop:MeixnerTail} to the case of step Bernoulli boundary conditions, i.e. the incoming arrows from the left are given by i.i.d. Bernoulli$(\rho)$ random variables, while the incoming positions from the bottom are all still empty.
 Denote these boundary conditions as $(\rho,0)$-Bernoulli boundary conditions.
 Even though we don't need this result to prove our main theorem, we state it as a corollary for completeness.

\begin{corollary}
\label{cor:upperTailBernoulli}
Fix $\varepsilon>0$. There exists a constant $c=c(\varepsilon)>0$ such that the following holds: Let $\rho \in[0, 1]$ and let $H^\rho(x,y)$ be the height function for the stochastic six-vertex model on the quadrant with $(\rho ; 0)$-Bernoulli boundary conditions. For any $\mu\in[\kappa^{-1}+\varepsilon,\kappa-\varepsilon]$ and for any $T\geq1$, $s  \geq 0$, 
\begin{equation}
\mathbb{P}\left[H^\rho(T\mu, T) \geq g(\mu)T + sT^{1/3}\right] \leq c^{-1} e^{-c s^\frac32} . 
\end{equation}

\end{corollary}

\begin{proof}
    This is a straightforward consequence of the attractivity of the stochastic six-vertex model, by which we can couple the model with $(\rho;0)$-Bernoulli initial data with the model with step initial data.
    In this coupling the height function of the model with $(\rho;0)$-Bernoulli initial data is smaller at every point, and thus the statement follows from Proposition~\ref{prop:MeixnerTail}.
\end{proof}

\begin{remark} \label{rmk:extendinglowertail}
Proposition \ref{prop:otherTail} can also be extended to the case of $(\rho,0)$-Bernoulli boundary conditions by restricting $\mu\in\left[\frac{\kappa}{(\kappa\rho-\rho+1)^2}+\varepsilon,\kappa-\varepsilon\right]$, since $\left[\frac{\kappa}{(\kappa\rho-\rho+1)^2},\kappa\right]$ is the rarefaction fan for $(\rho,0)$-Bernoulli boundary conditions. In fact, we actually prove this more general version in the proof of Proposition \ref{prop:otherTail}.
\end{remark}

\subsection{Proof of Proposition~\ref{prop:MeixnerTail}}\label{sec:Tail1}
To prove Proposition~\ref{prop:MeixnerTail}, we will make use of a remarkable exact identity that relates the height function of the stochastic six-vertex model to the holes of the Meixner ensemble. We can then reduce the question of studying the tail of the height function to studying the tail for the position of the smallest hole in this determinantal point process. We study this tail by taking asymptotics of the associated kernel. 
\red{This is an example of a general method of extracting asymptotics for determinantal point process via contour integrals and saddle point analysis.
For an introduction to this method, see for example \cite{okounkov2002symmetric, MR3468920}.}

In this subsection, we define the Meixner ensemble, which is a determinantal point process on $\Z$. We will then relate the $q$-Laplace transform of the stochastic six-vertex model height function to an expectation with respect to the Meixner ensemble. 

We give a brief introduction to the theory of discrete determinantal point processes. Let $\mathfrak{X}$ denote the state space of a single particle, which we will take to be a countable set (for the Meixner ensemble, we will take $\mathfrak{X} = \Z_{\geq 0}$). A subset $X \subseteq \mathfrak{X}$ is called a \textit{point configuration}, and we define $\text{Conf}(\mathfrak{X}) = 2^{\mathfrak{X}}$ to be the set of all possible point configurations. 

We define the following Borel sigma algebra for $\conf$:
$$\mathcal{B} : = \sigma\left(\left\{X \in \conf: |A \cap X|=n\right\}:n \in \mathbb{N}, A \subseteq \XX \text { compact }\right).$$
A probability measure $\Pb$ on $(\conf, \mathcal{B})$ is called a \textit{random point process}. From now on, we will use $X$ to denote this random point process by setting $X: \conf \to \conf$, $X(\omega) = \omega$.

We define the \textit{$n$-point correlation function} as follows: for $A = \{x_1, ..., x_n\} \subseteq \XX$, let 

$$\rho_n(A) = \rho_n(x_1, ..., x_n) := \Pb \left[A \subseteq X \right].$$ 

\begin{definition}[Determinantal Point Process]
A random point process $X$ is determinantal if there exists a kernel $K: \XX \times \XX \to \R$ such that for all $n \geq 1$ and for all $x_1, \ldots, x_n$ with $x_i \neq x_j$ for $i \neq j$, we have 

\begin{equation}
    \rho_n(x_1, ..., x_n) = \det \left(K(x_i, x_j)_{i,j =1}^n\right)
\end{equation}
\end{definition}

Let $W(x): \XX \to \R$ be a weight function, and let $P_0, P_1, \ldots$ be the family of orthonormal polynomials with respect to $W$, i.e.,

$$\int_{\XX} P_i(x)P_j(x)W(x)dx = \ind_{i = j}.$$ 
The corresponding \textit{$N$-point ensemble} (a random point process where $\Pb$ is supported on configurations with exactly $N$ particles) is given by 

$$\Pb(x_1, \ldots, x_N) \propto \det(V(x_1, \ldots, x_N))^2 \prod_{i=1}^NW(x_i),$$
where $V(x_1, \ldots, x_N) = (x_i^{j-1})_{i, j =1}^N$ is the Vandermonde matrix, and $\det(V(x_1, \ldots, x_N)) = \prod_{i < j}(x_j - x_i)$ is the Vandermonde determinant. An $N$-point ensemble generated in this way is determinantal with the \textit{Christoffel-Darboux} kernel
\begin{equation}\label{eq:orthogonalKernel}
     K_N(x,y) = (W(x) W(y))^{\frac{1}{2}} \sum_{n=0}^{N-1}P_n(x) P_n(y).
\end{equation} 

The Meixner polynomials are a family of orthogonal polynomials on $\Z_{\geq 0}$. We fix two parameters: $\beta > 0$ and $\xi \in (0,1)$, and then define the weight function $W: \Z_{\geq 0} \to \R$: 

\begin{equation}\label{eq:MeixnerWeight}
    W(x) = \frac{\Gamma(\beta + x)}{\Gamma(\beta)x!}\xi^x. 
\end{equation}
We can then define the Meixner polynomials to be the family of orthogonal polynomials with respect to the weight function $W$. Using these orthogonal polynomials, we can define the Meixner ensemble $\text{Meixner}(N, \beta,\xi)$ to be the corresponding $N$-point ensemble.

We use the term \textit{particles} to refer to the elements of a point process $X$ and use the term \textit{holes} to refer to elements of $\XX \setminus X$. \textit{Particle-hole involution} is an involution from $\conf \to \conf$ that exchanges particles with holes. In other words, $X \mapsto X^{\circ} := \XX \setminus X$. If we start with an $N$-point ensemble, then particle-hole involution yields a point process with infinitely many particles. Suppose that $X$ is a determinantal point process with kernel $K$. Then $X^{\circ}$ is a determinantal point process with kernel $1- K$.

Next, we give a brief overview of Schur measures. An integer partition is denoted as $\lambda = (\lambda_1, \lambda_2 ,\ldots )$ where $\lambda_1 \geq \lambda_2 \geq \cdots$ and $l(\lambda)$ denotes the number of nonzero $\lambda_i$ in the partition $\lambda$. Let $\mathbb{Y}$ denote the set of all integer partitions. Let $\textbf{x} = (x_1, x_2 \ldots), \textbf{y} = (y_1, y_2, \ldots)$ be two sets of nonnegative variables. For fixed $\textbf{x}$ and $\textbf{y}$, we define the Schur measure $\textbf{SM}(\textbf{x};\textbf{y})(\lambda)$ as a measure on partitions $\lambda$ as follows:

\begin{equation} \label{eq:SchurMeasure}
    \textbf{SM}(\textbf{x};\textbf{y})(\lambda):= \frac{s_{\lambda}(\textbf{x})s_{\lambda}(\textbf{y})}{\Pi (\textbf{x};\textbf{y})}
\end{equation}
where $s_{\lambda}$ is the Schur symmetric function indexed by $\lambda$ and $\Pi(\textbf{x};\textbf{y}) = \sum_{\lambda}s_{\lambda}(\textbf{x})s_{\lambda}(\textbf{y})$ is the partition function. We need to assume that $\Pi(\textbf{x};\textbf{y}) < \infty$ for our choice of $\textbf{x}$ and $\textbf{y}$ for this to define a valid probability measure.

The Meixner ensemble can be obtained as a pushforward of the Schur measure as follows: Consider the Schur measure of the form $\textbf{SM}(x, \ldots ,x; y\ldots, y)$ where we take $n$ copies of $x$ and $m$ copies of $y$. Using standard properties of Schur functions, it follows that this measure is supported on $\mathbb{Y}^{(\min \{m,n\})}$, which is the set of partitions with $l(\lambda)\leq \min \{m,n\}.$ Finally, consider the map from $\mathbb{Y}^{(\min \{m,n\})} \to \text{Conf}(\Z_{\geq 0})$ such that $\lambda \mapsto \{\min \{m,n\}  + \lambda_i - i\}_{i=1}^{ \min \{m,n\}}.$ Then the pushforward of $\textbf{SM}(x, \ldots ,x; y\ldots, y)$ to a measure on $\text{Conf}(\Z_{\geq 0})$ gives us the Meixner ensemble $\text{Meixner}(\min \{m,n\}, |m-n| +1, xy).$ This can be checked directly, see Proposition 8.2 in \cite{MR3649488}.

The following identity originates from \cite{MR3760963}, although we state a version written in \cite{MR3649488}: Let $\mathbb{E}_{\mathbf{6} \mathbf{v}}$ refer to the expectation with respect to the stochastic six-vertex model and let $\E_{\textbf{SM}}$ denote the expectation with respect to a specified Schur measure. 

\begin{proposition}[Proposition 8.4 in \cite{MR3649488}] \label{prop:schurMeasureId1}
Take any $0<q<1$ and $\kappa>1$ and consider the stochastic six-vertex model on the quadrant parameterized by $q$ and $\kappa$. Consider any integers $M, N \geq 1$. Then for any $\xi \notin-q^{\mathbb{Z}_{\leq 0}}$ we have
\begin{equation}\label{eq:schurIdentity}
\mathbb{E}_{\mathbf{6} \mathbf{v}} \prod_{i \geq 0} \frac{1}{1+\xi q^{H(M, N)+i}}=\mathbb{E}_{\mathbf{S M}} \prod_{j \geq 0} \frac{1+\xi q^{\lambda_{N-j}+j}}{1+\xi q^{j}}    
\end{equation}
where in the right-hand side we assume that $q^{\lambda_{-m}}=0$ for $m \geq 0$, and the expectation is with respect to the Schur measure $\mathbf{S M}(\underbrace{\kappa^{-1} q^{1/2}, \ldots, \kappa^{-1} q^{1/2}}_{N} ; \underbrace{q^{-1 / 2}, \ldots, q^{-1 / 2}}_{M-1})$.
\end{proposition}

If $M > N$, then the Schur measure in \eqref{eq:schurIdentity} is supported on $\mathbb{Y}^{(N)}$. We can obtain the $N$ particles of the $\text{Meixner}(N, M- N, \kappa^{-1})$ ensemble by taking the above-mentioned pushforward of the Schur measure so that the particles in the Meixner ensemble are given by $\{\lambda_i + N-i\}_{i=1}^N = \{\lambda_{N-j} + j\}_{j =0}^{N-1}.$ On the other hand, if $M \leq N$, then the Schur measure is supported on $\mathbb{Y}^{(M-1)}$. We now have that 
$$\{\lambda_{N-j} + j\}_{j =0}^{N-1} = \{0, \cdots, N-M\} \sqcup \{\lambda_i + (N-M) + M -i\}_{i =1}^{M-1}.$$
This gives us the $M-1$ particles in $\text{Meixner}(M-1, N - M + 2, \kappa^{-1})$ each shifted over deterministically by $N-(M-1)$ along with the addition of particles packed from $0$ to $N-M$. In either case, we can obtain the following identity: 
\begin{proposition}\label{prop:schurMeasureId2}
Take any $0<q<1$ and $\kappa>1$ and consider the stochastic six-vertex model on the quadrant parameterized by $q$ and $\kappa$. Consider any integers $M, N \geq 1$. Then for any $\xi \notin-q^{\mathbb{Z}_{\leq 0}}$ we have
\begin{align}
    \mathbb{E}_{\mathbf{S M}} \prod_{j \geq 0} \frac{1+\xi q^{\lambda_{N-j}+j}}{1+\xi q^{j}} = \mathbb{E}_{X} \prod_{x \in  X} \frac{1}{1+\xi q^{x}}. \label{eq:MeixnerHoles}\end{align}
where in the left-hand side we assume that $q^{\lambda_{-m}}=0$ for $m \geq 0$ and the right-hand expectation is with respect to the point process 
\begin{equation}\label{eq:Meixnerparameters}
     X \sim \begin{cases}
    \mathrm{Meixner}^{\circ}(N, M- N,\kappa^{-1}) & \text{if $M > N$} \\
    N - (M-1) + \mathrm{Meixner}^{\circ}(M-1, N - M + 2, \kappa^{-1}) & \text{if $M  \leq  N$}\,,
\end{cases}
\end{equation}
where for a point process $X$, $n+X$ denotes the point process obtained by deterministically shifting over each particle in $X$ by $n$, and as mentioned above, $X^{\circ}$ denotes the holes of $X$.  
\end{proposition}

\begin{proof}
  The proof of this follows from crossing out each term in the denominator that equals one of the nontrivial terms in the numerator (i.e., a term corresponding to one of the particles in the Meixner ensemble). All terms that remain in the denominator will correspond to holes of the associated Meixner ensemble.   
\end{proof}

We now explain how we go from Propositions \ref{prop:schurMeasureId1} and \ref{prop:schurMeasureId2} to proving Proposition \ref{prop:MeixnerTail}. We will first need the following definition and lemma:

\begin{definition}[q-Pochhammer symbol]For any complex numbers $q$ and $a$ such that $|q| < 1$, we define $(a; q)_{\infty} = \prod_{j=0}^{\infty} (1- aq^j)$. 
\end{definition}

The following Lemma is taken from \cite{ACG2023asepspeed} and it allows us to connect the $q$-Laplace transform of $H$ (the left-hand side of \eqref{eq:schurIdentity}) with the tail probability of $H$. 
\begin{lemma}[Lemma B.7 in \cite{ACG2023asepspeed}]
Let $\mathbf{A}$ be a real-valued random variable, $q \in[0,1)$ and $b \in \mathbb{R}$. Then,
\begin{align}
& \mathbb{P}[\mathbf{A} \leq 0] \leq 2 \cdot\left(1-\mathbb{E}\left[\left(-q^{\mathbf{A}} ; q\right)_{\infty}^{-1}\right]\right), \label{q1}\\
& \mathbb{E}\left[\left(-q^{\mathbf{A}} ; q\right)_{\infty}^{-1}\right] \geq e^{q^b /(q-1)} \cdot \mathbb{P}[\mathbf{A} \geq b], \label{q2} \\
& \mathbb{E}\left[\left(1+q^{\mathbf{A}}\right)^{-1}\right] \leq \mathbb{P}[\mathbf{A}>-b]+q^b \cdot \mathbb{P}[\mathbf{A} \leq-b] \label{q3} .
\end{align}
\end{lemma}

\begin{remark}
    While in the statement of Proposition~\ref{prop:MeixnerTail} we consider a height function of the form $H(T\mu, T)$, for the remainder of this section we will work with the more general form $H(T\mu, T\nu)$ in order to highlight that many of the formulas that we will use in our analysis will have some symmetries in $\mu$ and $\nu$. In the end, we will simply take $\nu = 1$.  This does not actually reduce generality, since any appropriate $H(M, N)$ can be obtained by taking $T=N$ and $\mu = M/T$.
    
    The constants in this section are allowed to depend on $\kappa$ freely, but can be chosen to be uniform in $\mu$ and $\nu$ as long as $\kappa^{-1}+\varepsilon\leq\mu\leq \kappa -\varepsilon$ and $\nu=1$ (Any other compact set bounded away from the two lines $\frac{\mu}{\nu}=\kappa^{-1}$ and $\frac{\mu}{\nu}=\kappa$ would also work).
    In particular, this will also be true for all implicit constants hidden in big $O$ notation terms.
\end{remark}

Take $M = T \mu, N = T \nu$, and $\xi = q^{g(\mu, \nu)T - sT^{1/3}}$. Then using \eqref{q2}, with  $\mathbf{A} = H(T\mu, T \nu) -g(\mu, \nu)T - sT^{1/3}$ and $b = 0$, we obtain 
\begin{align}
 \mathbb{P}\left[H(T\mu, T \nu) \geq g(\mu, \nu)T + sT^{1/3}\right] &\leq  e^{-1 /(q-1)} \cdot\mathbb{E}_{\mathbf{6} \mathbf{v}} \prod_{i \geq 0} \frac{1}{1+\xi q^{H(M, N)+i}} \\ \label{eq:comparison1}
 &= e^{-1 /(q-1)} \cdot \mathbb{E}_{X} \prod_{x \in  X} \frac{1}{1+\xi q^{x}}
\end{align}
where the point process $X$ is defined as in Proposition \ref{prop:schurMeasureId2}. We can estimate the last product by dropping all terms in the product except for that corresponding to the smallest hole $x_1$. More precisely, since all the terms in the product are at most $1$, we have 

\begin{align}\label{eq:comparison2}
 \mathbb{E}_{X} \prod_{x \in X} \frac{1}{1+\xi q^{x}}& \leq \frac{1}{1 + \xi q^{x_1}}.
\end{align}

Using \eqref{q3} with $\mathbf{A} = x_1 -g(\mu, \nu)T - sT^{1/3}$ and $b = \frac{sT^{1/3}}{2}$, we see that

\begin{align}
    \frac{1}{1 + \xi q^{x_1}} &\leq \Pb\left[x_1 > g(\mu, \nu)T +\frac{sT^{1/3}}{2}\right] + q^{\frac{sT^{1/3}}{2}} \Pb\left[x_1 \leq g(\mu, \nu)T +\frac{sT^{1/3}}{2}\right]\\
   \label{eq:comparison3} &\leq  \Pb\left[x_1 > g(\mu, \nu)T +\frac{sT^{1/3}}{2}\right] + q^{\frac{sT^{1/3}}{2}}.
\end{align} 

So in order to obtain an upper bound on $ \mathbb{P}\left[H(T\mu, T \nu) \geq g(\mu, \nu)T + sT^{1/3}\right]$, it will suffice to obtain an upper bound on $\Pb\left[x_1 > g(\mu, \nu)T +\frac{sT^{1/3}}{2}\right]$. Let us denote the holes of the Meixner ensemble by $x_1, x_2,...$. We know that $\lambda$ has at most $N$ nonzero parts. If there are only $k$ nonzero parts, then $\lambda_{k+1}, ... , \lambda_N = 0$, so there are $N-k$ Meixner particles at positions $0, ..., N-k-1$. Therefore, the smallest hole $x_1$ will occur at position $N-k$. It follows that $\ell(\lambda) = N- x_1$. Equivalently, we have 

\begin{equation}x_1 = N - \ell(\lambda). \label{eq:particle-length}\end{equation}

It follows from \eqref{eq:particle-length} that 
\begin{align}\label{eq:comparison4}
&\Pb\left[x_1 >g(\mu, \nu)T +\frac{sT^{1/3}}{2}\right] = \Pb\left[-\ell(\lambda) >\left(g(\mu, \nu) - \nu \right) T +\frac{sT^{1/3}}{2}  \right] .
\end{align} 

According to \cite[In the proof of Theorem 6.1]{MR3760963} we can represent the tail probability $\Pb \left[-\ell(\lambda) > h\right]$ as a Fredholm determinant. We first recall the definition of a Fredholm determinant:

\begin{definition}[Fredholm Determinant (Discrete)]
We define the Fredholm determinant of a kernel $K: \Z \times \Z \to \R$ as 
$$
\operatorname{det}(\operatorname{Id}+K)=1+\sum_{k=1}^{\infty} \frac{1}{k!}\sum_{x_1, \ldots, x_k \in \Z} \operatorname{det}\left[K\left(x_i, x_j\right)\right]_{i, j=1}^k.
$$
\end{definition}

We then have 
\begin{equation}\label{eq:determinantalForm}
\mathbb P[-\ell(\lambda) > h] = \operatorname{det}(\operatorname{Id}- \widetilde{K})_{\ell^2(h, h-1, ...)} \coloneqq \operatorname{det}(\operatorname{Id}-\Pi_h\widetilde{K}\Pi_h).
\end{equation} where $\widetilde{K}$ is a correlation kernel obtained as a dual of the kernel for the Meixner ensemble and $\Pi_h$ is the projection from $\ell^2(\Z)$ to $\ell^2(h,h-1,\dots)$.

We can write out $\widetilde{K}$ explicitly as (see \cite[Equation (6.1) and the subsequent paragraph in the reference]{MR3760963})
\begin{equation}\label{eq:FredholmMeixner}
    \widetilde{K}(x,y) = \frac{1}{(2 \pi i)^2} \oint \oint \frac{(\sqrt{\kappa} - z^{-1})^N}{(\sqrt{\kappa}- z)^{M-1}} \frac{(\sqrt{\kappa} -  w)^{M-1}}{(\sqrt{\kappa}- w^{-1})^{N} } \frac{dz dw}{(w-z) z^{x+1}w^{-y}}
\end{equation}
where  $x,y \in \Z$ and the integrals are taken over positively oriented circular contours with $1/{\sqrt{\kappa}} < |z| < 1 < |w| < \sqrt{\kappa}.$ Note that our integrand has poles at $0, 1/{\sqrt{\kappa}}$ and $\sqrt{\kappa}$, so the contours are chosen so that they do not pass through the poles.

To estimate the Fredholm determinant in \eqref{eq:determinantalForm} we will use a technique known as Widom's trick first used in \cite[Lemma 1]{widom2002convergence}.
It consists of the observation that for a kernel $K$ with eigenvalues in $[0,1]$ it holds that
\begin{equation}
   \operatorname{det}(\operatorname{Id}- K) \leq \exp(-\text{Tr}(K))\,.
\end{equation}
The following lemma checks that the operator $\Pi_h\widetilde{K}\Pi_h$ satisfies this condition.

\begin{lemma}\label{lem:projections}
    The operator $\Pi_h\widetilde{K}\Pi_h$ has real eigenvalues $(\mu_j)_{j\geq0}$ all of which are in $[0,1]$ and hence
    \begin{equation}\label{eq:widomtrick}
   \operatorname{det}(\operatorname{Id}- \Pi_h\widetilde{K}\Pi_h) \leq \exp(-\text{Tr}(\Pi_h\widetilde{K}\Pi_h))\,.
\end{equation}
\end{lemma}
\begin{proof}
    Let $I(x,y)=\bm{1}_{x=y}$ be the identity operator and $K$ as in \cite[Equation (6.1)]{MR3760963}. 
    As noted in \cite[Below Equation (6.1)]{MR3760963} the operator $\widetilde{K}$ satisfies $\widetilde{K}=\operatorname{Id}-K$.
    The operator $K$ is related to the Christoffel-Darboux kernel $K_N$ (see \eqref{eq:orthogonalKernel}) associated to the Meixner ensemble via a gauge transformation, see Theorem 3.3 and Lemma 3.5 in \cite{BorOlsh2006Meixner}.
    Since $\operatorname{Id}$ is invariant under gauge transformations, this also means that $\widetilde{K}$ is related to $\operatorname{Id} - K_N$ via a gauge transformation.
    The operator $K_N$ is a projection operator since it is a Christoffel-Darboux kernel.
    Therefore, $\operatorname{Id}-K_N$ is also a projection operator and finally $\widetilde{K}$ as well, since a gauge transform of a projection operator is a projection operator.
    After the gauge transformation, $\Pi_h\widetilde{K}\Pi_h$ becomes self-adjoint, so the eigenvalues are real and non-negative.
    Since we have now also seen that this operator is a composition of projections, the eigenvalues are at most $1$.
    Since $1-x\leq\exp(-x)$ for $x\in[0,1]$ this implies \eqref{eq:widomtrick}.
\end{proof}

Therefore, obtaining an upper bound on $\mathbb P[-\ell(\lambda) > h]$ reduces to obtaining a lower bound for $\text{Tr}(\Pi_h\widetilde{K}\Pi_h)$.
Similar kinds of bounds were obtained for other kernels in e.g. \cite{widom2002convergence,BaikFerrariPeche2014twopointTASEP}.
We first compute this trace in the following lemma.
\begin{lemma}
We have
\begin{align} \label{eq:contourIntegral}
    \mathrm{Tr}(\Pi_h\widetilde{K}\Pi_h) &= \frac{1}{(2 \pi i)^2} \oint \oint \exp (T (G_x(z) - G_x(w))\frac{dz dw}{(w-z)^2 }\,,
\end{align}
where
\begin{equation}
    G_x(z) = \nu \ln\left(\sqrt{\kappa} - z^{-1}\right)  - \mu \ln\left(\sqrt{\kappa} - z\right) - x \ln(z)
\end{equation}
and we have reparameterized $M,N$ and $h$ as
\begin{align*}
    \mu &= \frac{M-1}{T};  \quad \nu = \frac{N}{T}; \quad x = \frac{h}{T}\,.
\end{align*}    
\end{lemma}
\begin{proof}
The trace is given by
\begin{align*}
    \text{Tr}(\Pi_h\widetilde{K}\Pi_h) &= \sum_{j = -\infty}^h \widetilde{K}(j,j)\\
    &= \sum_{j = -\infty}^h \frac{1}{(2 \pi i)^2} \oint \oint \frac{(\sqrt{\kappa}  -  z^{-1})^N}{(\sqrt{\kappa} - z)^{M-1}} \frac{(\sqrt{\kappa}  -  w)^{M-1}}{(\sqrt{\kappa} - w^{-1})^{N} }\left(\frac{w}{z}\right)^j \frac{dz dw}{(w-z) z}.
\end{align*}
 Since $|w/z| > 1$ by our choice of contours, we can sum $w/z$ from $-\infty$ to $h$ which yields 
\begin{align*}
    \text{Tr}(\Pi_h\widetilde{K}\Pi_h) &= \frac{1}{(2 \pi i)^2} \oint \oint \frac{(\sqrt{\kappa} -  z^{-1})^N}{(\sqrt{\kappa} - z)^{M-1}} \frac{(\sqrt{\kappa}  - w)^{M-1}}{(\sqrt{\kappa} - w^{-1})^{N} } \left(\frac{w}{z}\right)^h\frac{dz dw}{(w-z)^2 }.
\end{align*}
Finally, we can rewrite the integrand in exponential form to obtain \eqref{eq:contourIntegral}.
\end{proof}

The function $G$ has the following two critical points:
\begin{align}
    z_c^{\pm} &= \frac{ \mu+\nu +(\kappa + 1)x  \pm \sqrt{-4\kappa(\mu+x)(\nu+x) + ( \mu+\nu +(\kappa + 1)x )^2}}{2\sqrt{\kappa}(\mu+x)} 
\end{align}
We can see that these two critical points are equal if we choose $x=x_c^{\pm}$ where 
\begin{align}
    x_c^{\pm} &= \frac{(\sqrt{\mu} \pm \sqrt{\kappa \nu})^2}{\kappa -1 } - \nu.
\end{align}
Note that $x_c^{-} = g(\mu, \nu) - \nu$. For $x = x_c^{-}$, we have 
\begin{align*}
    z_c^+ = z_c^{-} = \frac{\sqrt{ \mu} - \sqrt{\kappa \nu}}{\sqrt{\kappa \mu} - \sqrt{\nu}}.
\end{align*}
Denote this value as $z_c$.
For general $x$, we can rewrite the formula for $z_c^{\pm}$ as 
\begin{align}
    z_c^{\pm} &= \frac{ \mu+\nu +(\kappa + 1)x \pm \sqrt{(\kappa - 1)^2(x - x_c^{+})(x - x_c^{-})}}{2\sqrt{\kappa}(\mu+x)}.
\end{align}

If $x_c^{-} < x <  x_c^{+}$,
then the two critical points $z_c^{\pm}$ are not real. 
Then it holds that
\begin{equation}\label{eq:criticalpointsdistance}
    |z_c^{\pm}| = \sqrt{\frac{\nu + x}{\mu +x}}.
\end{equation}

The following lemma describes how this function behaves around $(x_c^-,z_c)$.
\begin{lemma}\label{lemma:taylor of G}
    The function $G_x(z)$ satisfies:
    \begin{equation*}
        G_{x_c^-}'''(z_c)=2\frac{\sqrt{\kappa\mu\nu}(\sqrt{\kappa}-\sqrt{\mu/\nu})^2(\sqrt{\kappa}-\sqrt{\nu/\mu})^2}{z_c^3(\kappa-1)^3}.
    \end{equation*}
    As $x\rightarrow x_c^-$ from above we have the following:
    \begin{align}
        z_c^+&=z_c+\frac{i\left(1-\kappa^{-1}\right)^\frac12(\kappa^{-1}\mu\nu)^\frac14\sqrt{x-x_c^-}}{\sqrt{\kappa}(\kappa \mu + \nu + 2 \sqrt{\kappa \mu \nu})} +O(x-x_c^-) \text{ and }\label{eq:critdistance}\\
        G_x''(z_c^+)&=\frac{iG_{x_c^-}'''(z_c)(1-\kappa^{-1})^\frac12(\kappa^{-1}\mu\nu)^\frac14\sqrt{x-x_c^-}}{\sqrt{\kappa}(\kappa \mu + \nu + 2 \sqrt{\kappa \mu \nu})}+O(x-x_c^-)\,,
    \end{align}
    where the implicit constant in the big $O$ term can be chosen independently of $\mu\in[\kappa^{-1}+\varepsilon,\kappa-\varepsilon]$ and $x\in [x_c^-,x_c^+]$, i.e. it depends only on $\kappa$ and $\varepsilon$.
\end{lemma}
\begin{proof}
    The first two equalities are calculations, the third one is the Taylor expansion of $G''$ in $x$ and $z$ around $(x_c^-,z_c)$.
\end{proof}

Now that we have established all the variables at play we can state an estimate on the trace. 
\begin{proposition} \label{prop:trace}
    Define $s=2(x-x_c^-)T^\frac23$.
    For any $\varepsilon$ there exist $s_0,T_0$ and $C$ such that for any $\mu\in[\kappa^{-1}+\varepsilon,\kappa-\varepsilon]$ and $\nu=1$, any $x\in \left[x_c^{-} + \frac{s_0}{2}T^{-2/3},0\right]$ and any $T>T_0$ it holds that:
    \begin{equation*}
        \mathrm{Tr}(\Pi_h\widetilde{K}\Pi_h)\geq Cs^{3/2}.
    \end{equation*}
    Here $s$ is seen as function of $h$ via the two equations $x=\frac{h}{T}$ and $s=\red{2}(x-x_c^-)T^\frac23$.
\end{proposition}

To prove Proposition~\ref{prop:trace}, we will deform the contours in \eqref{eq:contourIntegral} to make the asymptotic analysis simpler.
To do so we need to understand the level lines of $\mathrm{Re}(G_x(z))$ which pass through the critical points $z_c^\pm$. 
The following proposition describes the properties of these level lines, which are depicted in the left panel in Figure \ref{fig:contours1}: 

\begin{figure}
\centering
\def\svgscale{.5}
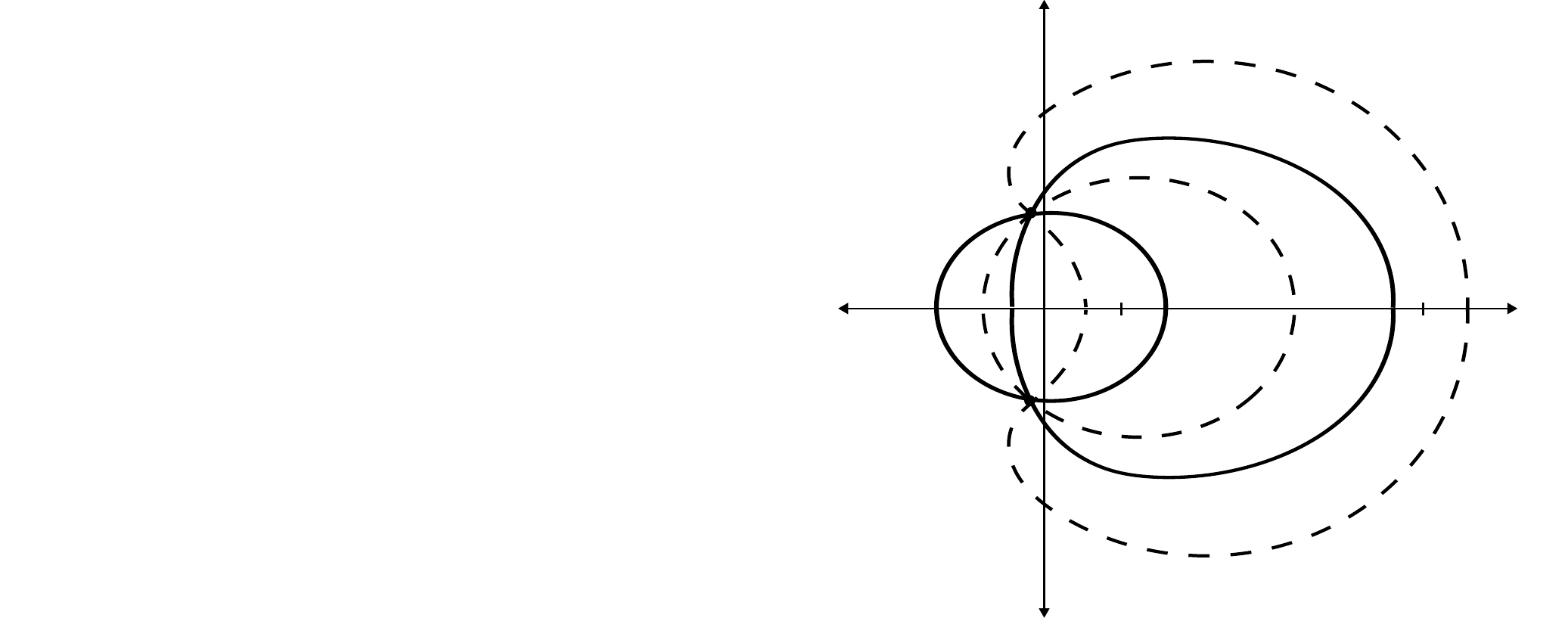
\caption{Left panel: the level lines $\mathcal{L}_1$ and $\mathcal{L}_2$. Right panel: the contours $\Gamma_1$ and $\Gamma_2$. } 
\label{fig:contours1}
\end{figure}

\begin{proposition}\label{prop:levellines}
    Let $x$ be such that $x_c^- < x < x_c^+$.
    Then there are two smooth curves $\mathcal L_1$ and $\mathcal L_2$ such that:
    \begin{enumerate}
        \item The two curves only intersect 
 at the critical points, i.e. $\mathcal{L}_1\cap\mathcal{L}_2=\{z_c^\pm\}$.
        \item The two curves are the level lines through the critical points, i.e. $\mathrm{Re}(G_x(z))=\mathrm{Re}(G_x(z_c^+))$ if and only if $z\in\mathcal{L}_1\cup\mathcal{L}_2$.
        \item Both curves are bounded simple loops.
        \item $\mathcal L_1$ contains $0$ and $1/{\sqrt{\kappa}}$ but not $\sqrt{\kappa}$, while $\mathcal L_2$ contains $\sqrt{\kappa}$ and $1/{\sqrt{\kappa}}$ but not $0$.
    \end{enumerate}
\end{proposition}
\begin{proof}
    Note that while the logarithms in the definition of $G(z)$ require a choice of branch cut, $\mathrm{Re}(\log(z))$ is defined and smooth everywhere except $0$ and does not depend on the choice of branch cut. Let us consider the level lines through the critical points.
    Since the function is critical at these points and the second derivative does not vanish, there are exactly two level lines emerging, which intersect each other at those points.
    These cannot intersect at any other points, since the intersection points would again be critical points of $G_x$.
    For $z$ with $|z|$ large it holds that
    \begin{equation}\label{eq:largeargument}
        \mathrm{Re}(G_x(z))=-(\mu+x)\ln(|z|)+\nu\ln\left(\sqrt{\kappa}\right)+o(1).
    \end{equation}
    Since for $x\geq x_c^-$ we have $\mu+x>0$, the level lines must be bounded.
    
    Considering $\mathrm{Re}(G_x(z))$ on the real line, we see poles at $0$, $1/{\sqrt{\kappa}}$, and $\sqrt{\kappa}$, where this function converges to $+$, $-$, and $+\infty$ respectively.
    Between two consecutive poles, the level lines can only cross once, since otherwise between two crossings there would be another critical point of $G$.
    This means there are exactly four points $d_1,\dots,d_4$ along the real line such that $\mathrm{Re}(G_x(d_i))=\mathrm{Re}(G(z_c^\pm))$ which satisfy $d_1<0<d_2<1/{\sqrt{\kappa}}<d_3<\sqrt{\kappa}<d_4$.
    Each of the four half-lines emanating from one of the critical points will intersect the real line at exactly one of those four points.
    Indeed the only other option would be for two of these lines to meet, but that would create a closed level-line loop containing no pole, which would force the function to be constant by harmonicity.
    A brief consideration shows that the only way to connect the half-lines gives the description in the fourth point.
    
    Finally, there cannot be any other points $z$ for which $\mathrm{Re}(G_x(z))=\mathrm{Re}(G(z_c^\pm))$ since each of those would need to lie on a closed level-line, and such a level line would need to surround a pole and therefore also intersect the real line.
    But all points on the real line with value  $\mathrm{Re}(G(z_c^\pm))$ already lie on the two level lines through the critical points.
    
\end{proof}

Using these properties of the level lines we can choose contours $\Gamma_1$ and $\Gamma_2$, as depicted in the right panel in Figure \ref{fig:contours1}:
\begin{proposition} \label{prop:contours}
    Let $x$ be such that $x_c^- < x < x_c^+$.
    Then there are two simple curves $\Gamma_1$ and $\Gamma_2$ such that:
    \begin{enumerate}
        \item The two curves only intersect 
        at the critical points, i.e. $\Gamma_1\cap\Gamma_2=\{z_c^+, z_c^-\}$.
        \item At the critical points the two curves intersect perpendicularly and in the direction of steepest ascent and descent respectively.
        \item The two curves only intersect the level lines $\mathcal L_1$ and $\mathcal L_2$ at the critical points.
        \item Both curves contain $0$ and $1/{\sqrt{\kappa}}$ but not $\sqrt{\kappa}$.
        \item On $\Gamma_1$, the function $G_x$ is always larger than $G_x(z_c^+(x))$. On $\Gamma_2$ it is always smaller.
        \item There exists an $r=r(\kappa)$ such that for $w\in\Gamma_1$ and $z\in\Gamma_2$, the inequality $|w-z|\leq r(x-x_c^-)^\frac{1}{2}$ implies that either 
        \[|w-z_c^+|<2 r(x-x_c^-)^\frac{1}{2}\text{ and }|z-z_c^+|<2 r(x-x_c^-)^\frac{1}{2}
        \]
         or 
         \[|w-z_c^-|<2 r(x-x_c^-)^\frac{1}{2}\text{ and }|z-z_c^-|<2 r(x-x_c^-)^\frac{1}{2}.
         \] 
         Furthermore $r$ can be chosen such that $4 r(x-x_c^-)^\frac{1}{2}<|z_c^+-z_c^-|$ for all $x\in[x_c^-,0]$.
    \end{enumerate}
\end{proposition}

\begin{proof}
Let us first consider the steepest descent/ascent curves through the critical points.
These are given by the level lines of $\mathrm{Im}(G_x(z))$.
They cannot cross $\mathcal L_1$ or $\mathcal L_2$ at points other than the critical points $z_c^\pm$.
Along these curves, the real part is strictly increasing/decreasing, therefore these curves must end at the poles of $G_x(z)$ which are at $0,\kappa^{-\frac12}$, and $\kappa^{\frac{1}{2}}$.
By considering the signs of the poles one can see that the steepest descent curve (which is in the region where $\mathrm{Re}(G_x)$ is positive) connects the pole at $0$ to the pole at $\kappa^{\frac{1}{2}}$.
The steepest ascent curve connects $\kappa^{-\frac12}$ to $\infty$.

By considering small circles $K_0,K_{\kappa^{-\frac12}},K_{\kappa^{\frac12}}$ around each pole and a large circle $K_\infty$ around the origin, we can construct the contours as follows:
The curve $\Gamma_1$ is given by the steepest descent curves through the critical points until those hit the circles $K_0$ and $K_{\kappa^{\frac12}}$.
Then it follows those circles such that it contains $0$ but not $\kappa^{\frac{1}{2}}$.
The curve $\Gamma_2$ is given by the steepest ascent curve until it hits $K_\infty$ and $K_{\kappa^{-\frac12}}$, where it similarly follows the circles such that it includes $0$ and $\kappa^{-\frac12}$.

By considering \eqref{eq:largeargument} and \eqref{eq:criticalpointsdistance}, one can see that the choice of circle can be made independently of $\mu$ and $x$.
Indeed one can see that $\mathrm{Re}(G(z_c))$ depends continuously on $\mu\in[\kappa^{-1}+\varepsilon,\kappa-\varepsilon]$ and $x\in[x_c^-,x_c^+]$ and is therefore bounded uniformly in absolute value, with the bound depending only on $\kappa$.
Around each of the poles, one can also find a uniform lower or upper bound depending on the sign of the pole.
For example, around $0$ one can bound:
\begin{equation}
    \mathrm{Re}(G_x(z))=\nu\ln(|\sqrt{\kappa}-z^{-1}|)-\mu\ln(|\sqrt{\kappa}-z|)-x\ln(|z|)\gtrsim (\nu+x)\ln(|z|^{-1})\,,
\end{equation}
where for $|z|$ small enough the implicit constant depends only on $\kappa$.
The prefactor $\nu+x=\frac{(\sqrt{\mu} \pm \sqrt{\kappa \nu})^2}{\kappa- 1}$ is bounded below by a constant which only depends on $\kappa$ and $\varepsilon$.
Therefore one can find a radius small enough, depending only on $\kappa$ and $\varepsilon$ such that for $K_0$, a circle of this radius, and $z\in K_0$, we have $\mathrm{Re}(G_x(z))>\mathrm{Re}(G(z_c))+1$ for all $\mu\in[\kappa^{-1}+\varepsilon,\kappa-\varepsilon]$ and $x\in[x_c^-,x_c^+]$.
In particular, this ensures that $z_c^+$ is outside this ball around $0$.
With very similar arguments one can determine the radii of $K_{\kappa^{-\frac12}},K_{\kappa^{\frac12}}$ and $K_\infty$, such that for all $x$ and $\mu$, the values of $\mathrm{Re}(G_x(z))$ on these circles is respectively larger, larger and smaller than the value of $\mathbb{Re}(G_x(z_c^+))$.

For the last point, we will actually show the following stronger statement: There exists an $r_0=r_0(\kappa)$ such that for all $r<r_0$,  $\mu\in[\kappa^{-1}+\varepsilon,\kappa-\varepsilon]$, and $x\in[x_c^-,0]$ it holds that for $w\in\Gamma_1$ and $z\in\Gamma_2$, the inequality $|w-z|\leq r$ implies that either 
        \[(|w-z_c^+|<2 r\text{ and }|z-z_c^+|<2 r)\text{ or }(|w-z_c^-|<2 r\text{ and }|z-z_c^-|<2 r)\,.
         \] 
This implies the desired statement.

To prove this, first consider fixed $x\in(x_c^-,0]$ and $\mu\in[\kappa^{-1}+\varepsilon,\kappa-\varepsilon]$. Since the two curves only intersect at $z_c^+$ and $z_c^-$ and intersect there perpendicularly, there exists an $\widetilde{r}_0=\widetilde{r}_0(x,\mu,\kappa)$ such that for all $r<\widetilde{r}_0$ the statement holds.
For $x=x_c^-$, the two critical points merge into a double critical point, and $\Gamma_1$ and $\Gamma_2$ deform in the following way.
The part of $\Gamma_1$ that connects the critical points to $K_{\kappa^\frac12}$ deforms into a piecewise continuous curve, which has a $2\pi/3$ angle at the $z_c$ and leaves this point in the directions $e^{\pi i/3}$ and $e^{-\pi i/3}$ 
The part of $\Gamma_1$ which connects the critical points to $K_0$ becomes straight lines connecting $z_c$ to $K_0$, parallel to the horizontal axis.
Similarly $\Gamma_2$ deforms into a straight line segment connecting $z_c$ to $K_{\kappa^{-1/2}}$ and a piece-wise continuous curve which goes through $z_c$ at a $2\pi/3$ angle, in the directions $e^{2\pi i/3}$ and $e^{4\pi i/3}$. 
See Figure~\ref{fig:contours2} for the level lines of a double critical point that appears in the proof of Proposition \ref{prop:otherTail} for an illustration of what this looks like.

Since these curves still only meet at $z_c$ at a $\pi/3$ angle, there is an $\widetilde{r}_0$ such that the statement holds for $x=x_c$.
\end{proof}

\begin{proof}[Proof of Proposition \ref{prop:trace}]
For clarity, we divide the proof into several steps:

\medskip

\textbf{Step 1. Decomposing the trace into two parts:} Recall that in our original definition of $\widetilde{K}$ in \eqref{eq:FredholmMeixner}, we started off with two positively oriented circular contours for $w$ and $z$ such that $1/{\sqrt{\kappa}} < |z| < 1 < |w| < \sqrt{\kappa}.$ We will now deform these two circular contours into our new choice of contours $\Gamma_1$ and $\Gamma_2$, respectively. 

\begin{figure}
\small
\centering
\def\svgscale{.34}
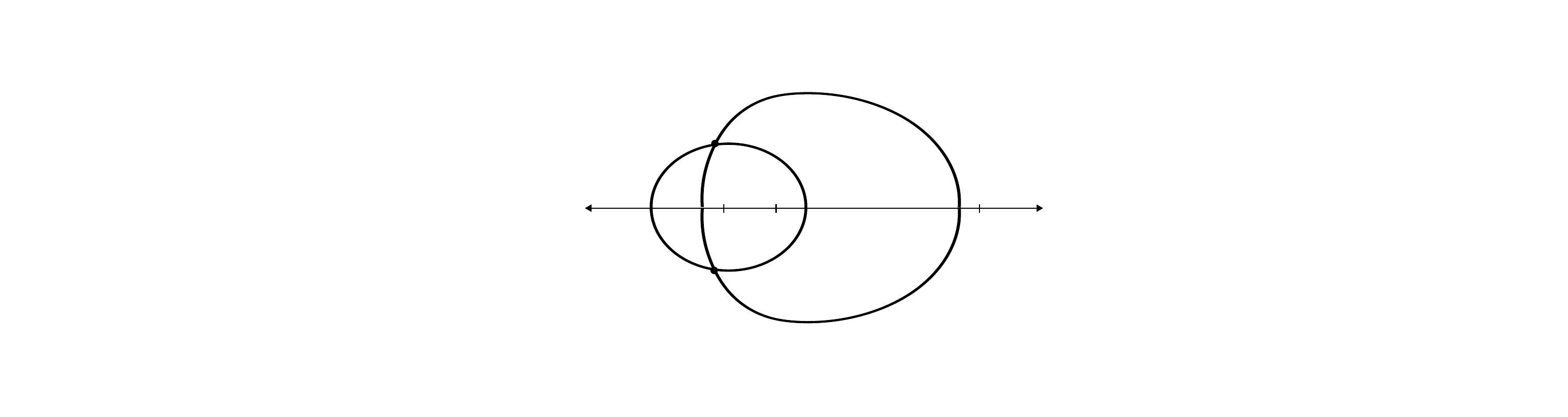
\caption{An illustration of how we pick up residues when deforming the original circular contours into $\Gamma_1$ and $\Gamma_2$.} 
\label{fig:deformation}
\end{figure}

Originally, the $z$ contour is nested inside of the $w$ contour. When we deform the $z$ contour into $\Gamma_2$, part of it will cross through the $w$  contour, see Figure \ref{fig:deformation}. Therefore we will pick up some residues since our integrand has a pole of order $2$ at $z = w$ due to the term $\frac{1}{(w-z)^2}.$ After doing this deformation, we can decompose $\text{Tr}(\Pi_h\widetilde{K}\Pi_h)$ into two parts as follows: 
\begin{multline} \label{eq:I12}
     \mathrm{Tr}(\Pi_h\widetilde{K}\Pi_h)=
   \left[ \frac{1}{(2 \pi i)^2} \oint_{\Gamma_1} \oint_{\Gamma_2} \exp (T (G_x(z) - G_x(w))\frac{dz dw}{(w-z)^2 }\right]\\
    +
   \left[ -\frac{1}{(2 \pi i)^2} \oint_\mathcal{C} \oint_{\Gamma_3} \exp (T (G_x(z) - G_x(w))\frac{dz dw}{(w-z)^2 }\right],
\end{multline}
where $\Gamma_1$ and $\Gamma_2$ are the contours given in Proposition \ref{prop:contours}, $\mathcal C$ is an arc connecting $z_c^+$ and $z_c^{-}$ and intersecting the real line between $\kappa^{-\frac12}$ and $\kappa^{\frac12}$, and $\Gamma_3$ is a contour around this arc, intersecting the real line only between $\kappa^{-\frac12}$ and $\kappa^{\frac12}$ (i.e. not including any pole other than $z=w$). 
Denote the first expression in \eqref{eq:I12} as $I_1$ and the second one as $I_2$, such that $\mathrm{Tr}(\Pi_h\widetilde{K}\Pi_h)=I_1+I_2$.
Here $I_1$ needs to be interpreted as a principal value integral due to the quadratic singularity at the intersection points. $I_2$ accounts for the residues picked up in the above-described deformation.

\textbf{Step 2: Estimating $I_2$:} In this step we show that there exists a constant $C=C(\varepsilon)$ , such that $I_2 > Cs_0^{3/2}.$ 
Let $f_w(z)  :=\frac{\exp (T(G_x(z) - G_x(w))}{(w-z)^2}$. 
Using Cauchy's residue theorem for a pole of order $2$, we can compute 
\begin{align}
    I_2= -\frac{1}{(2 \pi i)^2} \int_{\mathcal{C}} \oint_{\Gamma_3}f_w(z) dz dw 
    &=  -\frac{1}{2 \pi i} \int_{\mathcal{C}} \text{Res}(f_w, w) dw \\
    &= -\frac{1}{2 \pi i} \int_{\mathcal{C}} T G'_x(w) dw \\
    &= -\frac{T(G_x(z_c^{+}) - G_x(z_c^{-}))}{2 \pi i}. \label{eq:I2}
\end{align}

Since $z_c^{\pm}(y)$ are critical points, and since only one term in $G_y$ depends explicitly on $y$, we can compute
$$\frac{d}{dy} G_y(z_c^{\pm}(y)) = \ln (z_c^{\pm}(y)).$$
At $x_c^{-}$, we have $z_c^{+} = z_c^{-}$, so we can add and subtract $G_{x_c^{-}}(z_c^{+}) = G_{x_c^-}(z_c^{-})$ to get 

\begin{align}
    G_x(z_c^{+}) - G_x(z_c^{-}) &= (G_x(z_c^{+}) - G_{x_c^-}(z_c^{+})) - (G_x(z_c^{-}) - G_{x_c^-}(z_c^{-})) \\
    &= \int_{x_c^-}^x \ln\left(\frac{z_c^{+}(y)}{z_c^{-}(y)}\right) dy. \label{eq:I2log}
\end{align}

This integrand is purely imaginary since we are taking the log of the ratio of complex conjugates (which has modulus $1$). Therefore, when we divide by $2\pi i$ we will get something real. We now estimate the integrand. Since the modulus is $1$, the integrand is just the argument of $\frac{z_c^{+}(y)}{z_c^{-}(y)}$, which varies along the unit circle clockwise starting at $0$.

Let us define the new variable $v = \sqrt{(y - x_c^{-})}$ 
 and also define
\begin{align*}
  Z(v): = \frac{  z_c^{+}(y)}{z_c^{-}(y)} &= \frac{\mu+\nu+(\kappa + 1)y + \sqrt{(\kappa-1)^2(y - x_c^{+})(y - x_c^{-})}}{\mu+\nu+(\kappa+1)y - \sqrt{(\kappa-1)^2(y - x_c^{+})(y - x_c^{-})}} \\
  &=\frac{\mu+\nu+(\kappa+1)(v^2 + x_c^{-}) +v \sqrt{(\kappa-1)^2(v^2 + x_c^{-} - x_c^{+})}}{\mu+\nu+(\kappa+1)(v^2 + x_c^{-})  - v \sqrt{(\kappa-1)^2(v^2 + x_c^{-} - x_c^{+})}}\,.
\end{align*}

Note that  $Z(0) = 1$, so that 
\begin{align*}
    Z(v) = 1 + Z'(0)v + O(v^2)
\end{align*}
where 
$$Z'(0) = \frac{-i2(\kappa-1)^{3/2}\kappa^{1/4}(\frac{\mu}{\nu})^{1/4}}{\nu^{1/2}\left[\kappa^{1/2}(\frac{\mu}{\nu})^{1/2}(\kappa+1) - (1+\frac{\mu}{\nu})\kappa\right]}.$$
We can see that the numerator above is bounded for $\frac{\mu}{\nu} \in [\kappa^{-1}+ \varepsilon, \kappa - \varepsilon] $. 
We can also see that the denominator is zero precisely when $\frac{\mu}{\nu} \to \kappa^{-1}$ or $\frac{\mu}{\nu} \to \kappa$ and is positive between those two values. 
Finally, note that $Z'(0)$ is purely imaginary and $iZ'(0)>0$.
It follows that there exists $C = C(\varepsilon)$ such that $iZ'(0) > C$ for all $\mu  \in [\kappa^{-1}+ \varepsilon, \kappa - \varepsilon]$ and $\nu = 1$. 

It follows that 

\begin{align*}
    \ln\left(\frac{z_c^{+}(y)}{z_c^{-}(y)}\right) &= Z'(0)v + O(v^2)  \\
    &=Z'(0)\sqrt{(y - x_c^{-})} + O(y - x_c^{-})
\end{align*}
and that 
\begin{align*}
    \int_{x_c^-}^x \ln\left(\frac{z_c^{+}(y)}{z_c^{-}(y)}\right) dy = Z'(0)(x - x_c^{-})^{3/2} + O((x - x_c^{-})^2)\,,
\end{align*}
where the $O(x - x_c^{-})^2$ is uniform in $\mu$. 
Plugging this back into \eqref{eq:I2} and \eqref{eq:I2log}, we conclude that 

\begin{align}\label{eq:I2estimate}
    I_2 = -\frac{T}{2\pi i}\int_{x_c^-}^x \ln\left(\frac{z_c^{+}(y)}{z_c^{-}(y)}\right) dy&> CT(x - x_c^{-})^{3/2} \\\nonumber
    &= CT(s T^{-2/3})^{3/2} \\\nonumber
    &= Cs^{3/2}
\end{align}
for some positive real constant $C$ depending on $\varepsilon$, which changes from line to line. 

\textbf{Step 3: Estimating $I_1$.} In this section, we show that the integral $I_1$ defined above is bounded uniformly in $x$, i.e., there exists constants $C,T_0,s_0$ depending on $\varepsilon$, such that $|I_1|<C$ for all $x$ in $[x_c^-+ \frac{s_0}{\red{2}}T^{-\frac23},0]$ and $T>T_0$.

Consider first the part of the integral $I_1$ where $|z-w|\geq r(x-x_c)^\frac{1}{2}$, where $r$ is defined as in Proposition~\ref{prop:contours}:
\[
\frac{1}{(2 \pi i)^2} \oint_{\Gamma_1} \oint_{\Gamma_2} \exp (T (G_x(z) - G_x(w)))\frac{\bm{1}_{|z-w|\geq r(x-x_c)^\frac{1}{2}}dz dw}{(w-z)^2 }\,.
\]
This integral we can bound by taking absolute values and the triangle inequality to obtain
\begin{equation*}
    \frac{1}{(x-x_c)r^2(2 \pi)^2} \oint_{\Gamma_1} \oint_{\Gamma_2} \exp (T (\mathrm{Re}(G_x(z) - G_x(w)))) dz dw.
\end{equation*}
This integral has no singularities and can be split into the product of two integrals, each of which can be treated using the method of steepest descent.
Each of them gives a contribution $\frac{C}{\sqrt{TG''_x(z_c^+)}}$ for $T$ large enough.
Combined with the prefactor $\frac{1}{(x-x_c)r^2(2 \pi)^2}$, we obtain the upper bound $T^{-1}|G''_x(z_c^+)|^{-1}r^{-2}(x-x_c)^{-1}$.
By Lemma~\ref{lemma:taylor of G}, this is of order $\frac{(x-x_c)^{-\frac32}}{T}\lesssim s_0^{-3/2}$ and therefore $O(1)$.

Proposition~\ref{prop:contours} states the following.
For small but fixed $r$, all $w\in\Gamma_1$ and $z\in\Gamma_2$ such that $|w-z|\leq r(x-x_c)^\frac{1}{2}$, satisfy either 
\[|w-z_c^+|<2 r(x-x_c)^\frac{1}{2}\text{ and }|z-z_c^+|<2 r(x-x_c)^\frac{1}{2}
\]
 or 
 \[|w-z_c^-|<2 r(x-x_c)^\frac{1}{2}\text{ and }|z-z_c^-|<2 r(x-x_c)^\frac{1}{2},
 \] 
i.e. they are both close to the same critical point.
By symmetry it suffices to consider both $w$ and $z$ in the ball around $z_c^+$ of radius $2r(x-x_c)^\frac{1}{2}$.
Denote this ball by $B$.
Let us Taylor expand around our integrand.
\begin{equation}\label{eq:TaylorI1}
    \exp(T(G_x(z))-G_x(w))=\exp\left(\frac12 TG''_x(z_c^+)((z-z_c^+)^2-(w-z_c^+)^2)\right)(1+O(T|z-z_c^+|^3+|w-z_c^+|^3)),
\end{equation}
where the big $O$ constant depends on $\varepsilon$, but not on $x$ or $T$.
Let us first consider the contribution of the big $O$ term on the right.
After taking absolute values we have to bound
\[
\oint_{\Gamma_1\cap B} \oint_{\Gamma_2\cap B} T(|z-z_c^+|^3+|w-z_c^+|^3) \exp \left(\frac12 T\mathrm{Re}(G''_x(z_c^+)((z-z_c^+)^2-(w-z_c^+)^2))\right)\frac{dz dw}{|w-z|^2 }\,.
\]
Since our contours meet at a right angle at the critical point, we have $\frac1{|z-w|^2}= O\left(\frac{1}{|z-z_c^+|^2+|w-z_c^+|^2}\right)$ and the integral is bounded by
\begin{equation*}
    T\int_\R\int_\R(|z|+|w|)\exp\left(-\frac12T|G''_x(z_c^+)|(-z^2-w^2)\right)dzdw\lesssim\frac{T}{(TG''_x(z_c^+))^\frac32} =O(s_0^{-\frac34}).
\end{equation*}
For the other part of \eqref{eq:TaylorI1} we use a change of variables $z=z_c^++\xi_1/\sqrt{TG''_x(z_c^+)}$ and $w=z_c^++\xi_2/\sqrt{TG''_x(z_c^+)}$ which yields the principal value integral:
\begin{equation*}
    \int_{\widetilde{\Gamma}_1}\int_{\widetilde{\Gamma}_2} \frac{\exp(\xi_1^2-\xi_2^2)}{(\xi_1-\xi_2)^2}d\xi_1d\xi_2\,,
\end{equation*}
where $\widetilde{\Gamma}_1$ and $\widetilde{\Gamma}_2$ are contours crossing at the origin with $\widetilde{\Gamma}_1$ vertical and $\widetilde{\Gamma}_2$ horizontal there.
This is clearly bounded away from the origin. Close to the origin the exponential can be estimated by $1$ up to an $O(1)$ error, and the resulting principal value integral is also of order $O(1)$.

Combining the above estimates, one obtains that $I_1$ is $O(1)$.
Combining this with \eqref{eq:I2estimate} one obtains that
\[
\mathrm{Tr}(\Pi_h\widetilde{K}\Pi_h)=I_1+I_2>Cs^\frac32
\]
for a different constant $C$, using that $s\geq s_0$ to absorb the $O(1)$ term.
\end{proof}

We now prove Proposition \ref{prop:MeixnerTail}:

\begin{proof}[Proof of Proposition \ref{prop:MeixnerTail}]
Combining Equations \eqref{eq:comparison1}-\eqref{eq:determinantalForm} and \eqref{eq:widomtrick} and setting $\nu = 1$ and $h$ as 
\[
h=(g(\mu)-1)T+\frac{sT^{1/3}}{2}=x_c^{-}T+\frac{sT^{1/3}}{2}\,.
\]
we have that 
\begin{align}
 \mathbb{P}\left[H(T\mu, T) \geq g(\mu)T + sT^{1/3}\right] &\leq  e^{-1 /(q-1)}\left(\Pb\left[-\ell(\lambda) >\left(g(\mu) - 1\right) T +\frac{sT^{1/3}}{2}  \right] + q^{\frac{sT^{1/3}}{2}}\right)\\
 &= e^{-1 /(q-1)}\left(\operatorname{det}(\operatorname{Id}-\Pi_h\widetilde{K}\Pi_h) + q^{\frac{sT^{1/3}}{2}}\right) \\
 & \leq   e^{-1 /(q-1)}\left(\exp(-\text{Tr}(\Pi_h\widetilde{K}\Pi_h))+ q^{\frac{sT^{1/3}}{2}}\right)\,,\label{eq:estimatewithtrace}
\end{align}

By Proposition \ref{prop:trace} there exist constants $C,T_0$, and $s_0$ depending on $\varepsilon$ such that for $s>s_0, T>T_0$ and $x\in[x_c^{-} + \frac{s_0}{2}T^{-2/3},0]$:
    \begin{equation*}
        \mathrm{Tr}(\Pi_h\widetilde{K}\Pi_h)\geq Cs^{3/2}\,.
    \end{equation*}
Noting that $x$ is given by
\begin{align*}
    x = \frac{h}{T}= x_c^{-}+\frac{sT^{-2/3}}{2}\,,
\end{align*}
the restriction $x\in[x_c^{-} + \frac{s_0}{2}T^{-2/3},0]$ is equivalent to $s\in[s_0,-2T^\frac23x_c^-]$.
For $s$ in this range, applying Proposition~\ref{prop:trace} to \eqref{eq:estimatewithtrace} yields
\begin{equation}\label{eq:targetinproof}
\mathbb{P}\left[H(T\mu, T) \geq g(\mu)T + sT^{1/3}\right]\leq   e^{-1 /(q-1)}\left(\exp(-Cs^{\frac32})+ q^{\frac{sT^{1/3}}{2}}\right)\leq c^{-1}\exp(-cs^\frac32)\,.
\end{equation}
For $s>-2T^\frac23x_c^-$ we have
\[
g(\mu)T+sT^\frac{1}{3}>(x_c^-+1)T-2Tx_c^-=T(1-x_c^-)>T\,.
\]
Since $H(\mu T,T)$ can be at most $T$, for such $s$ the left-hand side of \eqref{eq:targetinproof} is $0$ and therefore \eqref{eq:tail_1} is trivially satisfied.

In summary, we have proved Proposition~\ref{prop:MeixnerTail} for all $s\geq s_0$ and all $T\geq T_0$.
For fixed $T < T_0$, the statement is trivial for $s$ large enough as the left-hand side is $0$ if $g(\mu)T+sT^\frac{1}{3}>T$.
Thus by increasing $s_0$, the statement holds for all $T\geq1$ and $s>s_0$.
By decreasing $c$ such that $c^{-1}e^{-cs}\leq1$ for $s<s_0$ the statement holds for all $T\geq 1$ and $s\geq0$.
That $c$ can be chosen weakly decreasing in $\varepsilon$ is easily checked by checking that all constants in the above estimates depend continuously on $\varepsilon$.
\end{proof}

\subsection{Proof of Proposition~\ref{prop:otherTail}}\label{sec:tail2}
The goal of this section is to prove Proposition \ref{prop:otherTail}.

\begin{remark}
Again, it suffices to prove the statement for $T \geq T_0$ and $s \geq s_0$ for some $T_0$ and $s_0$ large enough. For fixed $T$, the left-hand side becomes $0$ for $s$ large enough since the height function is always non-negative. We can also alter the constant $c$ to be small enough so that the right-hand side becomes greater than $1$ for all $s < s_0$. 
\end{remark}

In this section, we will closely follow the results in \cite{AB2019Aseps6vphase}. We start with an identity that relates the $q$-Laplace transform of the stochastic six-vertex model under step Bernoulli initial data, to a Fredholm determinant of some kernel. Recall that $(\rho,0)$-Bernoulli boundary conditions denotes the boundary condition in which the incoming arrows from the left are given by i.i.d. Bernoulli$(\rho)$ random variables, while the incoming positions from the bottom are all empty. 

We will need the following definition of a Fredholm determinant on a contour, see e.g. \cite[Definition A.1]{AB2019Aseps6vphase}
\begin{definition}[Fredholm Determinant on a Contour]
Fix a contour $\mathcal{C} \subset \mathbb{C}$ in the complex plane. Let $K: \mathcal{C} \times \mathcal{C} \rightarrow \mathbb{C}$ be a meromorphic function with no poles on $\mathcal{C} \times \mathcal{C}$. We define the Fredholm determinant
$$
\operatorname{det}(\operatorname{Id}+K)_{L^2(\mathcal{C})}=1+\sum_{k=1}^{\infty} \frac{1}{(2 \pi \mathrm{i})^k k!} \int_{\mathcal{C}} \cdots \int_{\mathcal{C}} \operatorname{det}\left[K\left(x_i, x_j\right)\right]_{i, j=1}^k \prod_{j=1}^k d x_j.
$$
\end{definition}

\begin{proposition}[Prop 5.1 in \cite{AB2019Aseps6vphase}] \label{prop:KpDef}
Fix $b_1, b_2 \in(0,1); \rho \in(0,1] ; x \in \mathbb{Z} ;$ and $p \in \mathbb{R}$. Denote $\beta=\rho /\left(1-\rho \right)$.

Let $\Gamma \subset \mathbb{C}$ be a positively oriented, star-shaped contour in the complex plane containing 0, but leaving outside $-q \kappa$ and $q \beta$. Let $\mathcal{C} \subset \mathbb{C}$ be a positively oriented, star-shaped contour contained inside $q^{-1} \Gamma$; that contains $0,-q$, and $\Gamma$; but that leaves outside $q \beta$. 

Let $\mathbb{E}_{\mathbf{6} \mathbf{v}}$ denote the expectation with respect to the stochastic six-vertex model with left jump probability $b_1$, right jump probability $b_2$, and $(\rho,0)$-Bernoulli initial data.
Then, we have that
\begin{equation}\label{eq:detIdentity}
\mathbb{E}_{\mathbf{6} \mathbf{v}}\left[\frac{1}{\left(-q^{H^{\rho}(X,T)+p} ; q\right)_{\infty}}\right]=\operatorname{det}\left(\mathrm{Id}+K^{(p)}\right)_{L^2\left(\mathcal{C}\right)}
\end{equation}
where
\begin{multline}\label{eq:kernel}
K^{(p)}\left(w, w^{\prime}\right)=\frac{1}{2 \mathbf{i} \log q} \sum_{j=-\infty}^{\infty} \oint_{\Gamma} \frac{\left(\kappa^{-1} w+q\right)^{X-1}}{\left(\kappa^{-1} v+q\right)^{X-1}}\frac{(v+q)^T}{(w+q)^T} \frac{\left(q^{-1} \beta^{-1} v ; q\right)_{\infty}}{\left(q^{-1} \beta^{-1} v ; q\right)_{\infty}}  \\ \cdot \frac{v^{p-1} w^{-p}}{\sin \left(\frac{\pi}{\log q}(\log v-\log w+2 \pi \mathbf{i} j)\right)} \frac{d v}{w^{\prime}-v}.
\end{multline}
\end{proposition}

While the above is stated for $\rho \in (0,1)$, we will ultimately need statements for $\rho = 1$, which is the case of step initial conditions. For ease of presentation, we will first obtain the tail bounds for $\rho <1$ (see Remark \ref{rmk:extendinglowertail}) and extend to $\rho = 1$ by attractivity. However, the above proposition and the estimates that follow do also extend directly to the case $\rho = 1$, see \cite[Theorem 4.16]{BCGStochasticSixVertex}. 

Once we choose $p$ appropriately, \eqref{eq:detIdentity} implies a bound on the tail probability $\mathbb{P}[ H^{\rho}(X, T) \leq  -p ]$ by applying \eqref{q1}. 
Let $\mu=\frac{X-1}{T}$ and for $\mu \in [\kappa^{-1}+\varepsilon, \kappa-\varepsilon]$  we define 
$$f_\mu=\frac{(\sqrt{\kappa \mu}-1)^{2 / 3}(\kappa-\sqrt{\kappa \mu})^{2 / 3}}{(\kappa-1) \kappa^{1 / 6} \mu^{1 / 6}}\,.$$
 The function $f_{\mu}$ appears as the scaling factor in the convergence to Tracy-Widom GUE fluctuations, see \cite[Theorem 1.2]{BCGStochasticSixVertex}. We then define for $s \geq 0$: 
$$p_T=s f_\mu T^{1 / 3}-g(\mu) T.$$ 

We will now study asymptotics of the kernel $K^{p_T}$, closely following section 6 of \cite{AB2019Aseps6vphase}, but adding in control of the decay in $s$ as well, since in \cite{AB2019Aseps6vphase} they treat $s$ as a constant. Our first step is to rewrite the formula for the kernel $K^{(p_T)}$ in an exponential form that utilizes the explicit form we chose for $p_T$. Plugging $p = p_T$ into \eqref{eq:kernel}, we obtain

\begin{align*} 
K^{(p_T)}\left(w, w^{\prime}\right)= & \frac{1}{2 \mathbf{i} \log q} \times \sum_{j \in \mathbb{Z}} \oint_{\Gamma} \frac{\exp \left(T\left(G(w)-G(v)\right)\right)}{\sin \left(\pi(\log q)^{-1}(2 \pi \mathbf{i} j+\log v-\log w)\right)} \\ & \times \frac{\left(q^{-1} \beta^{-1} v ; q\right)_{\infty}}{\left(q^{-1} \beta^{-1} w ; q\right)_{\infty}} \\ & \times\left(\frac{v}{w}\right)^{s f_{\mu} T^{1 / 3}} \frac{d v}{v\left(w^{\prime}-v\right)},
\end{align*}
where we define 

$$G(z)=\mu \log \left(\kappa^{-1} z+q\right)-\log (z+q)+g(\mu) \log z.$$

Next, we Taylor expand $G$ around its critical point: We compute its derivative 
$$G^{\prime}(z)=\frac{(\sqrt{\kappa \mu}-1)^2}{\kappa-1} \frac{\left(z-\psi \right)^2}{z(z+q \kappa)(z+q)} \quad$$ with $$\psi=\frac{q(\kappa-\sqrt{\kappa \mu})}{\sqrt{\kappa \mu}-1}.$$
Therefore, $\psi$ is a critical point of $G$, and we have $G''(\psi) = 0.$ We also have 

$$G^{\prime \prime \prime}\left(\psi\right)=\frac{2(\sqrt{\kappa \mu}-1)^5}{q^3(\kappa-1)^3(\kappa-\sqrt{\kappa \mu}) \sqrt{\kappa \mu}}=2\left(\frac{f_{\mu}}{\psi}\right)^3.$$
Putting this all together, the Taylor expansion of $G$ can be written as 
\begin{align}\label{eq:RtaylorExpansion}
G(z) & -G\left(\psi\right) =\frac{1}{3}\left(\frac{f_{\mu}\left(z-\psi\right)}{\psi}\right)^3+R\left(\frac{f_{\mu}\left(z-\psi\right)}{\psi}\right) \,,
\end{align}
where $R$ is the remainder.
By Taylor's remainder theorem, we have 
\begin{align*}
R\left(\frac{f_{\mu}\left(z-\psi\right)}{\psi}\right)&= O\left(\left|z-\psi\right|^4\right) \quad \text { as }\left|z-\psi\right| \rightarrow 0 .
\end{align*}
This remainder can be bounded uniformly for  $\kappa^{-1}+ \varepsilon  \leq \mu \leq  \kappa - \varepsilon$.
Note that $G'''(\psi)=0$ at $\mu=\kappa^{-1}$ and $G'''(\psi))=\infty$ at $\mu=\kappa^{-1}$, so $\mu$ being bounded away from $\kappa^{-1}$ and $\kappa$ is really necessary for uniformity.

\begin{figure}
\small
\centering
\def\svgscale{.5}
%% Creator: Inkscape 1.3.2 (091e20e, 2023-11-25), www.inkscape.org
%% PDF/EPS/PS + LaTeX output extension by Johan Engelen, 2010
%% Accompanies image file 'upper_tail.pdf' (pdf, eps, ps)
%%
%% To include the image in your LaTeX document, write
%%   \input{<filename>.pdf_tex}
%%  instead of
%%   \includegraphics{<filename>.pdf}
%% To scale the image, write
%%   \def\svgwidth{<desired width>}
%%   \input{<filename>.pdf_tex}
%%  instead of
%%   \includegraphics[width=<desired width>]{<filename>.pdf}
%%
%% Images with a different path to the parent latex file can
%% be accessed with the `import' package (which may need to be
%% installed) using
%%   \usepackage{import}
%% in the preamble, and then including the image with
%%   \import{<path to file>}{<filename>.pdf_tex}
%% Alternatively, one can specify
%%   \graphicspath{{<path to file>/}}
%% 
%% For more information, please see info/svg-inkscape on CTAN:
%%   http://tug.ctan.org/tex-archive/info/svg-inkscape
%%
\begingroup%
  \makeatletter%
  \providecommand\color[2][]{%
    \errmessage{(Inkscape) Color is used for the text in Inkscape, but the package 'color.sty' is not loaded}%
    \renewcommand\color[2][]{}%
  }%
  \providecommand\transparent[1]{%
    \errmessage{(Inkscape) Transparency is used (non-zero) for the text in Inkscape, but the package 'transparent.sty' is not loaded}%
    \renewcommand\transparent[1]{}%
  }%
  \providecommand\rotatebox[2]{#2}%
  \newcommand*\fsize{\dimexpr\f@size pt\relax}%
  \newcommand*\lineheight[1]{\fontsize{\fsize}{#1\fsize}\selectfont}%
  \ifx\svgwidth\undefined%
    \setlength{\unitlength}{399.87037959bp}%
    \ifx\svgscale\undefined%
      \relax%
    \else%
      \setlength{\unitlength}{\unitlength * \real{\svgscale}}%
    \fi%
  \else%
    \setlength{\unitlength}{\svgwidth}%
  \fi%
  \global\let\svgwidth\undefined%
  \global\let\svgscale\undefined%
  \makeatother%
  \begin{picture}(1,0.98004037)%
    \lineheight{1}%
    \setlength\tabcolsep{0pt}%
    \put(0,0){\includegraphics[width=\unitlength,page=1]{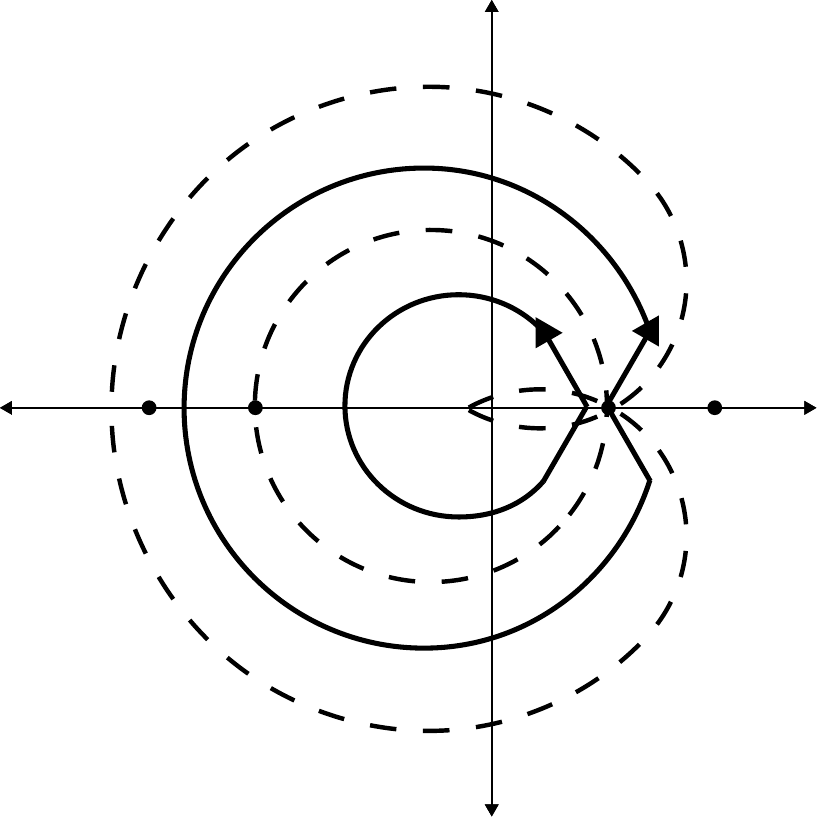}}%
    \put(0.26623823,0.42443557){\color[rgb]{0,0,0}\makebox(0,0)[lt]{\lineheight{1.25}\smash{\begin{tabular}[t]{l}$-q$\\\end{tabular}}}}%
    \put(0.13411999,0.42443557){\color[rgb]{0,0,0}\makebox(0,0)[lt]{\lineheight{1.25}\smash{\begin{tabular}[t]{l}$-q\kappa$\\\end{tabular}}}}%
    \put(0.49631767,0.55703033){\color[rgb]{0,0,0}\makebox(0,0)[lt]{\lineheight{1.25}\smash{\begin{tabular}[t]{l}$\Gamma$\\\end{tabular}}}}%
    \put(0.8491842,0.42443557){\color[rgb]{0,0,0}\makebox(0,0)[lt]{\lineheight{1.25}\smash{\begin{tabular}[t]{l}$q\beta$\\\end{tabular}}}}%
    \put(0.52181422,0.8033621){\color[rgb]{0,0,0}\makebox(0,0)[lt]{\lineheight{1.25}\smash{\begin{tabular}[t]{l}$\mathcal{C}$\\\end{tabular}}}}%
  \end{picture}%
\endgroup%

\caption{The solid lines represent the contours $\Gamma$ and $\mathcal{C}$, and the dashed curves represent the level lines of $\mathrm{Re}(G(z))$.} 
\label{fig:contours2}
\end{figure}

We now need to choose contours $\Gamma$ and $\mathcal{C}$.
We use the contours defined in \cite[Definitions 6.2-6.5]{AB2019Aseps6vphase}, which will take the following shape. The contour $\Gamma$ will consist of two parts: a piecewise linear part $\Gamma^{(1)}$ and a round part $\Gamma^{(2)}$ that connects the endpoints of $\Gamma^{(1)}$. Similarly $\mathcal{C}$ will consist of two parts $\mathcal{C}^{(1)}$ and $\mathcal{C}^{(2)}$, where $\mathcal{C}^{(1)}$ is piecewise linear and $\mathcal{C}^{(2)}$ is a round part connecting the end points of $\mathcal{C}^{(1)}$. 

\begin{definition}
For a real number $r \in \mathbb{R}$ and a positive real number $\varpi>0$ (possibly infinite), let $\mathfrak{W}_{r, \varpi}$ denote the piecewise linear curve in the complex plane that connects $r+\varpi e^{-\pi \mathbf{i} / 3}$ to $r$ to $r+\varpi e^{\pi \mathbf{i} / 3}$. Similarly, let $\mathfrak{V}_{r, \varpi}$ denote the piecewise linear curve in the complex plane that connects $r+\varpi e^{-2 \pi \mathbf{i} / 3}$ to $r$ to $r+\varpi^{2 \pi \mathbf{i} / 3}$.
    
\end{definition}

The contours $\mathcal{C}$ and $\Gamma$ look as follows:

\begin{itemize}
    \item $\mathcal{C}^{(1)}=\mathfrak{W}_{\psi,\varpi}$ and $\Gamma^{(1)}=\mathfrak{V}_{\psi-\psi f_\mu^{-1} T^{-1 / 3}, \varpi}$, for some sufficiently small $\varpi$ (independently of $T$ ).
    \item  $\mathcal{C}^{(2)}$ is a positively oriented contour from the top endpoint $\psi+\varpi e^{\pi \mathbf{i} / 3}$ of $\mathcal{C}^{(1)}$ to the bottom endpoint $\psi+\varpi e^{-\pi \mathbf{i} / 3}$ of $\mathcal{C}^{(1)}$, and  $\Gamma^{(2)}$ is a positively oriented contour from the top endpoint $\psi-\psi f_\mu^{-1} T^{-1 / 3}+$ $\varpi e^{2 \pi \mathbf{i} / 3}$ of $\Gamma^{(1)}$ to the bottom endpoint $\psi+\psi f_\mu^{-1} T^{-1 / 3}+\varpi e^{-2 \pi \mathbf{i} / 3}$ of $\Gamma^{(1)}$. 
    \item  We take $\mathcal{C}=\mathcal{C}^{(1)} \cup \mathcal{C}^{(2)}$ and $\Gamma=\Gamma^{(1)} \cup \Gamma^{(2)}$.
\end{itemize}
See Figure \ref{fig:contours2} for a depiction of these contours.

\begin{proposition}\label{prop:contourProperties}
The contours $\Gamma$ and $\mathcal C$ satisfy the following properties:
The contour $\Gamma$ is positively oriented and star-shaped; it contains 0, but leaves outside $-q \kappa$ and $q \beta$. Furthermore, $\mathcal{C}$ is a positively oriented, star-shaped contour that is contained inside $q^{-1} \Gamma$; that contains $0,-q$ and $\Gamma$; but that leaves outside $q \beta$. Furthermore, there exists some positive real number $c_1>0$, independent of $T$, such that
$$
\max \left\{\sup _{\substack{w \in \mathcal{C} \\ v \in \Gamma^{(2)}}} \mathrm{Re}(G(w)-G(v)), \sup _{\substack{w \in \mathcal{C}^{(2)} \\ v \in \Gamma}} \mathrm{Re}(G(w)-G(v))\right\}<-c_1\,,
$$
where $c_1$ depends on $\varepsilon$, but is uniform in $\mu$.
\end{proposition}

\begin{proof}
    These properties are all stated in \cite[Definition 6.3, Lemma 6.6, and Lemma 6.13]{AB2019Aseps6vphase} except for the uniformity of $c_1$. This uniformity follows from the uniformity of $G$ and $\psi$ in $\mu \in [\kappa^{-1} + \varepsilon, \kappa - \varepsilon].$ 
\end{proof}

By Proposition \ref{prop:contourProperties}, $\Gamma$ and $\mathcal{C}$ satisfy the necessary conditions stated in Proposition \ref{prop:KpDef}. We can now analyze the kernel $K: = K^{(p_T)}$ for the different cases where $v,w$ belong to the different components of these contours. The first case is where $w \in \mathcal{C}^{(1)}$ and  $v \in \Gamma^{(1)}$. In this case, both $w$ and $v$ are close to $\psi$. The second case is where either $w \in \mathcal{C}^{(2)}$ or $v \in \Gamma^{(2)}$. Let $\widetilde{K}(w, w')$ be the same kernel as $K(w, w')$, but where we replace the contour $\Gamma$ with $\Gamma^{(1)}$. 

We perform a change of variables to zoom in around $\psi$: Let $\sigma = \psi f_{\mu}^{-1}T^{-1/3}$ and set 

\begin{align*}
    &w = \psi + \sigma \widehat{w}; && w' = \sigma + \sigma \widehat{w}^{\prime}; \\
    &v = \psi + \sigma \widehat{v}; &&  \widehat{K}(\widehat{w}, \widehat{w}^{\prime}) = \sigma \widetilde{K}(w, w')\\
    &\overline{K}(\widehat{w}, \widehat{w}^{\prime}) = \sigma K(w, w').
\end{align*}

For any contour $\mathcal{D}$, set $\widehat{\mathcal{D}} = \sigma^{-1}(\mathcal{D} -\psi)$. In particular, we have 
$$\widehat{\mathcal{C}^{(1)}} = \mathfrak{W}_{0, \varpi/\sigma}, \hspace{15pt} \widehat{\Gamma^{(1)}} = \mathfrak{V}_{-1, \varpi/\sigma}.$$ 
Notice that as $T \to \infty$, we have $$\widehat{\mathcal{C}^{(1)}} \to \mathfrak{W}_{0, \infty}, \hspace{15pt} \widehat{\Gamma^{(1)}} \to  \mathfrak{V}_{-1, \infty/\sigma},$$
each of which consists of a pair of rays emanating from $0$ and $-1$ respectively. 

The only difference between $\widehat{K}$ and $\overline{K}$ is that for $\widehat{K}$ we are integrating over the contour $\widehat{\Gamma}^{(1)}$ and for $\overline{K}$ we are integrating over the contour $\widehat{\Gamma}$. The following lemma will deal with estimating these two kernels. 

\begin{lemma}
\label{lem:kernel_estimates}
There exist positive constants $c, C$ and $T_0$ all depending on $\varepsilon$ such that for $T \geq T_0$ we have

\begin{enumerate}
    \item \begin{equation*}
    \left|\overline{K}(\widehat{w}, \widehat{w}^{\prime}) - \widehat{K}(\widehat{w}, \widehat{w}^{\prime})\right| < c^{-1}\exp\left(-c\left(T+ |\widehat{w}|^3\right)\right)
\end{equation*}
for all $\widehat{w} \in \widehat{\mathcal{C}^{(1)}}$ and $\widehat{w}^{\prime} \in \widehat{\mathcal{C}} \cup  \mathfrak{W}_{0, \infty}$

    \item \begin{equation*}
    \left|\overline{K}(\widehat{w}, \widehat{w}^{\prime})\right| < c^{-1}\exp\left(-c\left(T+ |\widehat{w}|^3\right)\right)
\end{equation*}
for all $\widehat{w} \in \widehat{\mathcal{C}^{(2)}} $ and $\widehat{w}^{\prime} \in \widehat{\mathcal{C}} \cup  \mathfrak{W}_{0, \infty}$.

\item \label{item:k_hat} $$ |\widehat{K}(\widehat{w}, \widehat{w}^{\prime})| \leq \frac{C}{1 + |\widehat{w}^{\prime}|}\exp(-c |\hat{w}|^3 - cs)$$ for all $\widehat{w} \in \widehat{\mathcal{C}^{(1)}}$ and $\widehat{w}^{\prime} \in \widehat{\mathcal{C}}.$ 
\end{enumerate}

Combining the above items, we can conclude that 
\begin{equation}
 \left|\overline{K}(\widehat{w}, \widehat{w}^{\prime})\right| \leq c^{-1}\exp\left(-c|\widehat{w}|^3 - cs\right) + c^{-1}\exp\left(-c|\widehat{w}|^3 - cT\right)
\end{equation}

for all $\widehat{w} \in \widehat{\mathcal{C}}$ and $\widehat{w}^{\prime} \in \widehat{\mathcal{C}}.$ 
\end{lemma}

\begin{proof}
The proof of the first two items is the content of Corollary 6.14 in \cite{AB2019Aseps6vphase}. For the proof of Item \ref{item:k_hat} we write out the formula for kernel $\widehat{K}(\widehat{w}, \widehat{w}^{\prime})$ in terms of the variables $\widehat{w}, \widehat{w}^{\prime}$, and $\widehat{v}$, using the fact that $\sigma = \psi f_{\mu}^{-1}T^{-1/3}$ as well as the Taylor estimates in \eqref{eq:RtaylorExpansion}. We have

\begin{align*}
\widehat{K}(\widehat{w}, \widehat{w}^{\prime}) &= \sigma \widetilde{K}(w, w') \\
%%%%%%%%%%%%%%%%%%%%%%%%%%%%%%%%%%%%%%%%%%%%%
&=  \frac{\sigma}{2 \mathbf{i} \log q} \times \sum_{j \in \mathbb{Z}} \oint_{\widehat{\Gamma^{(1)}}} \frac{\exp \left(T\left(G(\psi + \sigma \widehat{w})-G(\psi + \sigma \widehat{v})\right)\right)}{\sin \left(\pi(\log q)^{-1}(2 \pi \mathbf{i} j+\log (\psi + \sigma \widehat{v})-\log (\psi + \sigma \widehat{w}))\right)} \\ 
%%%%
& \times \frac{\left(q^{-1} \beta^{-1} (\psi + \sigma \widehat{v}) ; q\right)_{\infty}}{\left(q^{-1} \beta^{-1} (\psi + \sigma \widehat{w}) ; q\right)_{\infty}} \\ 
%%%%%
& \times\left(\frac{\psi + \sigma \widehat{v}}{\psi + \sigma \widehat{w}}\right)^{s f_{\mu} T^{1 / 3}} \frac{\sigma d \widehat{v}}{(\psi + \sigma \widehat{v})\left( \sigma \widehat{w}^{\prime}- \sigma \widehat{v}\right)}\\
%%%%%%%%%%%%%%%%%%%%%%%%%%%%%%%%%%%%%%%%%%%%%%
&=  \oint_{\widehat{\Gamma^{(1)}}} \frac{\sigma \psi^{-1}}{2 \mathbf{i} \log q} \times \sum_{j \in \mathbb{Z}} \frac{\exp \left(\frac{\widehat{w}^3-\widehat{v}^3}{3}+T\left(R\left(T^{-1 / 3} \widehat{w}\right)-R\left(T^{-1 / 3} \widehat{v}\right)\right)\right)}{\sin \left(\pi(\log q)^{-1}(2 \pi \mathbf{i} j+\log (1 + \sigma \psi^{-1} \widehat{v})-\log (1 + \sigma \psi^{-1} \widehat{w}))\right)} \\
%%%%%
& \times \frac{\left(q^{-1} \beta^{-1} (\psi + \sigma \widehat{v}) ; q\right)_{\infty}}{\left(q^{-1} \beta^{-1} (\psi + \sigma \widehat{w}) ; q\right)_{\infty}} \\ 
%%%%%%
& \times\left(\frac{1 + \sigma \psi^{-1} \widehat{v}}{1 + \sigma \psi^{-1} \widehat{w}}\right)^{s \psi \sigma^{-1}} \frac{ d \widehat{v}}{(1 + \sigma\psi^{-1} \widehat{v})\left( \widehat{w}^{\prime}-\widehat{v}\right)}\,.
\end{align*}

Next, we estimate each of the terms in the integrand, using Proposition \ref{prop:integrandBounds}. Multiplying the seven bounds in Proposition  \ref{prop:integrandBounds} together, we obtain that the integrand is bounded in absolute value by $\frac{c^{-1}}{1+ |\widehat{w}^{\prime}|}\exp\left(cs\text{Re} \widehat{v} - c(|\widehat{w}|^3 + |\widehat{v}|^3)\right).$ Noting that $\text{Re} \widehat{v} \leq -1$ for $\widehat{v} \in \widehat{\Gamma^{(1)}}$, we obtain 

\begin{align*}
    |\widehat{K}(\widehat{w}, \widehat{w}^{\prime})| & \leq  \oint_{\widehat{\Gamma^{(1)}}} \frac{c^{-1}}{1+ |\widehat{w}^{\prime}|}\exp\left(cs\text{Re} \widehat{v} - c(|\widehat{w}|^3 + |\widehat{v}|^3)\right) d\widehat{v} \\
    & \leq  \oint_{\widehat{\Gamma^{(1)}}} \frac{c^{-1}}{1+ |\widehat{w}^{\prime}|}\exp\left(-cs - c(|\widehat{w}|^3 + |\widehat{v}|^3)\right) d\widehat{v} \\
    &\leq \frac{C}{1 + |\widehat{w}^{\prime}|}\exp(-c |\hat{w}|^3 - cs).
\end{align*}

\end{proof}

\begin{proposition}\label{prop:integrandBounds}
In this proposition, we prove estimates for the terms in the integrand of the kernel $\widehat{K}$. There exists $c = c(\varepsilon) > 0$ such that for $T$ large enough, the following seven bounds hold for all $\widehat{w} \in \widehat{\mathcal{C}^{(1)}} = \mathfrak{W}_{0, \varpi/\sigma} , \widehat{w}^{\prime} \in \widehat{\mathcal{C}}$, and $\widehat{v} \in \widehat{\Gamma^{(1)}} = \mathfrak{V}_{-1, \varpi/\sigma}$: 

\begin{align}
&\left|\frac{1}{1 + \sigma\psi^{-1} \widehat{v}}\right| \leq c^{-1}; \\
&\left|\frac{1}{\widehat{w}^{\prime}-\widehat{v}}\right| \leq \frac{c^{-1}}{1+ |\widehat{w}^{\prime}|}; \\
&\left|\frac{1 + \sigma \psi^{-1} \widehat{v}}{1 + \sigma \psi^{-1} \widehat{w}}\right|^{s \psi \sigma^{-1}} \leq c^{-1}\exp\left(cs \mathrm{Re} (\widehat{v})\right); \label{eq:ineq3}\\
%j neq 0
&\frac{\sigma \psi^{-1}}{|\log q|} \sum_{j  \neq 0} \left| \frac{1}{\sin \left(\pi(\log q)^{-1}(2 \pi \mathbf{i} j+\log (1 + \sigma \psi^{-1} \widehat{v})-\log (1 + \sigma \psi^{-1} \widehat{w}))\right)}\right| \leq c^{-1}T^{-1/3};\\
%j = 0
&\left | \frac{\sigma \psi^{-1}}{\log q \sin \left(\pi(\log q)^{-1}(\log (1 + \sigma \psi^{-1} \widehat{v})-\log (1 + \sigma \psi^{-1} \widehat{w}))\right)}\right| \leq c^{-1};\\
&\left|\frac{\left(q^{-1} \beta^{-1} (\psi + \sigma \widehat{v}) ; q\right)_{\infty}}{\left(q^{-1} \beta^{-1} (\psi + \sigma \widehat{w}) ; q\right)_{\infty}}\right| \leq c^{-1}\exp\left(c^{-1}\left(|\widehat{w}| + |\widehat{v}|\right)\right) \label{eq:beta}\\
&\left|\exp \left(\frac{\widehat{w}^3-\widehat{v}^3}{3}+T\left(R\left(T^{-1 / 3} \widehat{w}\right)-R\left(T^{-1 / 3} \widehat{v}\right)\right)\right)\right| \leq c^{-1}\exp\left(-\frac15(|\widehat{w}|^3 + |\widehat{v}|^3)\right)
\end{align}
\end{proposition}

\begin{proof}
The proof of all of these inequalities except for \eqref{eq:ineq3} can be found in \cite[Proof of Lemma 6.12]{AB2019Aseps6vphase} without uniformity in $\mu$. Uniformity in $\mu$ is checked in \cite[Proof of Lemma C.9]{ACG2023asepspeed}. In \cite{AB2019Aseps6vphase}, $s$ is fixed, and therefore they do not need estimates that depend on $s$. The proof of \eqref{eq:ineq3} is as follows:

First note that $|1 + \sigma \psi^{-1} \widehat{w}| \geq 1$  for all  $\widehat{w} \in \widehat{\mathcal{C}^{(1)}}$, so that we have 
\begin{align*}\left|\frac{1 + \sigma \psi^{-1} \widehat{v}}{1 + \sigma \psi^{-1} \widehat{w}}\right|^{s \psi \sigma^{-1}} \leq \left|1 + \sigma \psi^{-1} \widehat{v}\right|^{s \psi \sigma^{-1}}.
\end{align*}
For fixed $\widehat{v}$, we have 
$\lim_{T \to \infty} \left|1 + \sigma \psi^{-1} \widehat{v}\right|^{s \psi \sigma^{-1}}= \exp\left(s \text{Re}\widehat{v}\right).$ Without loss of generality let us suppose that we choose $\widehat{v}$ to be on the upper half plane so that $\widehat{v} = -1 + \frac{r\varpi}{\sigma}e^{2 \pi i/3} = -1 - \frac{r\varpi}{2\sigma} +  \frac{r\sqrt{3}\varpi i }{2\sigma}$ for some $ 0 \leq r \leq 1$. Then we have 
\begin{align*}
1 + \sigma \psi^{-1} \widehat{v} &= 1 - \sigma \psi^{-1} -\frac{r\varpi \psi^{-1}}{2} +\frac{r\sqrt{3}\varpi \psi^{-1} i }{2}.
\end{align*}
If we choose $\varpi$ small enough then there exists a $T_0$ such that for $T \geq T_0$, $\left|1 + \sigma \psi^{-1} \widehat{v}\right| < 1$. Therefore, $ \left|1 + \sigma \psi^{-1} \widehat{v}\right|^{s \psi \sigma^{-1}}$ is increasing monotonically in $T$ for $T \geq T_0$ and is therefore bounded by $\exp\left(s \text{Re}\widehat{v}\right)$. 
\end{proof}

We will need the following lemma from \cite{AB2019Aseps6vphase}: 

\begin{lemma}[Lemma A.4 in \cite{AB2019Aseps6vphase}] \label{lem:FredhomEstimate} We have 
\begin{equation}
\left|\operatorname{det}(\operatorname{Id}+K)_{L^2(\mathcal{C})}-1\right| \leq \sum_{k=1}^{\infty} \frac{2^k k^{k / 2}}{(k-1) !} \int_{\mathcal{C}} \cdots \int_{\mathcal{C}} \prod_{i=1}^k\left|\frac{1}{k} \sum_{j=1}^k| K\left(x_i, x_j\right)|^2\right|^{1 / 2} \prod_{i=1}^k d x_i .
\end{equation}
\end{lemma}

\begin{proof}[Proof of Proposition \ref{prop:otherTail}]

Combining Lemma \ref{lem:FredhomEstimate} with Lemma \ref{lem:kernel_estimates}, we obtain the following (allowing the constant $c$ to change between lines):
\begin{align*}
    &\left|\operatorname{det}(\operatorname{Id}+K^{(p_T)})_{L^2(\mathcal{C})}-1\right|  =\left|\operatorname{det}(\operatorname{Id}+\overline{K})_{L^2(\widehat{\mathcal{C}})}-1\right| \\
    & \leq \sum_{k=1}^{\infty} \frac{2^k k^{k / 2}}{(k-1) !} \int_{\widehat{\mathcal{C}}} \cdots \int_{\widehat{\mathcal{C}}} \prod_{i=1}^k\left|\frac{1}{k} \sum_{j=1}^k| \overline{K}\left(x_i, x_j\right)|^2\right|^{1 / 2} \prod_{i=1}^k d x_i \\
    & \leq \sum_{k=1}^{\infty} \frac{2^k k^{k / 2}}{(k-1) !} \int_{\widehat{\mathcal{C}}} \cdots \int_{\widehat{\mathcal{C}}} \prod_{i=1}^k\left|\frac{1}{k} \sum_{j=1}^k\left(c^{-1}\exp\left(-c|x_i|^3 - cs\right) + c^{-1}\exp\left(-c|x_i|^3 - cT\right)\right)^2\right|^{1 / 2} \prod_{i=1}^k d x_i \\
    & = \sum_{k=1}^{\infty} \frac{2^k k^{k / 2}}{(k-1) !}\left(\int_{\widehat{\mathcal{C}}} \left(c^{-1}\exp\left(-c|x|^3 - cs\right) + c^{-1}\exp\left(-c|x|^3 - cT\right)\right)d x \right)^k\\
    &= \sum_{k=1}^{\infty} \frac{2^k k^{k / 2}}{(k-1) !}\left(c^{-1}e^{-cs} + c^{-1}e^{-cT}\right)^k \left(\int_{\widehat{\mathcal{C}}} c^{-1}\exp\left(-c|x|^3\right)d x \right)^k.
\end{align*}

This last integral is bounded above by a constant, and we can also bound $(k-1)! \geq 8^{-k}k^k.$ We then obtain 
\begin{align} \label{eq:detBound}
    \left|\operatorname{det}(\operatorname{Id}+K^{(p_T)})_{L^2(\mathcal{C})}-1\right| & \leq \sum_{k=1}^{\infty} \frac{16^k}{ k^{k / 2}}\left(c^{-1}e^{-cs} + c^{-1}e^{-cT}\right)^k  \leq c^{-1}(e^{-cs} + e^{-cT}).
\end{align}

Combining \eqref{q1}, \eqref{eq:detIdentity}, and \eqref{eq:detBound} we have the following: There exists $c = c(\varepsilon)$ such that for all $T$ large enough and for all $s \geq 0$ we have that for $\mu \in \left[\frac{\kappa}{(\kappa\rho-\rho+1)^2}+\varepsilon,\kappa-\varepsilon\right]$:
    \begin{align*}
    \Pb[H^{\rho}(T\mu, T) \leq g(\mu)T - sf_{\mu}T^{1/3}] &\leq 2\left ( 1 -\operatorname{det}(\operatorname{Id}+K^{(p_T)})_{L^2(\mathcal{C})} \right)\\
    & \leq c^{-1}(e^{-cs} + e^{-cT}).
\end{align*}
Due to the fact that $f_{\mu}$ is bounded uniformly for $\mu \in [\kappa^{-1} + \varepsilon, \kappa  - \varepsilon]$, we can absorb the constant $f_{\mu}$ into $s$ by simply substituting $s\to f_\mu^{-1}s$. We then get the desired bound on $\mathbb{P}[H^{\rho}(T \mu, T)\leq g(\mu) T - sT^{1/3}  ]$ for all $\rho \in (0,1)$. Finally, we use attractivity to obtain the same bound for step initial data (i.e. $\rho =1$). Chose $\rho$ such that $[\kappa^{-1} + \varepsilon, \kappa  - \varepsilon] \subset \left[\frac{\kappa}{(\kappa\rho-\rho+1)^2}+\frac{\varepsilon}{2},\kappa-\frac{\varepsilon}{2}\right]$. Then we have that for all $\mu \in [\kappa^{-1} + \varepsilon, \kappa  - \varepsilon]$:
$$\mathbb{P}[H(T \mu, T)\leq g(\mu) T - sT^{1/3}  ] \leq \mathbb{P}[H^{\rho}(T \mu, T)\leq g(\mu) T - sT^{1/3}  ] \leq  c^{-1}(e^{-cs} + e^{-cT}).$$

Finally, the term $e^{-cT}$ can be absorbed into the term $e^{-cs}$ because the left-hand side becomes $0$ if $s > \max_{\mu \in [\kappa^{-1} + \varepsilon, \kappa - \varepsilon]}g(\mu)T^{2/3}$. We can choose $c$ to be weakly decreasing in $\varepsilon$ since all constants in the above estimates depend continuously on $\varepsilon$.
\end{proof}

\appendix

\section{Hydrodynamic limit and weak convergence}\label{ap:weak}

In this section we will give a summary of the hydrodynamic limit and local statistics for the stochastic six-vertex model proved in \cite{Aggarwal2020}, specialized to the step initial conditions on the corner.
We will then use these results to show the weak convergence of $\frac{\bm{X}_t}{t}$, which mirrors the arguments for ASEP from \cite{FerrariKipnis1995}.

We will work with the single-class stochastic six-vertex model on the quadrant, i.e. as process $(\eta_t(x))_{x,t\in\Z_{\geq0}}$ taking values in $\{0,1\}$. To state the hydrodynamic limit, we define the function $\varphi$ as
\[
\varphi(\rho):=\frac{\kappa\rho}{(\kappa-1)\rho+1}\,,
\]
where we recall that $\kappa=\frac{1-b_1}{1-b_2}$.
This function encodes the ``slope relation'' of the stochastic six-vertex model (see \cite{aggarwal2022nonexistence}).
In particular, if the stochastic six-vertex model is run from i.i.d. Bernoulli$(\rho)$ initial conditions on the bottom and i.i.d. Bernoulli$(\phi(\rho))$ from the left, the process is stationary under space-time shifts.
Furthermore, the asymptotic speed of a single second-class particle added to such initial conditions will be given by $\varphi'(\rho)$.

We can now state the general hydrodynamic limit in \cite[Theorem 1.1]{Aggarwal2020}.
This is stated for the stochastic six-vertex model on the Torus.
Let $\mathbb T^N$ be the discrete torus $\Z/N\Z$ and $\mathbb T$ the torus $\R/\Z$.
\begin{theorem}[Theorem 1.1 from \cite{Aggarwal2020}]\label{thm:apphydrolimit}
    Consider initial conditions $\eta^N_0(x):\mathbb T^N\to\{0,1\}$ which approximate a profile $\rho_0:\mathbb T\to[0,1]$, in the sense that
    \[
    \lim_{N\to\infty}\sup_{x,y}\left|\frac{1}{N}\sum_{i=\lfloor Nx\rfloor }^{\lfloor Ny\rfloor}\eta^N_0(i)-\int_x^y\rho_0(x)\right|=0\,.
    \]
    Let $(\rho_t(x))_{t\geq0,x\in\R}$ be the entropy solution of the partial differential equation
    \begin{equation}\label{eq:appPDE}
    \frac\partial{\partial_t}\rho_t(x) +\frac\partial{\partial_x}\varphi(\rho_t(x))=0\,
    \end{equation}
    with initial condition given by $\rho_0$.
    Then uniformly on compact sets, we have the following convergence in probability: 
    \[
    \lim_{N \to \infty}\frac{1}{N^2}\sum_{t = 0}^{NT} \sum_{x = 0}^{NX} \eta^N_t(x) = \int_0^T\int_0^X \rho_t(x) dx dt\,.
    \]
\end{theorem}
This theorem gives the expected density of particles at large times and scales.
It is supplemented by the following local statistics result, given in \cite[Theorem 1.3]{Aggarwal2020}, which states that whenever $\rho_t(x)$ is continuous at $(x,t)$, the system will be approximately at equilibrium. In other words, the microscopic behavior around $\eta_{Nt}(Nx)$ will approach a stationary measure.
\begin{theorem}[Theorem 1.3 from \cite{Aggarwal2020}]\label{thm:appLocal}
    In the setting of Theorem~\ref{thm:apphydrolimit}, consider a point $(x,t)$ such that $\rho_t(x)$ is continuous at $(x,t)$, and fix an integer $k\geq1$.
    Then the law of
    \[
    [\eta^N_{\lfloor Nt\rfloor+s}(\lfloor Nx\rfloor +y)]_{y\in\llbracket-k,k\rrbracket,s\in\llbracket0,k\rrbracket}
    \]
    converges in law to the stationary process started from i.i.d. Bernoulli$(\rho_t(x))$ random variables and restricted to the rectangle $\llbracket-k,k\rrbracket\times\llbracket0,k\rrbracket$. 
\end{theorem}

To prove the weak convergence we will now need the following consequence of these theorems.
\begin{corollary}\label{cor:weak}
    Let $\eta_t:\Z\to\{0,1\}$ be the stochastic six-vertex process on the line started from $\eta_0(x)=\bm{1}_{x<0}$, and let $\alpha$ be a positive real number.
    Then 
    \[
    \lim_{t\to\infty}\mathbb P[\eta_t(\alpha t)=1]=\rho_1(\alpha)\,,
    \]
    where
    \begin{equation}\label{eq:appSolution}
    \rho_t(x)=
    \begin{cases} 
        1 &\qquad \text{if $\frac{x}{t} < \kappa^{-1}$ }\\ 
        \frac{\sqrt{\kappa t/x} -1}{\kappa -1} &\qquad \text{if $\kappa^{-1} \leq \frac{x}{t} \leq\kappa  $}\\ 
        0 &\qquad \text{if $\frac{x}{t} >  \kappa$ }
    \end{cases}
    \end{equation}
    is the unique weak solution to \eqref{eq:appPDE}
    from $\rho_0(x)=\bm{1}_{x<0}$.
\end{corollary}
\begin{proof}
    Since Theorems~\ref{thm:apphydrolimit} and \ref{thm:appLocal} are stated on the torus, we need to connect the process on the line with the process on a torus.
    Let $B$ be some integer, which will later be chosen to be large.
    Let $\eta_t$ be the stochastic six-vertex model on the line started from step initial conditions and let $\eta^{N}_t:\mathbb T^{2B^3N+1}\to\{0,1\}$ be the stochastic six-vertex process on $\mathbb T^{2B^3N+1}$ started from the initial conditions
    \[
    \eta^{N,\mathbb T}_t(x)=
    \begin{cases}
        1\text{, for } -B^2N\leq x\leq 0\\
        0\text{ else,}
    \end{cases}
    \]
    where we identify the torus $\mathbb T^{2B^3N+1}$ with the set $\llbracket-B^3N,B^3N\rrbracket$.
    Then by \cite[Proposition 5.7]{Aggarwal2020} the processes $\eta$ and $\eta^N$ can be coupled to agree on the interval $\llbracket-BN,BN\rrbracket$ until time $BN$ with probability $1-\exp{-BN}$, if $B$ is large enough such that $B^2N-\frac{4BN}{1-b_2}>BN$.
    
    Assuming further that $B$ is large than $\alpha$, we obtain
    \[
    |\mathbb P[\eta_N(\alpha N)]-\mathbb P[\eta^{N,\mathbb T}_N(\alpha N)]|\leq \exp{-BN}\,.
    \]
    Therefore it remains to apply Theorem~\ref{thm:appLocal} to $\eta^N$. 
    Note that the PDE \eqref{eq:appPDE} is invariant under scaling space and time by the same factor, and we can therefore also consider $\rho_t$ defined on the torus of size $2B^3$, which simplifies the notation. 
    The initial condition of $\eta^N$ approximate initial conditions $\rho^{\mathbb T}_0(x)=\bm{1}_{x\in [-B^2,0]}$ on the torus of side-length $2B^3$ identified with $[-B^3,B^3]$.
    By Proposition 5.3 and Remark 5.4 of \cite{Aggarwal2020}, the solution $\rho^{\mathbb T}_t$ of \eqref{eq:appPDE} for these initial conditions and $\rho_t$ (the solution for step initial conditions on $\mathbb R$) agree at time $t$ on $[-B^2+ct,B^2-ct]$, where $c=\max_{x\in[\kappa^{-1},\kappa]}|\phi'(x)|$.
    Given $B$ large enough such that $B^2-c>\alpha$, this implies that
    \[
    \rho^{\mathbb T}_1(\alpha)=\rho_1(\alpha)\,.
    \]
    Applying Theorem~\ref{thm:appLocal} with $k=0,t=1$ and $x=\alpha$ gives the desired result.
\end{proof}

\begin{proposition}[Weak Convergence of the Speed of the Second-Class Particle]\label{prop:weak}
    Let $\bm{X}_t$ be the position of the second-class particle under step initial conditions with a single second-class particle at the origin, as in Theorem~\ref{thm:main} and let $\rho_t(x)$ be given by \eqref{eq:appSolution}.
    Then the asymptotic speed of the second-class particle $\frac{\bm{X}_t}{t}$ converges weakly to a random variable with density
    \begin{equation}
        \frac{\sqrt{\kappa}}{2(\kappa-1)}x^{-\tfrac32}\bm{1}_{\kappa^{-1}\leq x\leq \kappa}\,.
    \end{equation}
\end{proposition}
\begin{proof}
Consider step initial conditions (for the single-class model), i.e. $\eta_0(x) = \bm{1}_{x <  0},$
% \begin{equation}
%     \eta_0(x)=
%     \begin{cases}
%     1&\text{if }x<0\\
%     0&\text{if }x\geq 0\,.
%     \end{cases}
% \end{equation}
and let $_0$ be the same initial conditions shifted by $1$ to the right i.e. $\widetilde{\eta}_0(x)=\bm{1}_{x\leq 0}$.
Let $(\eta_t)_{t\geq0}$ and $(\widetilde{\eta}_t)_{t\geq0}$ be the two stochastic six-vertex processes started from these initial conditions, with the height functions made unique by the choice that $h_0(1;\eta)=h_0(1;\widetilde{\eta})=0$.
In particular $h_0(0;\eta)=0$, but $h_0(0;\widetilde{\eta})=1$.
We will couple them in two different ways.
The first coupling $\pi_1$ is given by the multi-class stochastic six-vertex model with step initial conditions and a single second-class particle at the origin, i.e. the setup of Theorem~\ref{thm:main} and of this proposition.
Denote the position of the second-class particle at time $t$ with $\bm{X}_t$, as above.
The second coupling $\pi_2$ is given by the deterministic shift i.e.
\[
\mathbb P_{\pi_2}[\widetilde{\eta}_t(x)=\eta_t(x-1)]=1.
\]

We will now calculate $\mathbb{E}[h_t(x;\widetilde{\eta})-h_t(x;\eta)]$ under both of these couplings.
Under $\pi_1$ the configurations are identical except for the second-class particle and so the height functions agree for all $(x,t)$ such that $x>\bm{X}_t$ and $h_t(x;\eta)=h_t(x;\widetilde{\eta})-1$ for $(x,t)$ such that $x\leq\bm{X}_t$.
This gives
\[
\mathbb{E}_{\pi_1}[h_t(x;\widetilde{\eta})-h_t(x;\eta)]=\mathbb P[x\leq\bm{X}_t]\,.
\]
Under $\pi_2$ the height functions are related by deterministic shift $h_t(x;\widetilde{\eta})=h_t(x-1;\eta)$ and therefore
\[
\mathbb{E}_{\pi_2}[h_t(x;\widetilde{\eta})-h_t(x;\eta)]=\mathbb{E}_{\pi_2}[h_t(x-1;\eta)-h_t(x;\eta)]=\mathbb E[\eta_t(x-1)]\,.
\]

By the linearity of expectations, $\mathbb{E}[h_t(x;\widetilde{\eta})-h_t(x;\eta)]$ does not depend on the coupling we take and therefore we obtain the following identity:
\[
\mathbb P[x\leq\bm{X}_t]=\mathbb E[\eta_t(x-1)]\,.
\]
Choosing $x=\lfloor\alpha t\rfloor$ we obtain
\[
\mathbb P\left[\frac{\bm{X}_t}{t}\geq\alpha\right]=\mathbb E[\eta_t(\lfloor\alpha t\rfloor-1)].
\]
By Corollary~\ref{cor:weak} the right hand side converges to $\rho_1(\alpha)$.
Therefore $\frac{\bm{X}_t}{t}$ converges weakly to the random variable with density
\[
-\rho_1'(x)=\frac{\sqrt{\kappa}}{2(\kappa-1)}x^{-\frac32}\bm{1}_{\kappa^{-1}\leq x\leq \kappa}\,,
\]
which proves the proposition.
\end{proof}

\section{Uniqueness of stationary measures}\label{ap:unique}

The goal of this section is to prove Proposition~\ref{prop:uniqueStationaryMeasures}.
The proof builds on ideas in \cite{liggett1976coupling} and \cite{ferrari1991microscopic}. First, we define some sets of measures that we will be working with:

\begin{definition}

\begin{enumerate}[label=(\alph*),leftmargin=18pt]
    \item[]
    \item Denote by $\mathcal{P}_n$ the set of measures on $\{1,2,\dots,n-1,\infty\}^\Z$.
    \item Denote by $\mathfrak M_t$ the map from $\mathcal{P}_n$to $\mathcal{P}_n$ given by the dynamics of the multi-class stochastic six-vertex model, i.e. if the initial conditions $\eta_0$ are sampled according to $\mu$, then $\mathfrak M_t$ is the law of $\eta_t$.
    % \item Denote by $\mathfrak S_k$ the shift operator on $\mathcal{P}_n$, i.e. if $\eta$ is sampled according to $\mu\in\mathcal P_n$, then $\mathfrak S_k$ is is the law of $(\eta(x+k))_{x\in\Z}$.
    \item Denote by $\mathcal{S}_n \subset \mathcal{P}_n$ the set of stationary measures with respect to the dynamics of the multi-class stochastic six-vertex model.
    \item Denote by $\mathcal{T}_n \subset \mathcal{P}_n$ the set of translation invariant measures with respect to shifts of $\mathbb{Z}$.
    \item Denote by $\mathcal{P}^{\otimes2}_n$ the set of measures on $\{1,2,\dots,n-1,\infty\}^\Z\times\{1,2,\dots,n-1,\infty\}^\Z$, i.e. each $\pi\in\mathcal P^{\otimes2}_n$ is a coupling of two measures in $\mathcal P_n$.
    Similarly, let $\mathcal T^{\otimes2}_n$ be the set of translation invariant measures in $\mathcal P^{\otimes2}_n$ and let $\mathcal S^{\otimes2}_n$ be the set of measures which are stationary under the dynamics given by the basic coupling.
\end{enumerate}
    
\end{definition}

\begin{definition}
    For a set of measures $\mathcal A$, denote by $\mathcal{A}_e$ the set of extremal measures in $\mathcal A$. In otherwords, if $\mu\in\mathcal{A}_e$, then $\mu=\lambda\mu_1+(1-\lambda)\mu_2$ for $\lambda\in[0,1]$ and $\mu_1,\mu_2\in\mathcal A$ implies that $\lambda$ is either $0$ or $1$.
    Measures in $(\mathcal T_n)_e$ are called ergodic.
\end{definition}

We will use the following two results from \cite{Aggarwal2020}, the second of which we need to modify slightly:
% First the statement for $n=1$.
\begin{theorem}[Theorem 3.6 from \cite{Aggarwal2020}]\label{thm:singleclassuniqueness}
    If $\mu\in(\mathcal{T}_1\cap\mathcal{S}_1)_e$, it is given by the law of i.i.d. Bernoulli$(\rho)$ random variables for some $\rho\in[0,1]$.
\end{theorem}

\begin{proposition}\label{prop:orderedStationaryprocesses}
    Let $\pi$ be in $\mathcal S^{\otimes2}_1 \cap \mathcal T_{1}^{\otimes 2}$ and $(\eta,\xi)$ be sampled according to $\pi$.
    Then $\mathbb P[\eta\leq\xi\text{ or } \xi\leq\eta]=1$.
\end{proposition}

\begin{proof}
    In \cite[Corollary 3.9]{Aggarwal2020} this is stated for a different coupling, namely the ``higher rank coupling".
    However, the proof generalizes to the basic coupling in a straightforward manner.
    Most of the proof only uses attractivity, the merging property, and finite speed of propagation, all of which hold for both the higher rank and the basic coupling (Lemmas \ref{prop:attractivity}, \ref{lem:merging}, and Proposition \ref{prop:finitespeed}, respectively).
    The only point where specific properties of the higher rank coupling are used is at the end of the proof of \cite[Proposition 3.7]{Aggarwal2020} to obtain the following fact.
    
    Let 
    $$R(I;\eta,\xi) :=\begin{cases}
        1 & \text{if there are $x,y\in I\subset\mathbb Z$, such that $\eta(x)=\xi(y)=1$ and $\eta(y)=\xi(x)=0$}\\
        0 & \text{otherwise}.
    \end{cases}
    $$
    % $R(I;\eta,\xi)$ be $1$ if there are $x,y\in I\subset\mathbb Z$, such that $\eta(x)=1=\xi(y)$ and $\eta(y)=0=\xi(x)$ and $0$ otherwise.
    In other words $R(I;\eta,\xi)=1$ if there are uncoupled arrows in both $\eta$ and $\xi$ and is $0$ otherwise.
    Given two coupled stochastic six-vertex processes $\eta$ and $\xi$ and an interval $I$ of length $k$, let $F(I;t)$ be the event that there are uncoupled arrows in both $\eta$ and $\xi$ in the interval $I$ at time $t-1$ which couple at time $t$.
    Then
    \begin{equation}\label{eq:arrowsCouple}
        \mathbb P[F(I;t)|R(I;\eta_{t-1},\xi_{t-1})=1]\geq c^k\,,
    \end{equation}
    for some $c=c(b_1,b_2)>0$ independent of $k,\eta$ and $\xi$.
    
    The inequality \eqref{eq:arrowsCouple} also holds for the basic coupling, and can be seen as follows:
    Let $x$ and $y$ be the positions of the two uncoupled arrows i.e. $\eta_t(x)=\xi_t(y)=1$ and $\eta_t(y)=\xi_t(x)=0$.
    Without loss of generality we assume that $x<y$ and $\eta_{t-1}(z)=\xi_{t-1}(z)$ for all $z\in[x+1,y-1]$.
    Using the Bernoulli random variables from Definition~\ref{def:basiccoupling}, the two arrows will couple if
    \begin{itemize}
        \item $\chi^2(t,x)=0$,
        \item $\chi^1(t,x)=0$,
        \item $\chi^2(t,z)=1$ for $z\in[x+1,y-1]$ and
        \item $\chi^2(t,y)=0$ and $\chi^1(t,y)=1$ (alternatively $\chi^2(t,y)=1$ and $\chi^1(t,y)=0$).
    \end{itemize}
    Indeed the first two conditions guarantee that the arrow started at $x$ will not couple with an arrow coming from the left, the second and third conditions guarantee that it will move to $y$, and the last condition guarantees that it will couple with the other uncoupled arrow at $y$.
    Since the $\chi$ are all independent, this implies that $\Pb[F(I;t)|R(I;\eta_{t-1},\xi_{t-1})=1]$ is at least $(1-b_2)^2b_2^{k-2}(1-b_1)b_1$, which clearly can be bounded below by $c^k$ for some $c$ depending only on $b_1$ and $b_2$.
\end{proof}

To prove Proposition~\ref{prop:uniqueStationaryMeasures} we will consider two types of projections from $\mathcal P_{n+1}$ to $\mathcal P_n$:
\begin{definition}
    Let $\phi:\{1,\dots,n,\infty\}\to\{1,\dots,n-1,\infty\}$ and $\psi:\{1,\dots,n,\infty\}\to\{1,\infty\}$ be given by
    \[
    \phi(x)=
    \begin{cases}
        x\text{, if } x<n\\
        \infty\text{, if } x=n,\infty
    \end{cases}
    \quad\text{and}\quad
    \psi(x)=
    \begin{cases}
        1\text{, if } x\leq n\\
        \infty\text{, if } x=\infty
    \end{cases}\,.
    \]
    For $\mu\in\mathcal P_{n+1}$ and $\eta$ distributed according to $\mu$, let $\Phi(\mu)\in\mathcal P_n$ be the law of $\phi\circ\eta$, and $\Psi(\mu)$ be the law of $\psi\circ\eta$.
    Similarly define $\Phi:\mathcal{P}^{\otimes2}_{n+1}\to\mathcal P^{\otimes2}_n$ by setting $\Phi(\pi)$ to be the law of $(\phi\circ\eta,\phi\circ\xi)$ where $(\eta,\xi)$ is sampled according to $\pi$.
\end{definition}
Both of these maps correspond to merging classes of particles.
The map $\Phi$ merges class $n$ with class $\infty$, i.e. it turns all particles of class $n$ into holes.
The map $\Psi$ merges the classes $1$ to $n$ into a single class to obtain a single-class stochastic six-vertex model.
\begin{lemma}\label{lem:Phi}
    The maps $\Phi:\mathcal P_{n+1}\to\mathcal P_n$ and $\Psi:\mathcal P_{n+1}\to \mathcal P_1$ satisfy
    \begin{align*}
        \Phi(\mathcal S_{n+1})\subset \mathcal S_n\,,\\
        \Phi(\mathcal T_{n+1})\subset \mathcal T_n\,,\\
        \red{\Phi((S_{n+1}\cap\mathcal T_{n+1})_e)\subset (S_n\cap\mathcal T_n)_e\,,}
    \end{align*}
    and the same holds for $\Psi$, with the sets on the right-hand side replaced with their single-class counterparts.
\end{lemma}
\begin{proof}
    The first statement is a consequence of the merging property, which implies that $\Phi$ commutes with the Markov kernel for the multi-class dynamics.
    The second one follows since $\Phi$ commutes with shifts.
    \red{For the third statement assume the contrary, i.e $\mu\in(S_{n+1}\cap\mathcal T_{n+1})_e$ but $\Phi(\mu)\notin(S_n\cap\mathcal T_n)_e$.
    Stationary, translation invariant measures can be identified with Gibbs measures on $(\eta_t(x))_{t,x\in\mathbb Z}$, which are translation-invariant with respect to shifts in both $t$ and $x$.
    To do this note that the law of $\eta_0$ determines the laws of $\eta_t$ for $t\geq0$, and by shifting $t$ backwards one obtains a measure on $(\eta_t(x))_{x\in\mathbb Z,t\geq -n,}$.
    These are consistent, and thus define a measure $\tilde{\mu}$ on $(\eta_t(x))_{t,x\in\mathbb Z}$.
    Stationarity and translation-invariance of $\mu$ give rise to translation-invariance of $\tilde{\mu}$ in $t$ and $x$ respectively.
    We can define $\Phi(\tilde{\mu})$ analogously to $\Phi(\mu)$.
    If $\Phi(\mu)$ is not extremal, then $\Phi(\tilde{\mu})$ is also not extremal and therefore there is an event $A$, which has $\Phi(\tilde{\mu})(A)\in(0,1)$, which is translation-invariant in both $x$ and $t$.
    However $\tilde{\mu}(\Phi^{-1}(A))=\Phi(\tilde{\mu})(A)\in(0,1)$ and therefore $\tilde{\mu}$ is not extremal either.
    This means that $\tilde\mu$ can be written as a convex combination of translation invariant Gibbs measure measures, which implies that $\mu$ is also not extremal (by restricting to $\eta_0$).
    This contradicts the assumption, and therefore $\Phi(\mu)$ is extremal.}
\end{proof}

\begin{proof}[Proof of Proposition~\ref{prop:uniqueStationaryMeasures}]
    The existence of such measures is shown in \cite[Theorem 7.6]{ANPstationarityBaxter}, see in particular Section 7.3 therein, where the specialization to the stochastic six-vertex model is discussed, and \cite[Equation (7.8)]{ANPstationarityBaxter} where the densities of vertical and horizontal arrows are calculated.
    From this equation, one can see that all possible densities $(\lambda_k)_{k\in\{0,\dots,n\}}$ satisfying $\sum_{k=0}^n\lambda_k=1$ can be obtained by correctly choosing the parameters $\alpha_i$.
    Note that one could also prove the existence of such measures abstractly by a compactness argument.
    However, the above construction allows us to calculate marginals for these measures.

    As remarked in \cite[Remark 7.9]{ANPstationarityBaxter}, the uniqueness of such measures does not seem to be present in the literature.
    We will now prove this by induction on $n$.
    For $n=1$ this is a consequence of \cite[Theorem 3.6]{Aggarwal2020} (restated above as Theorem \ref{thm:singleclassuniqueness}), which states that the only stationary ergodic measures for the single-class stochastic six-vertex process are the i.i.d. Bernoulli$(\rho)$ measures.
    
    Let us assume that the statement has already been proven for some $n\geq 1$.
    Consider two measures $\mu,\nu\in(\mathcal{T}_{n+1}\cap\mathcal{S}_{n+1})_e$, which satisfy
    \begin{equation}\label{eq:densities}
    \Pb_\mu[\eta(0)=k]=\Pb_\nu[\eta(0)=k]\text{ for }k=1,\dots,n\,.
     \end{equation}
    By Lemma~\ref{lem:Phi} both $\Phi(\mu)$ and $\Phi(\nu)$ are in $(\mathcal S_n\cap\mathcal T_n)_e$.
    By \eqref{eq:densities} they have the same densities and hence by the induction hypothesis, they are equal.

    We will now construct a coupling $\pi\in(\mathcal S^{\otimes2}_{n+1}\cap\mathcal T^{\otimes2}_{n+1})_e$ of the measures $\mu$ and $\nu$ such that
    \begin{equation}\label{eq:condition}
        \text{if $(\eta,\xi)$ are sampled according to $\pi$, then $\phi\circ\eta=\phi\circ\xi$ almost surely.}
    \end{equation}
    To do so, first let $\pi_0$ be the coupling obtained by first sampling $\phi\circ\eta=\phi\circ\xi$ according to $\Phi(\mu)=\Phi(\nu)$ and then $\eta$ according to $\mu$ conditioned on $\phi\circ\eta$ and $\xi$ according to $\nu$ conditioned on $\phi\circ\xi$.
    Since $\mu$ and $\nu$ are translation invariant, so is $\pi_0$.

    Next, let $\pi_t$ be the law given by running the basic coupling started from initial conditions $(\eta,\xi)$ sampled according to $\pi_0$.
    By the construction of the basic coupling and the stationarity of $\mu$ and $\nu$, $\pi_t$ will be translation invariant and satisfy the condition \eqref{eq:condition}.
   
    Since $\{1,2,\dots,n,\infty\}^\Z\times\{1,2,\dots,n,\infty\}^\Z$ is compact, there is a subsequence $t_n$ such that $\frac{1}{t_n}\sum_{i=0}^{t_n}\pi_i$ converges weakly to a measure $\tilde{\pi}\in\mathcal S^{\otimes2}_{n+1}\cap\mathcal T^{\otimes2}_{n+1}$, which satisfies \eqref{eq:condition}.
    The set of measures in $\mathcal S^{\otimes2}_{n+1}\cap\mathcal T^{\otimes2}_{n+1}$ that satisfy \eqref{eq:condition} is convex and by the above argument non-empty.
    Taking an extremal measure in this set we obtain the desired measure $\pi$.

    If $\pi$ were not in $(\mathcal S^{\otimes2}_{n+1}\cap T^{\otimes2}_{n+1})_e$, that would mean that there are $\mu_1\neq\mu_2\in\mathcal S^{\otimes2}_{n+1}\cap T^{\otimes2}_{n+1}$ and $\lambda\in(0,1)$ such that $\pi=\lambda\mu_1+(1-\lambda)\mu_2$.
    However since $\pi$ is extremal in the set of measures in $\mathcal S^{\otimes2}_{n+1}\cap T^{\otimes2}_{n+1}$ satisfying condition \eqref{eq:condition}, this means that either $\mu_1$ or $\mu_2$ does not satisfy \eqref{eq:condition}.
    This cannot be, since for any event $A$ the equality $\pi(A)=1$ implies $\mu_1(A)=1$ and $\mu_2(A)=1$.
    Therefore $\pi\in(\mathcal S^{\otimes2}_{n+1}\cap T^{\otimes2}_{n+1})_e$.

    Consider now $\Psi(\pi)$, which is a coupling of $\Psi(\mu)$ and $\Psi(\nu)$. Both $\Psi(\mu)$ and $\Psi(\nu)$ are in $(\mathcal S_1\cap\mathcal T_1)_e$ and have the same particle density, therefore they are both given by the law of i.i.d. Bernoulli$(\rho)$ random variables with the same $\rho$, by Theorem~\ref{thm:singleclassuniqueness}.
    Furthermore, $\Psi(\pi)$ satisfies the conditions of Proposition~\ref{prop:orderedStationaryprocesses} and therefore if $(\tilde{\eta},\tilde{\xi})$ is sampled according to $\Psi(\pi)$, either almost surely $\tilde{\eta}\leq\tilde{\xi}$ or $\tilde{\xi}\leq\tilde{\eta}$, since it is also extremal.
    Since for two Bernoulli$(\rho)$ random variables $\eta(0), \xi(0)$  $\mathbb P(\eta(0)>\xi(0))>0$ implies $\mathbb P(\eta(0)<\xi(0))>0$, this means that under $\Psi(\pi)$ it holds that $\tilde{\eta}=\tilde{\xi}$ almost surely.
    This means that for $(\eta, \xi)$ sampled according to $\pi$, $\psi\circ\eta=\psi\circ\xi$ almost surely. By the definition of $\pi$, we also have that $\phi\circ\eta=\phi\circ\xi$ almost surely. Note that any two configurations that are equal under both projections $\psi$ and $\phi$ are themselves equal. Therefore, we conclude that almost surely $\eta=\xi$, and hence $\mu=\nu$.
\end{proof}

\bibliographystyle{alpha}
\bibliography{refs.bib}

@article {AB2019Aseps6vphase,
    AUTHOR = {Aggarwal, Amol and Borodin, Alexei},
     TITLE = {Phase transitions in the {ASEP} and stochastic six-vertex
              model},
   JOURNAL = {Ann. Probab.},
  FJOURNAL = {The Annals of Probability},
    VOLUME = {47},
      YEAR = {2019},
    NUMBER = {2},
     PAGES = {613--689},
      ISSN = {0091-1798},
   MRCLASS = {60K35 (82C22)},
  MRNUMBER = {3916931},
MRREVIEWER = {Jianping Jiang},
       DOI = {10.1214/17-AOP1253},
       URL = {https://doi.org/10.1214/17-AOP1253},
}

@article{ACG2023asepspeed,
title = {The {ASEP} speed process},
journal = {Advances in Mathematics},
volume = {422},
pages = {109004},
year = {2023},
issn = {0001-8708},
doi = {https://doi.org/10.1016/j.aim.2023.109004},
url = {https://www.sciencedirect.com/science/article/psii/S0001870823001470},
author = {Amol Aggarwal and Ivan Corwin and Promit Ghosal}
}

@article {MR3649488,
    AUTHOR = {Borodin, Alexei and Olshanski, Grigori},
     TITLE = {The {ASEP} and determinantal point processes},
   JOURNAL = {Comm. Math. Phys.},
  FJOURNAL = {Communications in Mathematical Physics},
    VOLUME = {353},
      YEAR = {2017},
    NUMBER = {2},
     PAGES = {853--903},
      ISSN = {0010-3616},
   MRCLASS = {60K35 (33C45 60G55)},
  MRNUMBER = {3649488},
       DOI = {10.1007/s00220-017-2858-1},
       URL = {https://doi.org/10.1007/s00220-017-2858-1},
}

@article {MR3760963,
    AUTHOR = {Borodin, Alexei},
     TITLE = {Stochastic higher spin six vertex model and {M}acdonald
              measures},
   JOURNAL = {J. Math. Phys.},
  FJOURNAL = {Journal of Mathematical Physics},
    VOLUME = {59},
      YEAR = {2018},
    NUMBER = {2},
     PAGES = {023301, 17},
      ISSN = {0022-2488},
   MRCLASS = {60K35 (60B20 60E10 81R12)},
  MRNUMBER = {3760963},
       DOI = {10.1063/1.5000046},
       URL = {https://doi.org/10.1063/1.5000046},
}

@article{FerrariKipnis1995,
author = {Ferrari, P. A. and Kipnis, C.},
journal = {Annales de l'I.H.P. Probabilités et statistiques},
language = {eng},
number = {1},
pages = {143-154},
publisher = {Gauthier-Villars},
title = {Second class particles in the rarefaction fan},
url = {http://eudml.org/doc/77499},
volume = {31},
year = {1995},
}

@Article{Aggarwal2020,
author={Aggarwal, Amol},
title={Limit Shapes and Local Statistics for the Stochastic Six-Vertex Model},
journal={Communications in Mathematical Physics},
year={2020},
month={May},
day={01},
volume={376},
number={1},
pages={681-746},
issn={1432-0916},
doi={10.1007/s00220-019-03643-w},
url={https://doi.org/10.1007/s00220-019-03643-w}
}

@article{Aggarwal2016CurrentFO,
  title={Current fluctuations of the stationary {ASEP} and six-vertex model},
  author={Amol Aggarwal},
  journal={Duke Mathematical Journal},
  year={2016},
  volume={167},
  pages={269-384}
}

@article{AmirAngelValko2008TasepSpeed,
  title={The {TASEP} speed process.},
  author={Gideon Amir and Omer Angel and Benedek Valk\'o},
  journal={Annals of Probability},
  year={2008},
  volume={39},
  pages={1205-1242}
}

@article{BorodinBufetov2019ColorSymmetry,
  title={Color-position symmetry in interacting particle systems},
  author={Alexei Borodin and Alexey Bufetov},
  journal={The Annals of Probability},
  year={2019}
}

@article{GwaSpohnStochasticSixVertex,
  title = {Six-vertex model, roughened surfaces, and an asymmetric spin {H}amiltonian},
  author = {Gwa, Leh-Hun and Spohn, Herbert},
  journal = {Phys. Rev. Lett.},
  volume = {68},
  issue = {6},
  pages = {725--728},
  numpages = {0},
  year = {1992},
  month = {Feb},
  publisher = {American Physical Society},
  doi = {10.1103/PhysRevLett.68.725},
  url = {https://link.aps.org/doi/10.1103/PhysRevLett.68.725}
}

@misc{ANPstationarityBaxter,
      title={{Colored Interacting Particle Systems on the Ring: Stationary Measures from {Y}ang-{B}axter Equation}}, 
      author={Amol Aggarwal and Matthew Nicoletti and Leonid Petrov},
      year={2023},
      eprint={2309.11865},
      archivePrefix={arXiv},
      primaryClass={math.PR}
}

@article{ABGMtazrpspeed,
author = {Gideon Amir and Ofer Busani and Patr{\'i}cia Gon{\c{c}}alves and James B. Martin},
title = {{The {TAZRP} speed process}},
volume = {57},
journal = {Annales de l'Institut Henri Poincaré, Probabilités et Statistiques},
number = {3},
publisher = {Institut Henri Poincaré},
pages = {1281 -- 1305},
year = {2021},
doi = {10.1214/20-AIHP1117},
URL = {https://doi.org/10.1214/20-AIHP1117}
}

@article{widom2002convergence,
  title={On convergence of moments for random {Y}oung tableaux and a random growth model},
  author={Widom, Harold},
  journal={International Mathematics Research Notices},
  volume={2002},
  number={9},
  pages={455--464},
  year={2002},
  publisher={OUP}
}

@article{BCGStochasticSixVertex,
author = {Alexei Borodin and Ivan Corwin and Vadim Gorin},
title = {{Stochastic six-vertex model}},
volume = {165},
journal = {Duke Mathematical Journal},
number = {3},
publisher = {Duke University Press},
pages = {563 -- 624},
year = {2016},
doi = {10.1215/00127094-3166843},
URL = {https://doi.org/10.1215/00127094-3166843}
}

@article{liggett1976coupling,
  title={Coupling the simple exclusion process},
  author={Liggett, Thomas M},
  journal={The Annals of Probability},
  pages={339--356},
  year={1976},
  publisher={JSTOR}
}

@article{MartinStationaryASEP,
author = {James B. Martin},
title = {{Stationary distributions of the multi-type ASEP}},
volume = {25},
journal = {Electronic Journal of Probability},
number = {none},
publisher = {Institute of Mathematical Statistics and Bernoulli Society},
pages = {1 -- 41},
year = {2020},
doi = {10.1214/20-EJP421},
URL = {https://doi.org/10.1214/20-EJP421}
}

@article{aggarwal2024scaling,
  title={Scaling limit of the colored {ASEP} and stochastic six-vertex models},
  author={Aggarwal, Amol and Corwin, Ivan and Hegde, Milind},
  journal={arXiv preprint arXiv:2403.01341},
  year={2024}
}

@article {drillickLinS6V,
    AUTHOR = {Drillick, Hindy and Lin, Yier},
     TITLE = {Strong law of large numbers for the stochastic six vertex
              model},
   JOURNAL = {Electron. J. Probab.},
  FJOURNAL = {Electronic Journal of Probability},
    VOLUME = {28},
      YEAR = {2023},
     PAGES = {Paper No. 148, 21},
      ISSN = {1083-6489},
   MRCLASS = {60K35 (60F15 82C22)},
  MRNUMBER = {4669747},
       DOI = {10.1214/23-ejp1041},
       URL = {https://doi.org/10.1214/23-ejp1041},
}

@article {PeresWinkler2013CensoringInequality,
    AUTHOR = {Peres, Yuval and Winkler, Peter},
     TITLE = {Can extra updates delay mixing?},
   JOURNAL = {Comm. Math. Phys.},
  FJOURNAL = {Communications in Mathematical Physics},
    VOLUME = {323},
      YEAR = {2013},
    NUMBER = {3},
     PAGES = {1007--1016},
      ISSN = {0010-3616,1432-0916},
   MRCLASS = {82C20 (60K35)},
  MRNUMBER = {3106501},
MRREVIEWER = {Rapha\"{e}l\ Cerf},
       DOI = {10.1007/s00220-013-1776-0},
       URL = {https://doi.org/10.1007/s00220-013-1776-0},
}

@article {REzakhanlou1995Shocks,
    AUTHOR = {Rezakhanlou, Fraydoun},
     TITLE = {Microscopic structure of shocks in one conservation laws},
   JOURNAL = {Ann. Inst. H. Poincar\'{e} C Anal. Non Lin\'{e}aire},
  FJOURNAL = {Annales de l'Institut Henri Poincar\'{e} C. Analyse Non
              Lin\'{e}aire},
    VOLUME = {12},
      YEAR = {1995},
    NUMBER = {2},
     PAGES = {119--153},
      ISSN = {0294-1449,1873-1430},
   MRCLASS = {60K35 (35L67 35Q99 82C21)},
  MRNUMBER = {1326665},
MRREVIEWER = {Pablo\ A.\ Ferrari},
       DOI = {10.1016/S0294-1449(16)30161-5},
       URL = {https://doi.org/10.1016/S0294-1449(16)30161-5},
}

@article{pauling1935structure,
  title={The structure and entropy of ice and of other crystals with some randomness of atomic arrangement},
  author={Pauling, Linus},
  journal={Journal of the American Chemical Society},
  volume={57},
  number={12},
  pages={2680--2684},
  year={1935},
  publisher={ACS Publications}
}

@article {Aggarwal2017S6VtoASEP,
    AUTHOR = {Aggarwal, Amol},
     TITLE = {Convergence of the stochastic six-vertex model to the {ASEP}:
              stochastic six-vertex model and {ASEP}},
   JOURNAL = {Math. Phys. Anal. Geom.},
  FJOURNAL = {Mathematical Physics, Analysis and Geometry. An International
              Journal Devoted to the Theory and Applications of Analysis and
              Geometry to Physics},
    VOLUME = {20},
      YEAR = {2017},
    NUMBER = {2},
     PAGES = {Paper No. 3, 20},
      ISSN = {1385-0172,1572-9656},
   MRCLASS = {82C22 (60K35)},
  MRNUMBER = {3589552},
MRREVIEWER = {B.\ L.\ Granovsky},
       DOI = {10.1007/s11040-016-9235-8},
       URL = {https://doi.org/10.1007/s11040-016-9235-8},
}

@article {BorOlsh2006Meixner,
    AUTHOR = {Borodin, Alexei and Olshanski, Grigori},
     TITLE = {Meixner polynomials and random partitions},
   JOURNAL = {Mosc. Math. J.},
  FJOURNAL = {Moscow Mathematical Journal},
    VOLUME = {6},
      YEAR = {2006},
    NUMBER = {4},
     PAGES = {629--655, 771},
      ISSN = {1609-3321,1609-4514},
   MRCLASS = {60C05 (05E10 33C45)},
  MRNUMBER = {2291156},
MRREVIEWER = {Martin\ V.\ Hildebrand},
       DOI = {10.17323/1609-4514-2006-6-4-629-655},
       URL = {https://doi.org/10.17323/1609-4514-2006-6-4-629-655},
}

@incollection {BaikFerrariPeche2014twopointTASEP,
    AUTHOR = {Baik, Jinho and Ferrari, Patrik Lino and P\'{e}ch\'{e},
              Sandrine},
     TITLE = {Convergence of the two-point function of the stationary
              {TASEP}},
 BOOKTITLE = {Singular phenomena and scaling in mathematical models},
     PAGES = {91--110},
 PUBLISHER = {Springer, Cham},
      YEAR = {2014},
      ISBN = {978-3-319-00785-4; 978-3-319-00786-1},
   MRCLASS = {60K35},
  MRNUMBER = {3205038},
MRREVIEWER = {Benedikt\ Jahnel},
       DOI = {10.1007/978-3-319-00786-1\_5},
       URL = {https://doi.org/10.1007/978-3-319-00786-1_5},
}

@article{lowe2001moderateupper,
  title={Moderate deviations for longest increasing subsequences: the upper tail},
  author={L{\"o}we, Matthias and Merkl, Franz},
  journal={Communications on Pure and Applied Mathematics: A Journal Issued by the Courant Institute of Mathematical Sciences},
  volume={54},
  number={12},
  pages={1488--1519},
  year={2001},
  publisher={Wiley Online Library}
}

@article{lowe2002moderatelower,
  title={Moderate deviations for longest increasing subsequences: the lower tail},
  author={L{\"o}we, Matthias and Merkl, Franz and Rolles, Silke},
  journal={Journal of Theoretical Probability},
  volume={15},
  pages={1031--1047},
  year={2002},
  publisher={Springer}
}

@article {Yier2020KPZlimit,
    AUTHOR = {Lin, Yier},
     TITLE = {K{PZ} equation limit of stochastic higher spin six vertex
              model},
   JOURNAL = {Math. Phys. Anal. Geom.},
  FJOURNAL = {Mathematical Physics, Analysis and Geometry. An International
              Journal Devoted to the Theory and Applications of Analysis and
              Geometry to Physics},
    VOLUME = {23},
      YEAR = {2020},
    NUMBER = {1},
     PAGES = {Paper No. 1, 118},
      ISSN = {1385-0172,1572-9656},
   MRCLASS = {82C22 (60K35 82B23)},
  MRNUMBER = {4046044},
       DOI = {10.1007/s11040-019-9325-5},
       URL = {https://doi.org/10.1007/s11040-019-9325-5},
}

@article {CorwinTsai2017KPZlimit,
    AUTHOR = {Corwin, Ivan and Tsai, Li-Cheng},
     TITLE = {K{PZ} equation limit of higher-spin exclusion processes},
   JOURNAL = {Ann. Probab.},
  FJOURNAL = {The Annals of Probability},
    VOLUME = {45},
      YEAR = {2017},
    NUMBER = {3},
     PAGES = {1771--1798},
      ISSN = {0091-1798,2168-894X},
   MRCLASS = {60K35 (82C22 82C23)},
  MRNUMBER = {3650415},
MRREVIEWER = {Patr\'{\i}cia\ Gon\c{c}alves},
       DOI = {10.1214/16-AOP1101},
       URL = {https://doi.org/10.1214/16-AOP1101},
}

@article {CorwinGhosalShenTsai2020SPDElimit,
    AUTHOR = {Corwin, Ivan and Ghosal, Promit and Shen, Hao and Tsai,
              Li-Cheng},
     TITLE = {Stochastic {PDE} limit of the six vertex model},
   JOURNAL = {Comm. Math. Phys.},
  FJOURNAL = {Communications in Mathematical Physics},
    VOLUME = {375},
      YEAR = {2020},
    NUMBER = {3},
     PAGES = {1945--2038},
      ISSN = {0010-3616,1432-0916},
   MRCLASS = {60K35 (35R60 60H15 82C22)},
  MRNUMBER = {4091508},
       DOI = {10.1007/s00220-019-03678-z},
       URL = {https://doi.org/10.1007/s00220-019-03678-z},
}

@article {BorodinGorin2019StochasticTelegraph,
    AUTHOR = {Borodin, Alexei and Gorin, Vadim},
     TITLE = {A stochastic telegraph equation from the six-vertex model},
   JOURNAL = {Ann. Probab.},
  FJOURNAL = {The Annals of Probability},
    VOLUME = {47},
      YEAR = {2019},
    NUMBER = {6},
     PAGES = {4137--4194},
      ISSN = {0091-1798,2168-894X},
   MRCLASS = {60G60 (35R60 60H15 82B20)},
  MRNUMBER = {4038051},
       DOI = {10.1214/19-aop1356},
       URL = {https://doi.org/10.1214/19-aop1356},
}

@article {BorodinPetrov2018higherspinSymmetric,
    AUTHOR = {Borodin, Alexei and Petrov, Leonid},
     TITLE = {Higher spin six vertex model and symmetric rational functions},
   JOURNAL = {Selecta Math. (N.S.)},
  FJOURNAL = {Selecta Mathematica. New Series},
    VOLUME = {24},
      YEAR = {2018},
    NUMBER = {2},
     PAGES = {751--874},
      ISSN = {1022-1824,1420-9020},
   MRCLASS = {60K35 (05E05 82B23)},
  MRNUMBER = {3782413},
       DOI = {10.1007/s00029-016-0301-7},
       URL = {https://doi.org/10.1007/s00029-016-0301-7},
}

@article {CorwinPetrov2016higherspinLine,
    AUTHOR = {Corwin, Ivan and Petrov, Leonid},
     TITLE = {Stochastic higher spin vertex models on the line},
   JOURNAL = {Comm. Math. Phys.},
  FJOURNAL = {Communications in Mathematical Physics},
    VOLUME = {343},
      YEAR = {2016},
    NUMBER = {2},
     PAGES = {651--700},
      ISSN = {0010-3616,1432-0916},
   MRCLASS = {82B23 (60K35 82C23)},
  MRNUMBER = {3477349},
MRREVIEWER = {Aernout\ C. D. van Enter},
       DOI = {10.1007/s00220-015-2479-5},
       URL = {https://doi.org/10.1007/s00220-015-2479-5},
}

@article {CorwinDimitrov2018Transversalfluctuations,
    AUTHOR = {Corwin, Ivan and Dimitrov, Evgeni},
     TITLE = {Transversal fluctuations of the {ASEP}, stochastic six vertex
              model, and {H}all-{L}ittlewood {G}ibbsian line ensembles},
   JOURNAL = {Comm. Math. Phys.},
  FJOURNAL = {Communications in Mathematical Physics},
    VOLUME = {363},
      YEAR = {2018},
    NUMBER = {2},
     PAGES = {435--501},
      ISSN = {0010-3616,1432-0916},
   MRCLASS = {60K35 (82C22)},
  MRNUMBER = {3851820},
MRREVIEWER = {Sunder\ Sethuraman},
       DOI = {10.1007/s00220-018-3139-3},
       URL = {https://doi.org/10.1007/s00220-018-3139-3},
}

@article{aggarwal2022nonexistence,
  title={Nonexistence and uniqueness for pure states of ferroelectric six-vertex models},
  author={Aggarwal, Amol},
  journal={Proceedings of the London Mathematical Society},
  volume={124},
  number={3},
  pages={387--425},
  year={2022},
  publisher={Wiley Online Library}
}

@article{landonSosoe2023tail,
  title={Tail estimates for the stationary stochastic six vertex model and {ASEP}},
  author={Landon, Benjamin and Sosoe, Philippe},
  journal={arXiv preprint arXiv:2308.16812},
  year={2023}
}

@article{dasLiaoMucciconi2024lower,
  title={Lower tail large deviations of the stochastic six vertex model},
  author={Das, Sayan and Liao, Yuchen and Mucciconi, Matteo},
  journal={arXiv preprint arXiv:2407.08530},
  year={2024}
}

@article{busani2022scaling,
  title={Scaling limit of the {TASEP} speed process},
  author={Busani, Ofer and Sepp{\"a}l{\"a}inen, Timo and Sorensen, Evan},
  journal={arXiv preprint arXiv:2211.04651},
  year={2022}
}

@article {MountfordGuiol05,
    AUTHOR = {Mountford, Thomas and Guiol, Herv\'e},
     TITLE = {The motion of a second class particle for the {TASEP} starting
              from a decreasing shock profile},
   JOURNAL = {Ann. Appl. Probab.},
  FJOURNAL = {The Annals of Applied Probability},
    VOLUME = {15},
      YEAR = {2005},
    NUMBER = {2},
     PAGES = {1227--1259},
      ISSN = {1050-5164,2168-8737},
   MRCLASS = {60K35 (60F15 82C22)},
  MRNUMBER = {2134103},
MRREVIEWER = {Zhong\ Gen\ Su},
       DOI = {10.1214/105051605000000151},
       URL = {https://doi.org/10.1214/105051605000000151},
}

@article {FerrariPimentel05,
    AUTHOR = {Ferrari, Pablo A. and Pimentel, Leandro P. R.},
     TITLE = {Competition interfaces and second class particles},
   JOURNAL = {Ann. Probab.},
  FJOURNAL = {The Annals of Probability},
    VOLUME = {33},
      YEAR = {2005},
    NUMBER = {4},
     PAGES = {1235--1254},
      ISSN = {0091-1798,2168-894X},
   MRCLASS = {60K35 (82B21 82B24)},
  MRNUMBER = {2150188},
MRREVIEWER = {Nicoletta\ Cancrini},
       DOI = {10.1214/009117905000000080},
       URL = {https://doi.org/10.1214/009117905000000080},
}

@article {FerrariMartinPimentel09,
    AUTHOR = {Ferrari, Pablo A. and Martin, James B. and Pimentel, Leandro
              P. R.},
     TITLE = {A phase transition for competition interfaces},
   JOURNAL = {Ann. Appl. Probab.},
  FJOURNAL = {The Annals of Applied Probability},
    VOLUME = {19},
      YEAR = {2009},
    NUMBER = {1},
     PAGES = {281--317},
      ISSN = {1050-5164,2168-8737},
   MRCLASS = {60K35 (82B26)},
  MRNUMBER = {2498679},
       DOI = {10.1214/08-AAP542},
       URL = {https://doi.org/10.1214/08-AAP542},
}

@article {Goncalves14,
    AUTHOR = {Gon\c{c}alves, Patr\'icia},
     TITLE = {On the asymmetric zero-range in the rarefaction fan},
   JOURNAL = {J. Stat. Phys.},
  FJOURNAL = {Journal of Statistical Physics},
    VOLUME = {154},
      YEAR = {2014},
    NUMBER = {4},
     PAGES = {1074--1095},
      ISSN = {0022-4715,1572-9613},
   MRCLASS = {82C22},
  MRNUMBER = {3164603},
       DOI = {10.1007/s10955-013-0892-8},
       URL = {https://doi.org/10.1007/s10955-013-0892-8},
}

@article {busani2024scalinglimitmultitypeinvariant,
    AUTHOR = {Busani, Ofer and Sepp\"al\"ainen, Timo and Sorensen, Evan},
     TITLE = {Scaling {L}imit of {M}ulti-{T}ype {I}nvariant {M}easures via
              the {D}irected {L}andscape},
   JOURNAL = {Int. Math. Res. Not. IMRN},
  FJOURNAL = {International Mathematics Research Notices. IMRN},
      YEAR = {2024},
    NUMBER = {17},
     PAGES = {12382--12432},
      ISSN = {1073-7928,1687-0247},
   MRCLASS = {60K35 (28A33)},
  MRNUMBER = {4795007},
       DOI = {10.1093/imrn/rnae168},
       URL = {https://doi.org/10.1093/imrn/rnae168},
}

@misc{ghosal2019limitingspeedsecondclass,
      title={Limiting speed of a second class particle in {ASEP}}, 
      author={Promit Ghosal and Axel Saenz and Ethan C. Zell},
      year={2019},
      eprint={1903.09615},
      archivePrefix={arXiv},
      primaryClass={math.PR},
      url={https://arxiv.org/abs/1903.09615}, 
}

@article {MR4010933,
    AUTHOR = {Ferrari, P. L. and Ghosal, P. and Nejjar, P.},
     TITLE = {Limit law of a second class particle in {TASEP} with
              non-random initial condition},
   JOURNAL = {Ann. Inst. Henri Poincar\'e{} Probab. Stat.},
  FJOURNAL = {Annales de l'Institut Henri Poincar\'e{} Probabilit\'es et
              Statistiques},
    VOLUME = {55},
      YEAR = {2019},
    NUMBER = {3},
     PAGES = {1203--1225},
      ISSN = {0246-0203,1778-7017},
   MRCLASS = {82C22 (60B20 60K35)},
  MRNUMBER = {4010933},
       DOI = {10.1214/18-aihp916},
       URL = {https://doi.org/10.1214/18-aihp916},
}

@article {MR3112921,
    AUTHOR = {Cator, Eric and Pimentel, Leandro P. R.},
     TITLE = {Busemann functions and the speed of a second class particle in
              the rarefaction fan},
   JOURNAL = {Ann. Probab.},
  FJOURNAL = {The Annals of Probability},
    VOLUME = {41},
      YEAR = {2013},
    NUMBER = {4},
     PAGES = {2401--2425},
      ISSN = {0091-1798,2168-894X},
   MRCLASS = {60K35 (82C21)},
  MRNUMBER = {3112921},
MRREVIEWER = {Tertuliano\ Franco},
       DOI = {10.1214/11-AOP709},
       URL = {https://doi.org/10.1214/11-AOP709},
}

@article {MR2312944,
    AUTHOR = {Coletti, Cristian F. and Pimentel, Leandro P. R.},
     TITLE = {On the collision between two {PNG} droplets},
   JOURNAL = {J. Stat. Phys.},
  FJOURNAL = {Journal of Statistical Physics},
    VOLUME = {126},
      YEAR = {2007},
    NUMBER = {6},
     PAGES = {1145--1164},
      ISSN = {0022-4715,1572-9613},
   MRCLASS = {82D60 (60K35)},
  MRNUMBER = {2312944},
MRREVIEWER = {Nicoletta\ Cancrini},
       DOI = {10.1007/s10955-006-9272-y},
       URL = {https://doi.org/10.1007/s10955-006-9272-y},
}

@article {Harris1978Graphical,
    AUTHOR = {Harris, T. E.},
     TITLE = {Additive set-valued {M}arkov processes and graphical methods},
   JOURNAL = {Ann. Probability},
  FJOURNAL = {The Annals of Probability},
    VOLUME = {6},
      YEAR = {1978},
    NUMBER = {3},
     PAGES = {355--378},
      ISSN = {0091-1798},
   MRCLASS = {60K35},
  MRNUMBER = {488377},
MRREVIEWER = {Richard\ Holley},
       DOI = {10.1214/aop/1176995523},
       URL = {https://doi.org/10.1214/aop/1176995523},
}

@article {busani2023diffusivescalinglimitbusemann,
    AUTHOR = {Busani, Ofer},
     TITLE = {Diffusive scaling limit of the {B}usemann process in last
              passage percolation},
   JOURNAL = {Ann. Probab.},
  FJOURNAL = {The Annals of Probability},
    VOLUME = {52},
      YEAR = {2024},
    NUMBER = {5},
     PAGES = {1650--1712},
      ISSN = {0091-1798,2168-894X},
   MRCLASS = {60K35 (60K37)},
  MRNUMBER = {4791418},
       DOI = {10.1214/24-aop1681},
       URL = {https://doi.org/10.1214/24-aop1681},
}

@article {MR2630064,
    AUTHOR = {Bal\'azs, M\'arton and Sepp\"al\"ainen, Timo},
     TITLE = {Order of current variance and diffusivity in the asymmetric
              simple exclusion process},
   JOURNAL = {Ann. of Math. (2)},
  FJOURNAL = {Annals of Mathematics. Second Series},
    VOLUME = {171},
      YEAR = {2010},
    NUMBER = {2},
     PAGES = {1237--1265},
      ISSN = {0003-486X,1939-8980},
   MRCLASS = {60K35 (60J25)},
  MRNUMBER = {2630064},
MRREVIEWER = {Daniel\ Remenik},
       DOI = {10.4007/annals.2010.171.1237},
       URL = {https://doi.org/10.4007/annals.2010.171.1237},
}

@article {MR2318311,
    AUTHOR = {Quastel, Jeremy and Valko, Benedek},
     TITLE = {{$t^{1/3}$} {S}uperdiffusivity of finite-range asymmetric
              exclusion processes on {$\mathbb Z$}},
   JOURNAL = {Comm. Math. Phys.},
  FJOURNAL = {Communications in Mathematical Physics},
    VOLUME = {273},
      YEAR = {2007},
    NUMBER = {2},
     PAGES = {379--394},
      ISSN = {0010-3616,1432-0916},
   MRCLASS = {60K35 (60J25 82B41 82C22)},
  MRNUMBER = {2318311},
MRREVIEWER = {C\'edric\ Bernardin},
       DOI = {10.1007/s00220-007-0242-2},
       URL = {https://doi.org/10.1007/s00220-007-0242-2},
}

@article {MR4709103,
    AUTHOR = {Corwin, Ivan and Hegde, Milind},
     TITLE = {The lower tail of {$q$}-push{TASEP}},
   JOURNAL = {Comm. Math. Phys.},
  FJOURNAL = {Communications in Mathematical Physics},
    VOLUME = {405},
      YEAR = {2024},
    NUMBER = {3},
     PAGES = {Paper No. 64, 55},
      ISSN = {0010-3616,1432-0916},
   MRCLASS = {60K35 (82C22)},
  MRNUMBER = {4709103},
       DOI = {10.1007/s00220-024-04944-5},
       URL = {https://doi.org/10.1007/s00220-024-04944-5},
}

@unpublished{GhosalSilva2024,
  author = "Ghosal, Promit and Silva, Guilherme L. F.",
  year = {2024},
  title  = "Six-vertex model and the {M}eixner Ensemble",
note="In preparation",
}

@article{ferrari1991microscopic,
  title={Microscopic structure of travelling waves in the asymmetric simple exclusion process},
  author={Ferrari, Pablo A and Kipnis, Claude and Saada, Ellen},
  journal={The Annals of Probability},
  volume={19},
  number={1},
  pages={226--244},
  year={1991},
  publisher={Institute of Mathematical Statistics}
}

@article {MR4193889,
    AUTHOR = {Lin, Yier},
     TITLE = {The stochastic telegraph equation limit of the stochastic
              higher spin six vertex model},
   JOURNAL = {Electron. J. Probab.},
  FJOURNAL = {Electronic Journal of Probability},
    VOLUME = {25},
      YEAR = {2020},
     PAGES = {Paper No. 148, 30},
      ISSN = {1083-6489},
   MRCLASS = {60H15 (60F17 82B20)},
  MRNUMBER = {4193889},
       DOI = {10.1214/20-ejp552},
       URL = {https://doi.org/10.1214/20-ejp552},
}

@article {MR3951443,
    AUTHOR = {Shen, Hao and Tsai, Li-Cheng},
     TITLE = {Stochastic telegraph equation limit for the stochastic six
              vertex model},
   JOURNAL = {Proc. Amer. Math. Soc.},
  FJOURNAL = {Proceedings of the American Mathematical Society},
    VOLUME = {147},
      YEAR = {2019},
    NUMBER = {6},
     PAGES = {2685--2705},
      ISSN = {0002-9939,1088-6826},
   MRCLASS = {60H15 (82B20)},
  MRNUMBER = {3951443},
       DOI = {10.1090/proc/14415},
       URL = {https://doi.org/10.1090/proc/14415},
}

@article {MR4561796,
    AUTHOR = {Dimitrov, Evgeni},
     TITLE = {Two-point convergence of the stochastic six-vertex model to
              the {A}iry process},
   JOURNAL = {Comm. Math. Phys.},
  FJOURNAL = {Communications in Mathematical Physics},
    VOLUME = {398},
      YEAR = {2023},
    NUMBER = {3},
     PAGES = {925--1027},
      ISSN = {0010-3616,1432-0916},
   MRCLASS = {82B20 (60K35 82B44)},
  MRNUMBER = {4561796},
MRREVIEWER = {Leonid\ Petrov},
       DOI = {10.1007/s00220-022-04499-3},
       URL = {https://doi.org/10.1007/s00220-022-04499-3},
}

@article {MR2930377,
    AUTHOR = {Corwin, Ivan},
     TITLE = {The {K}ardar-{P}arisi-{Z}hang equation and universality class},
   JOURNAL = {Random Matrices Theory Appl.},
  FJOURNAL = {Random Matrices. Theory and Applications},
    VOLUME = {1},
      YEAR = {2012},
    NUMBER = {1},
     PAGES = {1130001, 76},
      ISSN = {2010-3263,2010-3271},
   MRCLASS = {82B31 (60B20 60K35 60K37)},
  MRNUMBER = {2930377},
       DOI = {10.1142/S2010326311300014},
       URL = {https://doi.org/10.1142/S2010326311300014},
}

@article{basu2014last,
  title={Last passage percolation with a defect line and the solution of the slow bond problem},
  author={Basu, Riddhipratim and Sidoravicius, Vladas and Sly, Allan},
  journal={arXiv preprint arXiv:1408.3464},
  year={2014}
}

@article {MR3304747,
    AUTHOR = {Ferrari, Patrik L. and Nejjar, Peter},
     TITLE = {Anomalous shock fluctuations in {TASEP} and last passage
              percolation models},
   JOURNAL = {Probab. Theory Related Fields},
  FJOURNAL = {Probability Theory and Related Fields},
    VOLUME = {161},
      YEAR = {2015},
    NUMBER = {1-2},
     PAGES = {61--109},
      ISSN = {0178-8051,1432-2064},
   MRCLASS = {60K35 (82C22)},
  MRNUMBER = {3304747},
MRREVIEWER = {Adrian\ Muntean},
       DOI = {10.1007/s00440-013-0544-6},
       URL = {https://doi.org/10.1007/s00440-013-0544-6},
}

@article {MR4620410,
    AUTHOR = {Emrah, Elnur and Janjigian, Christopher and Sepp\"al\"ainen,
              Timo},
     TITLE = {Optimal-order exit point bounds in exponential last-passage
              percolation via the coupling technique},
   JOURNAL = {Probab. Math. Phys.},
  FJOURNAL = {Probability and Mathematical Physics},
    VOLUME = {4},
      YEAR = {2023},
    NUMBER = {3},
     PAGES = {609--666},
      ISSN = {2690-0998,2690-1005},
   MRCLASS = {60K35 (60K37 82B43)},
  MRNUMBER = {4620410},
MRREVIEWER = {Andrew\ R.\ Wade},
       DOI = {10.2140/pmp.2023.4.609},
       URL = {https://doi.org/10.2140/pmp.2023.4.609},
}

@article {MR4610276,
    AUTHOR = {Landon, Benjamin and Sosoe, Philippe},
     TITLE = {Upper tail bounds for stationary {KPZ} models},
   JOURNAL = {Comm. Math. Phys.},
  FJOURNAL = {Communications in Mathematical Physics},
    VOLUME = {401},
      YEAR = {2023},
    NUMBER = {2},
     PAGES = {1311--1335},
      ISSN = {0010-3616,1432-0916},
   MRCLASS = {82D60 (60K35)},
  MRNUMBER = {4610276},
MRREVIEWER = {G\"okhan\ Y\i ld\i r\i m},
       DOI = {10.1007/s00220-023-04669-x},
       URL = {https://doi.org/10.1007/s00220-023-04669-x},
}

@misc{emrah2020righttailmoderatedeviationsexponential,
      title={Right-tail moderate deviations in the exponential last-passage percolation}, 
      author={Elnur Emrah and Chris Janjigian and Timo Seppäläinen},
      year={2020},
      eprint={2004.04285},
      archivePrefix={arXiv},
      primaryClass={math.PR},
      url={https://arxiv.org/abs/2004.04285}, 
}

@article{Kulish:1981aa,
	author = {Kulish, P. P. and Reshetikhin, N. Yu. and Sklyanin, E. K.},
	date = {1981/09/01},
	doi = {10.1007/BF02285311},
	id = {Kulish1981},
	isbn = {1573-0530},
	journal = {Letters in Mathematical Physics},
	number = {5},
	pages = {393--403},
	title = {Yang-Baxter equation and representation theory: I},
	url = {https://doi.org/10.1007/BF02285311},
	volume = {5},
	year = {1981}}

@article {MR806529,
    AUTHOR = {Bazhanov, V. V.},
     TITLE = {Trigonometric solutions of triangle equations and classical
              {L}ie algebras},
   JOURNAL = {Phys. Lett. B},
  FJOURNAL = {Physics Letters. B. Particle Physics, Nuclear Physics and
              Cosmology},
    VOLUME = {159},
      YEAR = {1985},
    NUMBER = {4-6},
     PAGES = {321--324},
      ISSN = {0370-2693,1873-2445},
   MRCLASS = {82A15 (17B10 81E15 82A68)},
  MRNUMBER = {806529},
MRREVIEWER = {F.\ W.\ Nijhoff},
       DOI = {10.1016/0370-2693(85)90259-X},
       URL = {https://doi.org/10.1016/0370-2693(85)90259-X},
}

@article {MR824090,
    AUTHOR = {Jimbo, Michio},
     TITLE = {Quantum {$R$} matrix for the generalized {T}oda system},
   JOURNAL = {Comm. Math. Phys.},
  FJOURNAL = {Communications in Mathematical Physics},
    VOLUME = {102},
      YEAR = {1986},
    NUMBER = {4},
     PAGES = {537--547},
      ISSN = {0010-3616,1432-0916},
   MRCLASS = {58F07 (82A15)},
  MRNUMBER = {824090},
MRREVIEWER = {L.\ A.\ Takhtajan (Takhtadzhyan)},
       URL = {http://projecteuclid.org/euclid.cmp/1104114539},
}

@article{kuniba2016stochastic,
	author = {Kuniba, Atsuo and Mangazeev, Vladimir V and Maruyama, Shouya and Okado, Masato},
	journal = {Nuclear Physics B},
	pages = {248--277},
	publisher = {Elsevier},
	title = {Stochastic {$R$} matrix for {$U_q (A_n^{(1)})$}},
	volume = {913},
	year = {2016}}

@article {MR4518478,
    AUTHOR = {Borodin, Alexei and Wheeler, Michael},
     TITLE = {Colored stochastic vertex models and their spectral theory},
   JOURNAL = {Ast\'erisque},
  FJOURNAL = {Ast\'erisque},
    NUMBER = {437},
      YEAR = {2022},
     PAGES = {ix+225},
      ISSN = {0303-1179,2492-5926},
      ISBN = {978-2-85629-963-0},
   MRCLASS = {05E05 (05E10 60K35 82B23)},
  MRNUMBER = {4518478},
MRREVIEWER = {Leonid\ Petrov},
       DOI = {10.24033/ast.1180},
       URL = {https://doi.org/10.24033/ast.1180},
}

@misc{aggarwal2024kpzfixedpointconvergence,
      title={{KPZ} fixed point convergence of the {ASEP} and stochastic six-vertex models}, 
      author={Amol Aggarwal and Ivan Corwin and Milind Hegde},
      year={2024},
      eprint={2412.18117},
      archivePrefix={arXiv},
      primaryClass={math.PR},
      url={https://arxiv.org/abs/2412.18117}, 
}

@article {MR4094738,
    AUTHOR = {Corwin, Ivan and Ghosal, Promit},
     TITLE = {Lower tail of the {KPZ} equation},
   JOURNAL = {Duke Math. J.},
  FJOURNAL = {Duke Mathematical Journal},
    VOLUME = {169},
      YEAR = {2020},
    NUMBER = {7},
     PAGES = {1329--1395},
      ISSN = {0012-7094,1547-7398},
   MRCLASS = {35R60 (60B20 60H25 82C22)},
  MRNUMBER = {4094738},
MRREVIEWER = {Le\ Chen},
       DOI = {10.1215/00127094-2019-0079},
       URL = {https://doi.org/10.1215/00127094-2019-0079},
}

@article {MR4373176,
    AUTHOR = {Cafasso, Mattia and Claeys, Tom},
     TITLE = {A {R}iemann-{H}ilbert approach to the lower tail of the
              {K}ardar-{P}arisi-{Z}hang equation},
   JOURNAL = {Comm. Pure Appl. Math.},
  FJOURNAL = {Communications on Pure and Applied Mathematics},
    VOLUME = {75},
      YEAR = {2022},
    NUMBER = {3},
     PAGES = {493--540},
      ISSN = {0010-3640,1097-0312},
   MRCLASS = {60B20 (60F10 60G55)},
  MRNUMBER = {4373176},
       DOI = {10.1002/cpa.21978},
       URL = {https://doi.org/10.1002/cpa.21978},
}

@article {MR4318370,
    AUTHOR = {Kim, Yujin H.},
     TITLE = {The lower tail of the half-space {KPZ} equation},
   JOURNAL = {Stochastic Process. Appl.},
  FJOURNAL = {Stochastic Processes and their Applications},
    VOLUME = {142},
      YEAR = {2021},
     PAGES = {365--406},
      ISSN = {0304-4149,1879-209X},
   MRCLASS = {60H15 (45M05 60B20 60F10 60G55 60H25)},
  MRNUMBER = {4318370},
       DOI = {10.1016/j.spa.2021.09.001},
       URL = {https://doi.org/10.1016/j.spa.2021.09.001},
}

@article {MR4294287,
    AUTHOR = {Cafasso, Mattia and Claeys, Tom and Ruzza, Giulio},
     TITLE = {Airy kernel determinant solutions to the {K}d{V} equation and
              integro-differential {P}ainlev\'e{} equations},
   JOURNAL = {Comm. Math. Phys.},
  FJOURNAL = {Communications in Mathematical Physics},
    VOLUME = {386},
      YEAR = {2021},
    NUMBER = {2},
     PAGES = {1107--1153},
      ISSN = {0010-3616,1432-0916},
   MRCLASS = {35Q53 (33C10 34M50 34M55 45K05)},
  MRNUMBER = {4294287},
       DOI = {10.1007/s00220-021-04108-9},
       URL = {https://doi.org/10.1007/s00220-021-04108-9},
}

@misc{shen2024timecorrelationsinversegammapolymer,
      title={Time correlations in the inverse-gamma polymer with flat initial condition}, 
      author={Xiao Shen},
      year={2024},
      eprint={2408.07928},
      archivePrefix={arXiv},
      primaryClass={math.PR},
      url={https://arxiv.org/abs/2408.07928}, 
}

@article {MR4771973,
    AUTHOR = {Landon, Benjamin and Sosoe, Philippe},
     TITLE = {Tail bounds for the {O}'{C}onnell-{Y}or polymer},
   JOURNAL = {Electron. J. Probab.},
  FJOURNAL = {Electronic Journal of Probability},
    VOLUME = {29},
      YEAR = {2024},
     PAGES = {Paper No. 99, 47},
      ISSN = {1083-6489},
   MRCLASS = {82D60},
  MRNUMBER = {4771973},
       DOI = {10.1214/24-ejp1162},
       URL = {https://doi.org/10.1214/24-ejp1162},
}

@article {MR3306003,
    AUTHOR = {Romik, Dan and \'Sniady, Piotr},
     TITLE = {Jeu de taquin dynamics on infinite {Y}oung tableaux and second
              class particles},
   JOURNAL = {Ann. Probab.},
  FJOURNAL = {The Annals of Probability},
    VOLUME = {43},
      YEAR = {2015},
    NUMBER = {2},
     PAGES = {682--737},
      ISSN = {0091-1798,2168-894X},
   MRCLASS = {60C05 (05E10 37A05 60K35 82C22)},
  MRNUMBER = {3306003},
MRREVIEWER = {Boyka\ L.\ Aneva},
       DOI = {10.1214/13-AOP873},
       URL = {https://doi.org/10.1214/13-AOP873},
}

@incollection{okounkov2002symmetric,
  title={Symmetric functions and random partitions},
  author={Okounkov, Andrei},
  booktitle={Symmetric functions 2001: surveys of developments and perspectives},
  pages={223--252},
  year={2002},
  publisher={Springer}
}

@book {MR3468920,
    AUTHOR = {Baik, Jinho and Deift, Percy and Suidan, Toufic},
     TITLE = {Combinatorics and random matrix theory},
    SERIES = {Graduate Studies in Mathematics},
    VOLUME = {172},
 PUBLISHER = {American Mathematical Society, Providence, RI},
      YEAR = {2016},
     PAGES = {xi+461},
      ISBN = {978-0-8218-4841-8},
   MRCLASS = {60B20 (30E25 33E17 41A60 47B35 82C23)},
  MRNUMBER = {3468920},
MRREVIEWER = {Terence\ Tao},
       DOI = {10.1090/gsm/172},
       URL = {https://doi.org/10.1090/gsm/172},
}
 
\end{document}